\documentclass[11pt]{article}
\pdfoutput=1
\usepackage{epsfig,graphics,graphicx,amsmath,amssymb,amsfonts}
\usepackage{epstopdf}
\usepackage{color}
\usepackage{amsmath}
\usepackage{cases}
\usepackage[position=bottom]{subfig}
\usepackage{paralist}
\usepackage{enumitem}
\usepackage{graphicx}
\usepackage[percent]{overpic}
\usepackage{tikz}

\newsavebox{\largestimage}

\setlist{nolistsep}

\newcommand{\bm}[1]{\mbox{\boldmath{$#1$}}}

\def\u{{\bm u}}
\def\v{{\bm{v}}}
\def\x{{\bm{x}}}

\def\0{{\bm 0}}

\def\cl {\nonumber \\}
\def\el {\nonumber}

\DeclareMathOperator*{\argmin}{arg\,min}

\begin{document}

\title{A kinetic theory approach to model pedestrian dynamics in bounded domains with obstacles}

\author{Daewa Kim$^*$ and Annalisa Quaini$^*$ \\
\footnotesize{$^*$Department of Mathematics, University of Houston, 3551 Cullen Blvd, Houston TX 77204}\\
\footnotesize{daewakim@math.uh.edu; quaini@math.uh.edu}
}

\maketitle

\noindent{\bf Abstract}
We consider a kinetic theory approach to model the evacuation of a crowd from bounded domains. The interactions of a person with other pedestrians and the environment, which includes walls, exits, and obstacles, are modeled by using tools of game theory and are transferred to the crowd dynamics. The model allows to weight between two competing behaviors: the search for less congested areas and the tendency to follow the stream unconsciously in a panic situation. 
For the numerical approximation of the solution to our model, we apply an operator splitting scheme
which breaks the problem into two pure advection problems and a problem involving the interactions.
We compare our numerical results against the data reported in a recent empirical study
on evacuation from a room with two exits.  
For medium and medium-to-large groups of people we achieve 
good agreement between the computed average people
density and flow rate and the respective measured quantities.
Through a series of numerical tests we also show that our approach is 
capable of handling evacuation from a room with one or more exits
with variable size, with and
without obstacles, and can reproduce lane formation in bidirectional flow in a corridor.

\section{Introduction} \label{sec:Intro}

The complex dynamical behavior of pedestrian crowds has fascinated researchers from various scientific fields
since the early 1950's. 
Academic studies started with empirical observations and continued with the development of models in the field of applied physics
and mathematics. The simulation of pedestrian flow has attracted increasing research attention in
recent years since a reliable simulation model for pedestrian flow
may greatly benefit engineers in mass transportation
management, and designers in urban planning and architecture.

A very large variety of models have been developed over the years. 
The different mathematical models can be divided into three main categories depending
on the scale of observation \cite{Bellomo2011}. 
A first approach corresponds to the macroscopic description:
evolution equations are derived for mass density 
and linear momentum, which are regarded as macroscopic observables of pedestrian
flow. See, e.g., \cite{HUGHES2002507,5773492}. Such an approach is 
suitable for high density, large-scale systems, which are not the focus of our work. 

A second approach looks at the problem at the microscopic level.
Microscopic models can be further divided into models which are grid-based or grid-free. 
Cellular Automata \cite{Blue1999, Blue2000, BURSTEDDE2001507, Dijkstra, Li_2012} models belong to the first category. 
They describe pedestrian flow in space-time by assigning discrete states to a grid of space-cells. 
These cells can be occupied by a pedestrian or be empty. 
Thus, the movement of pedestrians in space is done by passing them from cell to cell (discrete space) in discrete time. 
Grid-free methods {can be based on second order models (forces-based), 
first order models (vision-based or speed-based) or
zeroth order models (rule-based or decision-based).
Force-based models}
use Newtonian mechanics to interpret pedestrian movement 
as the physical interaction between the people and the environment, i.e. the action of other people
and the environment on a given pedestrian is modeled with forces.
These models are one of the most popular modeling paradigms of continuous models
because they describe the movement of pedestrians qualitatively well.
See, e.g.,
\cite{Chraibi2011425, BS:BS3830360405, Helbing1995, 1367-2630-1-1-313, 6701214, Moussad2755, TurnerPenn, 6248013} and references therein.
Collective phenomena, like unidirectional or bidirectional flow in a corridor, 
lane formation, oscillations at bottlenecks, the faster-is-slower effect,
and emergency evacuation from buildings, are well reproduced.
Agent-based models allow for flexibility, extensibility, and capability to realize heterogeneity
in crowd dynamics.
{For examples of vision-based, speed-based, rule-based, and decision-based models
we refer to \cite{ANTONINI2006667,ASANO2010842,bandini2009, CHOORAMUN20121685, Chraibi2019, DAI20132202}
and references therein.}
Both force-based and agent-based models may introduce artifacts due to the force
representation of human behavior, leading to unrealistic backward movement
or oscillating trajectories.
These artifacts can be reduced by incorporating extra rules and/or elaborate calibrations,
at the price of an increased computational cost.

{ 
The scale of observation for the third approach is between the previous two.
Introduced in \cite{Bellomo2011383} and further developed in 
\cite{Agnelli2015,  Bellomo2013_new,Bellomo2017_book,Bellomo2015_new,Bellomo2016_new,Bellomo2019_new,Bellomo2013},
this approach derives a Boltzmann-type evolution equation for the statistical
distribution function of the position and velocity of the pedestrians,
in a framework close to that of the kinetic theory of gases.
See also \cite{Bellomo2012} for a literature review on this approach.
The model proposed in \cite{Bellomo2011383,Bellomo2013_new, Bellomo2013} is valid in
unbounded domains and with a homogeneous distribution of walking ability for the pedestrians, 
while the extension to bounded domains is presented in \cite{Agnelli2015}
and further explored in \cite{Bellomo2015_new,Bellomo2016_new,Bellomo2019_new}. 
In \cite{Bellomo2015_new}, more general features of behavioral-social dynamics are
taken into account. In \cite{Bellomo2016_new}, Monte Carlo simulations are introduced 
to study pedestrians behavior in complex scenarios. The methodology in \cite{Bellomo2016_new} is 
validated by comparing the computed fundamental density-velocity diagrams
with empirically observed ones and by checking that well known emerging properties
are reproduced.
A kinetic theory approach for modeling pedestrian dynamics in presence of social phenomena,
such as the propagation of stress conditions, is presented in \cite{Bellomo2019_new}.
Finally, we refer to \cite{Bellomo2017_book} for a thorough description of
how kinetic theory and evolutionary game theory can be used to understand the dynamics of living systems.
}

The scale of observation for the third approach is between the previous two. 
In a framework close to that of the kinetic theory of gases,
this approach derives a Boltzmann-type evolution equation for the statistical
distribution function of the position and velocity of the pedestrians.
The kinetic theory approach was introduced in \cite{Bellomo2011383}
and further developed in \cite{Bellomo2013}.
The model in \cite{Bellomo2011383,Bellomo2013} is valid in unbounded domains.
The extension to bounded domains is presented in \cite{Agnelli2015}.
Further literature review on this approach can in found in \cite{Bellomo2012}.

In this work, we consider the model proposed in \cite{Agnelli2015}. We first validate
it against experimental data and then extend it to bounded domains with obstacles. 
It is worth noticing that most of the models and methods in the references cited so far 
have been shown to reproduce phenomena of pedestrian movement qualitatively
through analysis and/or numerical simulations. However, 
before using a model to predict quantitative results like, e.g., the total evacuation time,
the mathematical models have to be validated and the numerical methods
have to be verified \cite{Einarsson2005}. In the context of pedestrian dynamics, 
this is still difficult due to a lack of reliable experimental data. In addition, the few available datasets show 
large differences \cite{Schadschneider2011,Seyfried2009,Zhang2011}.
In this paper, we compare our numerical results against the data 
reported in a recent empirical study \cite{Kemloh}. We have selected this study
because it deals with egressing from a facility and thus it is the most 
directly related to our focus. With the model under consideration, 
for medium and medium-to-large groups of people we manage to 
achieve good agreement between the computed average people
density and flow rate and the respective measured quantities.
Finally, we mention that in order to make the models more reliable, 
evolutionary adjustment of the parameters and data assimilations have
also been proposed in \cite{JOHANSSON_HELBING,Ward150703},
respectively.

The strategy we propose to handle obstacles within the computational domain 
makes use of an effective obstacle area, which is an enlarged area that encloses the real obstacle, 
and a model parameter used to describe the quality of the environment.  
Thanks to those two ingredients, we can successfully exclude from the walkable
area square and rectangular obstacles. We test our strategy in a
square room that contains one or two obstacles and has one exit. In addition, 
through a series of numerical tests we show that our approach is 
capable of handling evacuation from a room with one or more exits
with variable size, with and
without obstacles, and can reproduce lane formation in bidirectional flow in a corridor.

The paper is organized as follows. Section~\ref{sec:mathematicl_model} introduces the representation of the
system and the modeling of interactions with pedestrians and with the envornment. 
In Section~\ref{sec:Splitting}, we apply the Lie splitting algorithm to the model described in 
Section~\ref{sec:mathematicl_model}. 
Numerical results are presented in  Section~\ref{sec:results} and conclusions are drawn in Sectioin~\ref{sec:concl}.

\section{Mathematical model}\label{sec:mathematicl_model}

The model we consider is based on the model proposed in \cite{Agnelli2015}.
Let  $\Omega \subset \mathbb{R}^2$ denote a bounded domain. 
We assume that the boundary $\partial \Omega$ includes an exit $E$
which could be the finite union of disjoint sets, and walls $W$. Here, $\overline{E}\cup \overline{W}
=\overline{\partial \Omega}$ and $E \cap W =\emptyset$. 
Let $\x=(x,y) \in \Omega$ denote position and ${\v}=v (\cos \theta, \sin \theta) \in \Omega_{\v}$ denote  
velocity, where $v$ is the velocity modulus, $\theta$ is the velocity direction, and
$\Omega_{\v} \subset \mathbb{R}^2$ is the velocity domain. 
For a system composed by a large number of pedestrians distributed inside  
$\Omega$, the distribution function is given by 
\[ f= f(t, \x, \v)\quad \text{for all} \,\,\, t \ge 0,  \,\, \x \in \Omega, \,\, \v \in \Omega_{\v}. \]
Under suitable integrability conditions, $f(t, \x, \v)d \x d \v$ represents the number of individuals who, at time $t$, 
are located in the infinitesimal rectangle $[x, x+dx] \times [y, y+dy]$ and have a velocity belonging to $[v, v + dv] 
\times [\theta, \theta+d\theta]$. 
Since we use polar coordinates for the velocity, we can write the distribution function as $f= f(t, \x, v, \theta)$. 

Following \cite{Agnelli2015}, we assume that variable $\theta$ is discrete. This assumption 
is motivated by the granularity of pedestrian dynamics when the crowd size is not enough to
justify the continuity of the distribution function over the variable $\theta$. For simplicity, 
we assume $\theta$ can take values in the set:
\[ I_{\theta}=\left \{ \theta_{i}= \frac{i-1}{N_d} 2\pi : i = 1, \dots, N_d \right \}, \]
where $N_d$ is the maxim number of possible directions. 
As for the velocity magnitude $v$, we model it as a continuous deterministic variable 
which evolves in time and space according to macroscopic effects determined by the overall dynamics.
In fact, experimental studies show that in practical
situations the speed of pedestrians depends mainly on the level of congestion around
them.

Due to the deterministic nature of the variable $v$, the kinetic type representation is given by the reduced distribution function
\begin{equation}\label{eq:f}
f(t, \x, \theta)= \sum_{i=1}^{N_d} f^i(t, \x)\delta(\theta - \theta_i),
\end{equation}
where $f^i(t, \x)=f(t, \x, \theta_i)$ represents the active particles that, at time $t$  and position $\x$, move with direction $\theta_i$. In equation~\eqref{eq:f}, $\delta$ denotes the Dirac delta function.

Let us introduce some reference quantities that will be use to make the variable dimensionless. 
We define: 
\begin{itemize}
\item[-] $D$: the largest distance a pedestrian can cover in domain $\Omega$;
\item[-] $V_M$: the highest velocity modulus a pedestrian can reach in low density and optimal environmental conditions;
\item[-] $T$: a reference time given by $D/V_M${;}
\item[-] $\rho_M$: the maximal admissible number of pedestrians per unit area.
\end{itemize}
The dimensionless variables are then: position $\hat{\x}=\x/D$, time $\hat{t}=t/T$, velocity modulus $\hat{v}=v/V_M$ and distribution function $\hat{f}=f/ \rho_M$. From now on, all the variables will be 
dimensionless and hats will be omitted to simplify notation. 

Due to the normalization of $f$, and of each $f^i$, the dimensionless local density is obtained by summing the distribution functions over the set of directions:
\begin{align}\label{eq:rho}
\rho(t, \x)=\sum_{i=1}^{N_d}f^i(t, \x) .
\end{align}
As mentioned above, we assume that pedestrians adjust their speed depending on {the} level of congestion around them.
This means that the velocity modulus depends formally on the local density, i.e. $v=v[\rho](t, \x)$, where square brackets are used to denote that $v$ depends on $\rho$ in a functional way. 
For instance, $v$ can depend on $\rho$ and on its gradient. 

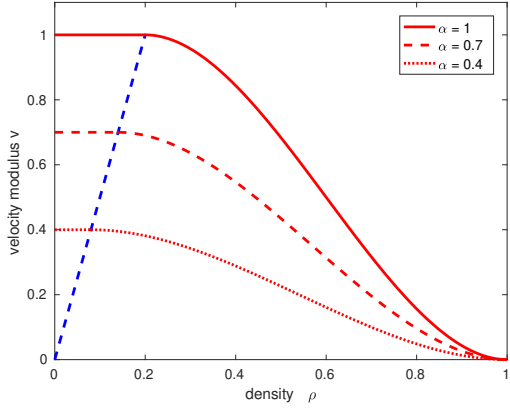
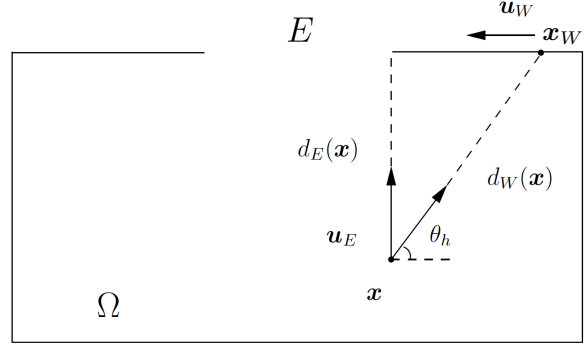
\begin{figure}[h!]
\centering
\subfloat[$v$ as a function of $\rho$]{
\begin{overpic}[width=0.47\textwidth,grid=false]{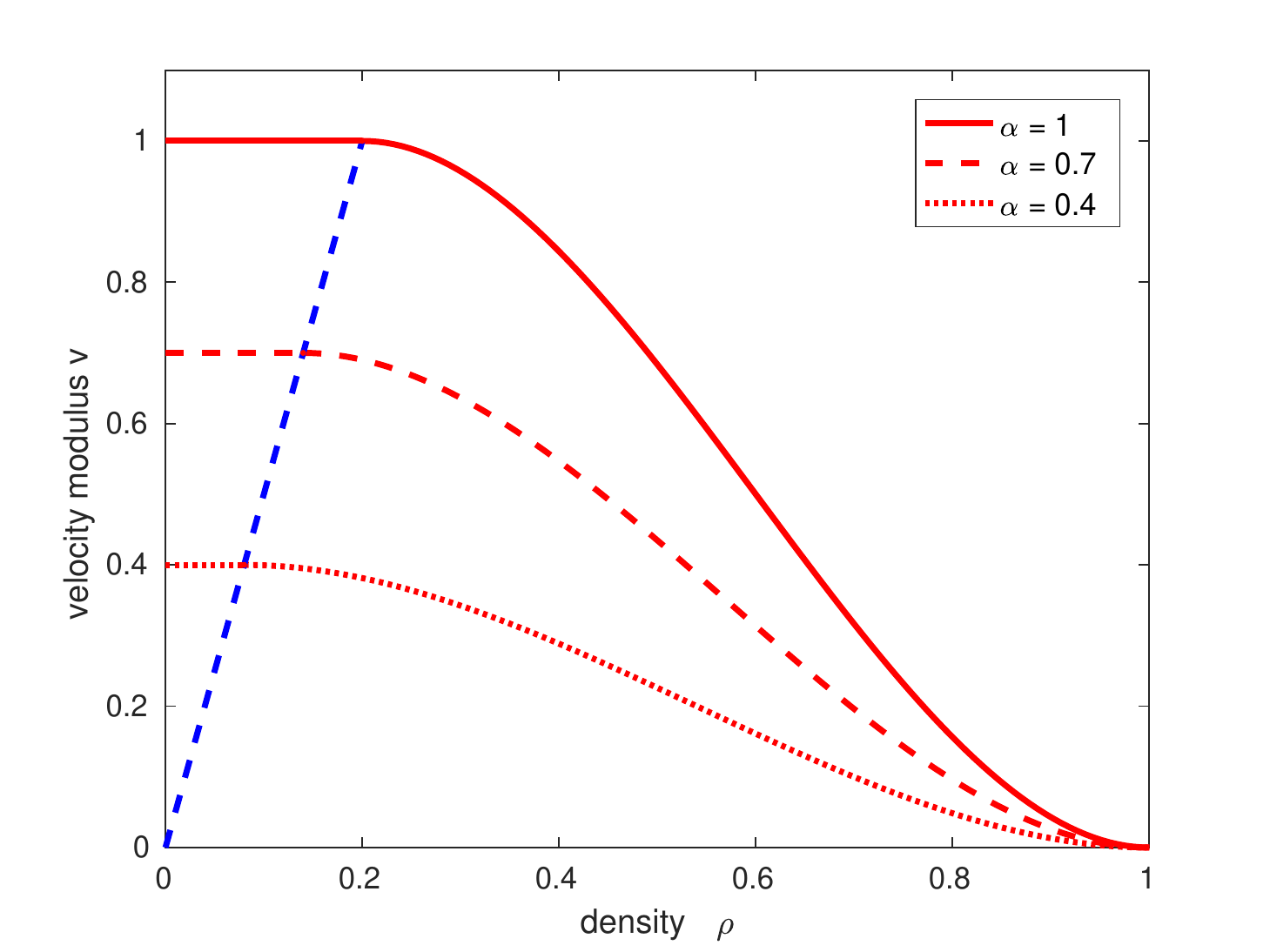}
\end{overpic}
}
\subfloat[Notation in a  sample computational domain]{
\begin{overpic}[width=0.47\textwidth,grid=false]{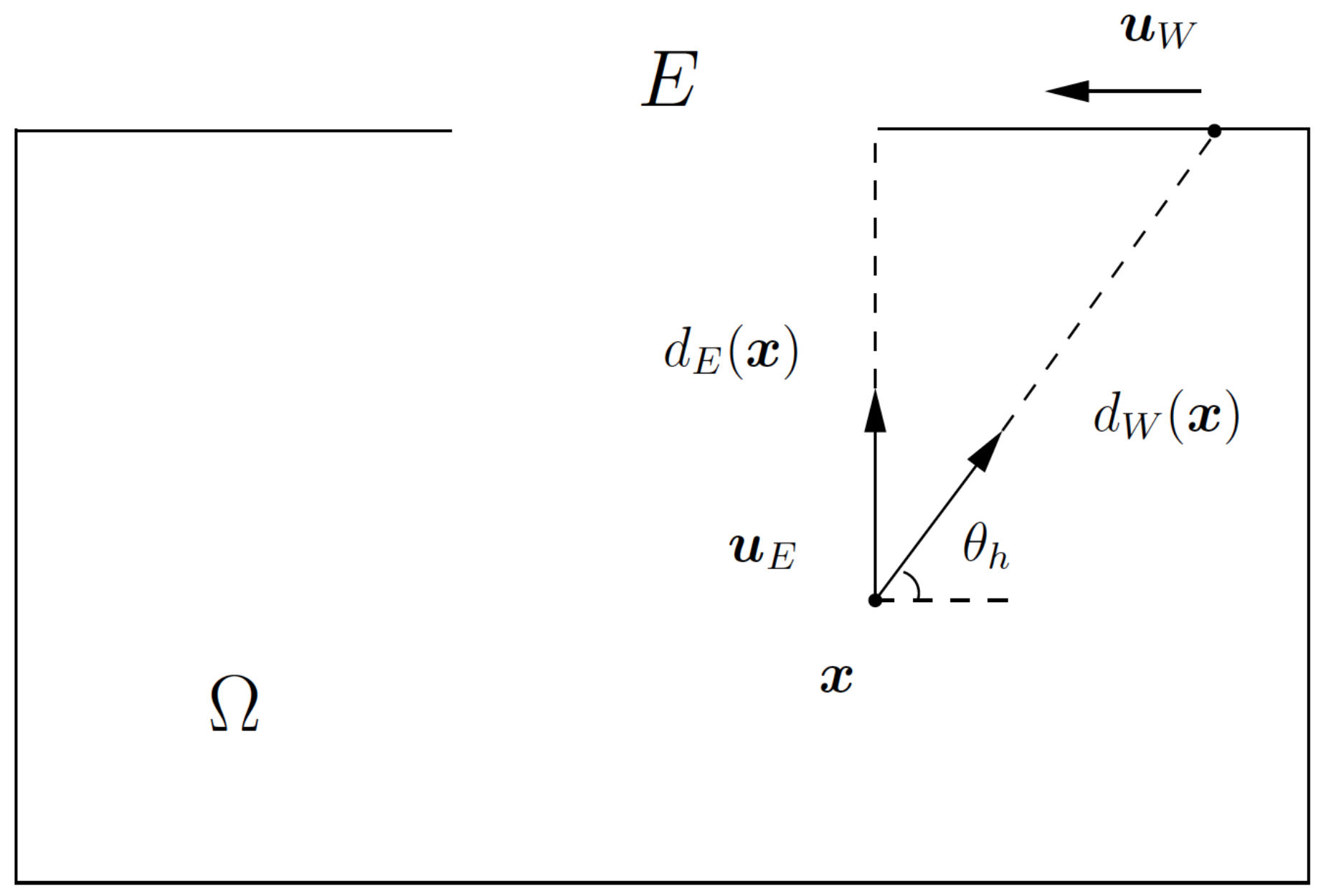}
\put(92, 60){\small{\text{$\x_{W}$}}}
\linethickness{0.45pt}
\put(1.2, 8.1){\line(1,0){97.5}}
\linethickness{3pt}
\put(1, 1){\color{white}{\line(1,0){98}}}
\put(1, 1){\color{white}{\line(0,1){6.98}}}
\put(99, 1){\color{white}{\line(0,1){6.98}}}
\end{overpic}
}
\caption{(A) Dependence of the dimensionless velocity modulus $v$ on the dimensionless density $\rho$ for different values of the parameter $\alpha$, which represents the quality of the environment.
(B) Sketch of computational domain $\Omega$ with exit $E$ and a pedestrian located at $\x$, moving
with direction $\theta_h$. The pedestrian should choose direction $\u_E$ to reach the exit, 
while direction $\u_W$ is to avoid collision with the wall. The distances form the exit and from the wall
are $d_E$ and $d_W$, respectively.}\label{velocity}
\end{figure}

A parameter $\alpha \in [0, 1]$ is introduced to represent the quality of the domain where pedestrians move:
$\alpha=0$ corresponds to the worst quality which forces pedestrians to slow down or stop,  
while $\alpha=1$ corresponds to the best quality, which allows pedestrians to walk at the desired speed.
We assume that the maximum dimensionless speed $v_M$ a pedestrian can reach depends linearly on the quality
of the environment. For simplicity, we take $v_M = \alpha$. Let $\rho_c$ be a critical density 
value such that for $\rho < \rho_c$ we have free flow regime (i.e.,~low density condition), while for $\rho > \rho_c$
we have a slowdown zone (i.e.,~high density condition). We set $\rho_c = \alpha/5$. Note that this choice is 
compatible with the experimentally measured values of $\rho_c$ reported in \cite{Schadschneider2011545}.  
Next, we set the velocity magnitude $v$ equal to $v_M$ in the free flow regime
and equal to a heuristic third-order polynomial in the slowdown zone:
\begin{align}\label{eq:v}
v=v(\rho)=
\begin{cases}
\alpha \quad & \text{for} \quad \rho \leq \rho_c(\alpha)= \alpha/5  \\
a_3\rho^3+a_2\rho^2+a_1\rho+a_0 \quad & \text{for} \quad \rho > \rho_c(\alpha)=\alpha/5   ,     
\end{cases}
\end{align}
where $a_i$ is constant for $i = 0,1,2,3$.
To set the value of these constants, we impose the following 
conditions: $v(\rho_c) = v_M$,  $\partial_{\rho} v(\rho_c) = 0$, $v(1)=0$ and $\partial_{\rho} v(1) = 0$. 
This leads to:
\begin{align}\label{eq:coeff}
\begin{cases}
a_0 &= (1/(\alpha^3-15\alpha^2+75\alpha-125))(75\alpha^2-125\alpha) \\
a_1 &= (1/(\alpha^3-15\alpha^2+75\alpha-125))(-150\alpha^2)\\
a_2 &= (1/(\alpha^3-15\alpha^2+75\alpha-125))(75\alpha^2+375\alpha) \\
a_3 &= (1/(\alpha^3-15\alpha^2+75\alpha-125))(-250\alpha). 
\end{cases}
\end{align}
Figure~\ref{velocity} (A) reports $v$ as a function of $\rho$ for $\alpha = 0.4, 0.7, 1$.

\subsection{Modeling interactions} \label{modelinginteractions}

Each pedestrian is modeled as a particle. Interactions involve three types of particles: 

\begin{itemize}
\item[-] \textit{test particles} with state $(\x, \theta_i)$: they are representative of the whole system;
\item[-] \textit{candidate particles} with state $(\x, \theta_h)$: they can reach in probability the state of the test particles after individual-based interactions with the environment or with field particles; 
\item[-] \textit{field particles} with state $(\x, \theta_k)$: their presence triggers the interactions of the candidate particles.
\end{itemize}
The process through which a pedestrian decides the direction to take is the results of several factors. 
We take into account four factors:

\begin{enumerate}[label={(F\arabic*})]
\item \textit{The goal to reach the exit.}\\
Given a candidate particle at the point $\x$, we define its distance to the exit as
\[d_E(\x)= \min_{\x_E \in E} || \x-\x_E ||,\]
and we consider the unit vector $\u_E(\x)$, pointing from $\x$ to the exit. See Figure~\ref{velocity} (B).

\item \textit{The desire to avoid the collision with walls.}\\
Given a candidate particle at the point $\x$ moving with direction $\theta_h$, we define the distance $d_W(\x, \theta_h)$ from the particle to a wall at a point $\x_W(\x, \theta_h)$ where the particle is expected to collide with the wall.
The unit tangent vector $\u_W(\x, \theta_h)$ to $\partial \Omega$ at $\x_W$
points to the direction of the the exit. Vector $\u_W$ is used to avoid a collision with the walls. See Figure~\ref{velocity} (B).

\item \textit{The tendency to look for less congested area.}\\
A candidate particle $(\x, \theta_h)$ may decide to change direction in order to avoid congested areas. 
This is achieved with the direction that gives the minimal directional derivative of the density at the point $\x$. 
We denote such direction by unit vector $\u_C(\theta_h, \rho)$.

\item \textit{The tendency to follow the stream.}\\
A candidate particle modifies its state, in probability, into that of the test particle due to interactions with field particles, 
while the test particle loses its state as a result of these interactions. A candidate particle $h$ interacting with a field particle $k$  may decide to follow it and thus adopt its direction, denoted with unit vector $\u_F=(\cos\theta_k, \sin\theta_k)$.
\end{enumerate}
\bigskip

Factors (F1) and (F2) are related to geometric aspects of the domain, while factors (F3) and (F4)
consider that people's behavior is strongly affected by surrounding crowd. 
Note that the effects related to factors (F3) and (F4) compete with each other: (F4) is dominant
in a panic situation, while (F3) characterizes rational behavior. As a weight between (F3) and (F4),
 we introduce parameter $\varepsilon \in [0,1]$: $\varepsilon=0$ corresponds to the situation 
 in which only the research of less congested areas is considered (rational behavior), 
 while $\varepsilon=1$ corresponds to the situation in which only the tendency to follow the stream is taken into account (panic behavior). 

\subsubsection{Interaction with the bounding walls}

The interaction with the bounding walls is modeled with two terms:
\begin{compactitem}
\item[-] $\mu[\rho]$: the \textit{interaction rate} models the frequency of interactions between candidate particles and 
the boundary of the domain. If the local density is getting lower, it is easier for pedestrians to see the walls and doors. Thus, we set $\mu[\rho] =1-\rho$.  

\item[-] $\mathcal{A}_h(i)$: the \textit{transition probability} gives the probability that a candidate particle $h$, i.e.~with direction  $\theta_h$, adjusts its direction into that of the test particle $\theta_i$ due to the presence of the walls and/or an exit. 
The following constraint for $\mathcal{A}_h(i)$ has to be satisfied:
\[
\sum_{i=1}^{N_d} \mathcal{A}_h(i)=1 \quad \text{for all} \,\, h \in \{1, \dots, N_d\}.
\]
\end{compactitem}

We assume that particles change direction, in probability, only to an adjacent clockwise or counterclockwise direction in the discrete set $I_\theta$. This means a candidate particle $h$ may end up into the states $h-1, h+1$ or remain in the state $h$. In the case $h=1$, we set $\theta_{h-1}=\theta_{N_d}$ and, in the case $h=N_d$, we set $\theta_{h+1}=\theta_1$.
The set of all transition probabilities $\mathcal{A}=\{\mathcal{A}_h(i) \}_{h,i= 1, \dots, N_d}$ forms the so-called \textit{table of games} that models the game played by active particles interacting with the bounding walls.

To take into account factors (F1) and (F2), we define the vector
\begin{align}\label{eq:uG}
\u_G(\x, \theta_h) &= \frac{(1-d_E(\x))\u_E(\x) + (1-d_W(\x, \theta_h))\u_W(\x, \theta_h)}{|| (1-d_E(\x))\u_E(\x) + (1-d_W(\x, \theta_h))\u_W(\x, \theta_h) ||} \cl
&= (\cos \theta_G, \sin \theta_G).
\end{align}

Here $\theta_G$ is the  \textit{geometrical preferred direction}, which is the ideal direction that a pedestrian should take in order to reach the exit and avoid the walls in an optimal way. Notice that the closer a pedestrian is
to an exit (resp., a wall), the more direction $\u_E$ (resp., $\u_W$) weights. 

A candidate particle $h$ will update its direction by choosing
the angle closest to $\theta_G$ among the three allowed 
angles $\theta_{h-1}, \theta_{h}$ and $\theta_{h+1}$. The transition probability is given by:
\begin{equation}\label{eq:A}
\mathcal{A}_h(i)=\beta_h(\alpha)\delta_{s,i} + (1-\beta_h(\alpha))\delta_{h, i}, \quad i=h-1, h, h+1,
\end{equation}
where
\[s=\argmin_{j \in \{h-1,h+1\}}\{d(\theta_G, \theta_j)\},\]
with
\begin{equation}\label{eq:distance}
 d(\theta_p, \theta_q)=
\begin{cases}
|\theta_p - \theta_q|  & \text{if} \,\,\, |\theta_p - \theta_q | \leq \pi, \\
2\pi - |\theta_p- \theta_q|  & \text{if} \,\,\, |\theta_p- \theta_q| > \pi .             
\end{cases}
\end{equation}
In \eqref{eq:A}, $\delta$ denotes the Kronecker delta function. Coefficient $\beta_h$, proportional to 
parameter $\alpha$, is defined by:
\[
\beta_h(\alpha)=
\begin{cases}
\alpha & \text{if} \,\,\, d(\theta_h, \theta_G) \geq \Delta\theta, \\
\alpha \dfrac{d(\theta_h, \theta_G)}{\Delta\theta}& \text{if} \,\,\, d(\theta_h, \theta_G)< \Delta \theta ,             
\end{cases}
 \]
where $\Delta\theta=2\pi/{N_d}$. The role of $\beta_h$ is to allow for a transition to $\theta_{h-1}$ or
$\theta_{h+1}$ even in the case that the 
geometrical preferred direction $\theta_G$ is closer to $\theta_h$. Such a transition is more likely 
to occur the more distant $\theta_h$ and $\theta_G$ are.
Notice that if $\theta_G=\theta_h$, then $\beta_h=0$ and $\mathcal{A}_h(h)=1$, meaning that a 
pedestrian keeps the same direction (in the absence of interactions other than with the 
environment) with probability 1.

\subsubsection{Interaction with obstacles}\label{sec:obstacles}

The strategy reported in the previous section to avoid collisions with the walls works well
when the pedestrian is sufficiently far from the walls. If pedestrians get too close to the 
bounding walls, and in particular if they are close to an exit, the definition of $\u_G$ in \eqref{eq:uG}
does not prevent collisions with the walls. Thus, obstacles within the domain $\Omega$
cannot be avoided just by adjusting $\u_W$. In this section, we report
an effective strategy to handle obstacles. 

{Four} ingredients are needed to exclude the real obstacle area from the 
walkable domain:
\begin{enumerate}
\item An effective area: an enlarged area that encloses the real obstacle.
\item A definition of $\u_W$ to account for the effective area. 
\item A setting of the parameter $\alpha$ in the effective area depending on the shape of the obstacle.
\item {A dynamic setting of $\u_E$.} 
\end{enumerate}
{The efficacy of this approach is demonstrated numerically in Sec.~\ref{sec:res_obstacle}.
Given the complexity of the pedestrian dynamics model under consideration, 
a formal demonstration that this approach guarantees exclusion
of the obstacle from the walkable area is less straightforward and we shall address it elsewhere.
We note that all the ingredients are already available within the model. 
}

\begin{figure}
\centering
\begin{overpic}[width=0.50\textwidth,grid=false]{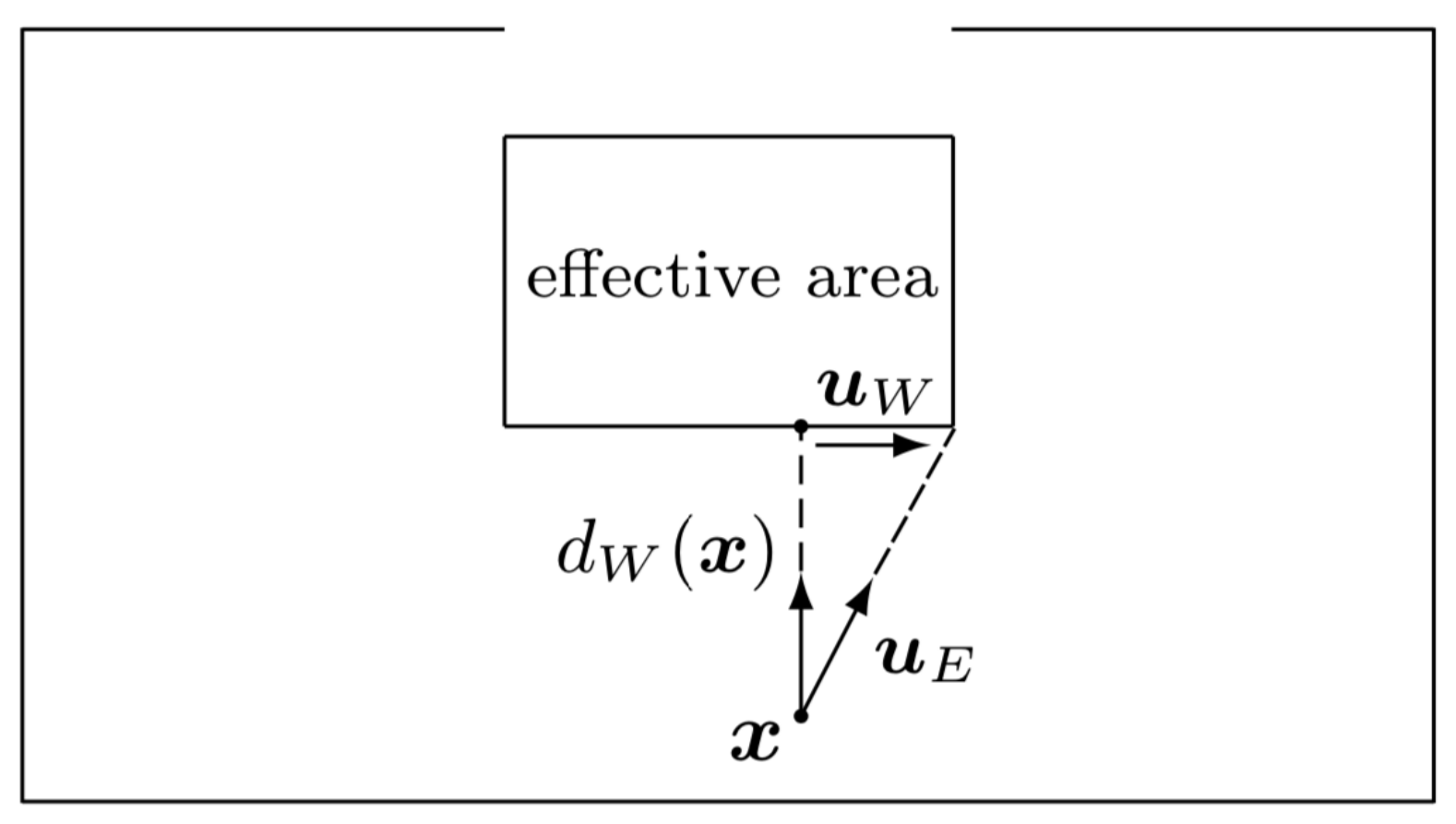}
\put(15,12){\textcolor{black}{\Large$\Omega$}}
\put(48, 55){\textcolor{black}{\Large$E$}}
\end{overpic}
\caption{{Definition of $\u_W$ and $\u_E$ with respect to the effective area.}} \label{fig:uW_obstacle}
\end{figure}

The effective area is necessary especially if the obstacle is close to an exit: it allows to 
define $\u_W$ with respect to a larger area than the one occupied by the obstacle to 
achieve the goal of having no pedestrian walking on the real obstacle area. 
{See Fig.~\ref{fig:uW_obstacle}.}
In the numerical results reported in Section~\ref{sec:res_obstacle}, we used 
an effective area that is four times bigger than the real obstacle area.
{The size and shape of the effective area has been determined heuristically.}

Since some pedestrians will walk on part of the effective area, one needs to set parameter $\alpha$. 
By setting $\alpha=1$ (i.e., best quality of the environment) in the effective area, 
pedestrians can move with the maximal velocity modulus as they approach the obstacle
and thus they quickly adapt to the effective area through $\u_W$. However, some pedestrians
will walk close to the top, bottom, and rear (with respect to the pedestrian motion) 
boundary of the effective area. Thus, the real obstacle is located at the front of the effective area.
From the numerical results reported in Section~\ref{sec:res_obstacle}, we also see that the shape
of the obstacle is square. By setting $\alpha=0$ (i.e., worst quality of the environment) in the effective area, 
pedestrians are forced to slow down at the front part of the effective area. The slow down leads
to higher densities in the front part of the effective area, therefore direction $\u_W$ competes
with direction $\u_C$. As a result some pedestrians walk on the front part of the effective area.
However, as the congestion decreases pedestrians avoid the rear part of the effective area.
From the numerical results shown in Section~\ref{sec:res_obstacle}, we see that the shape
of the obstacle for $\alpha=0$ in the effective area is slender. 

{Finally, the dynamic setting of $\u_E$ is needed for certain pedestrians,
depending on their position with respect to the obstacle. 
As a pedestrian approaches the effective area of an obstacle,
$\u_E$ connects his/her position to the closest corner of the effective area in 
the direction of the final target (i.e., the exit). See Fig.~\ref{fig:uW_obstacle}.
Notice that $\u_E$
can be regarded as the direction a pedestrian needs to take to reach his/her target. 
Thus, it is reasonable that a pedestrian adjusts the direction when in proximity 
of an obstacle. 
}

\subsubsection{Interactions between pedestrians}

The interaction with other pedestrians is modeled with two terms:
\begin{itemize}
\item[-] $\eta[\rho]$: the \textit{interaction rate} defines the number of binary encounters per unit time. If the local density increases, then the interaction rate also increases. For simplicity, we take $\eta[\rho]= \rho$. 
Notice that unlike the case of classical gas dynamics, this rate is not related to the relative particle velocity.

\item[-] $\mathcal{B}_{hk}(i)[\rho]$: the \textit{transition probability} gives the probability that a candidate particle 
$h$ modifies its direction $\theta_h$ into that of the test particle $i$, i.e. $\theta_i$, due to the research 
of less congested areas and the interaction with a field particle $k$ that moves with direction $\theta_k$. 
The following constrain for $\mathcal{B}_{hk}(i)$ has to be satisfied:
\[
\sum_{i=1}^{N_d} \mathcal{B}_{hk}(i)[\rho]=1 \quad \text{for all} \,\, h, k  \in \{1, \dots, N_d\},
\]
where again the square brackets denote the dependence on the density $\rho$.
\end{itemize}

The game consists in choosing the less congested direction among the three admissible ones.
This direction can be computed for a candidate pedestrian $h$ situated at $\x$, by taking 
\[C=\argmin_{j \in \{h-1, h, h+1\}}\{\partial_j\rho(t, \x)\},\]
where $\partial_j\rho$ denotes the directional derivative of $\rho$ in the direction given by angle $\theta_j$. 
We have $\u_C(\theta_h, \rho)=(\cos\theta_C, \sin\theta_C)$.
As for the tendency to follow the crowd, we set $\u_F=(\cos\theta_k, \sin\theta_k)$. This means that  a
candidate particle follows the direction of a field particle.

To take into account (F3) and (F4), we define the vector
\[\u_P(\theta_h, \theta_k, \rho)= \frac{\varepsilon\u_F+(1-\varepsilon)\u_C(\theta_h, \rho)}
{||\varepsilon\u_F+(1-\varepsilon)\u_C(\theta_h, \rho)||} = (\cos \theta_P, \sin \theta_P),
\]
where the subscript $P$ stands for \textit{pedestrians}. Direction $\theta_P$ is the \textit{interaction-based perferred direction}, obtained as a weighted combination between the trendency to follow the stream and the tendency to avoid crowded zones. 

The transition probability is given by:
\[\mathcal{B}_{hk}(i)[\rho]=\beta_{hk}(\alpha)\rho\delta_{r, i} + (1-\beta_{hk}(\alpha)\rho)\delta_{h,i}, \quad i=h-1,h, h+1,\]
where $r$ and $\beta_{hk}$ are defined by:
\[r=\argmin_{j \in \{h-1, h+1\}} \{d(\theta_P, \theta_j)\},\]
\[
\beta_{hk}(\alpha)=
\begin{cases}
\alpha & \text{if} \,\,\, d(\theta_h, \theta_P) \geq \Delta\theta \\
\alpha \dfrac{d(\theta_h, \theta_P)}{\Delta\theta}& \text{if} \,\,\, d(\theta_h, \theta_P)< \Delta \theta.             
\end{cases}
 \]
We recall that $d(\cdot, \cdot)$ is defined in \eqref{eq:distance}.

\subsection{Equation of balance}

The derivation of the mathematical model can be obtained by a suitable balance of particles in an elementary volume
of the space of microscopic states, considering the net flow into such volume due to
transport and interactions \cite{Agnelli2015}.
Taking into account factors (F1)-(F4), we obtain:
\begin{align}
\frac{\partial f^i}{\partial t} &+ \nabla \cdot \left( \v^i [\rho] (t, \x) f^i(t, \x) \right) \cl
& = \mathcal{J}^i[f](t, \x) \cl
& = \mathcal{J}^i_G[f](t, \x) + \mathcal{J}^i_P[f](t, \x) \cl
& = \mu[\rho] \left( \sum_{h = 1}^n \mathcal{A}_h(i) f^h(t, \x) - f^i(t,\x) \right) \cl
& \quad + \eta[\rho] \left( \sum_{h,k = 1}^n \mathcal{B}_{hk}(i) [\rho] f^h(t, \x)f^k(t, \x) - f^i(t,\x) \rho(t, \x)\right)\label{eq:model}
\end{align}
for $i= 1,2, \dots, N_d$. Functional $\mathcal{J}^i[f]$ represents the net balance of particles 
that move with direction $\theta_i$ due to interactions.
As explained in the previous subsection, we consider both the interaction with the environment and with the surrounding people. 
Thus, we can write $\mathcal{J}^i$ as $\mathcal{J}^i=\mathcal{J}^i_G + \mathcal{J}^i_P$, where $\mathcal{J}^i_G$ is an interaction between candidate particles and the environment and 
 $\mathcal{J}^i_P$ is an interaction between candidate and field particles.
 
Equation~\eqref{eq:model} is completed with equation~\eqref{eq:rho}
for the density and equation~\eqref{eq:v},\eqref{eq:coeff} for the velocity. 
In the next section, we will discuss a numerical method for the solution of
problem \eqref{eq:rho},\eqref{eq:v},\eqref{eq:coeff},\eqref{eq:model}.

\section{Numerical method} \label{sec:Splitting}

The approach we consider is based on a splitting method that decouples
the treatment the transport term and the interaction term in equation~\eqref{eq:model}.
As usual with splitting methods, the idea is to split the model into a set of subproblems that are easier to solve
and for which practical algorithms are readily available. The numerical method is then completed
by picking an appropriate numerical scheme for each subproblem.
Among the available operator-splitting methods, we
chose the Lie splitting scheme because it provides
a good compromise between accuracy and robustness, 
as shown in \cite{glowinski2003finite}.

\subsection{The Lie operator-splitting scheme}

Although the Lie splitting scheme is quite well-known, it may be useful to present briefly this scheme before applying it to the solution of problem \eqref{eq:rho},\eqref{eq:v},\eqref{eq:coeff},\eqref{eq:model}.

Let us consider a first-order system in time:
\begin{eqnarray}\label{LieProblem}
   \frac{\partial \phi}{\partial t} + A(\phi) &=& 0, \quad \textrm{in} \ (0,T), \cl
\phi(0) &=& \phi_0, \el
\end{eqnarray}
where A is an operator from a Hilbert space into itself. Operator A is then split, in a non-trivial decomposition, as
\begin{equation*}
 A = \sum\limits_{i=1}^I A_i.
\end{equation*}
The Lie scheme consists of the following. Let $\Delta t>0$ be a time discretization step for the time interval $[0, T]$. Denote $t^k=k\Delta t$, with $k = 0, \dots, N_t$ and let $\phi^k$
be an approximation of $\phi(t^k).$ Set $\phi^0=\phi_0.$ For $n \geq 0$, compute $\phi^{k+1}$ by solving
\begin{eqnarray}
   \frac{\partial \phi_i}{\partial t} + A_i(\phi_i) &=& 0 \quad \textrm{in} \; (t^k, t^{k+1}), \\
\phi_i(t^k) &=& \phi^{k+(i-1)/I}, 
\end{eqnarray}
and then set $\phi^{k+i/I} = \phi_i(t^{k+1}),$ for $i=1, \dots. I.$

This method is first-order accurate in time. More precisely, if~\eqref{LieProblem} 
is defined on a finite-dimensional space, and if the operators $A_i$ are smooth enough, 
then $\| \phi(t^k)-\phi^k \| = O(\Delta t)$~\cite{glowinski2003finite}.

In the next section, we will apply Lie splitting to problem \eqref{eq:model}. The
whole problem will be split into three subproblems: 
\begin{enumerate}
 \item A pure advection problem in the $x$ direction.
 \item A pure advection problem in the $y$ direction.
 \item A problem involving the interaction with the environment and other pedestrians.
\end{enumerate}

\subsection{Lie scheme applied to problem \eqref{eq:model}}\label{sec:Lie_applied}

Let us apply the Lie operator-splitting scheme described in the previous section to problem \eqref{eq:model}.
Given an initial condition $f^{i,0}=f^i(0, \x)$, for $i = 1, \dots, N_d$, the algorithm reads:
For $k=0,1,2, \dots, N_t-1,$ perform the following steps:
\begin{itemize}
\item[-] {\bf Step 1}: Find $f^i$, for $i = 1, \dots, N_d$, such that\\
\begin{equation}
\begin{cases}
 \dfrac{\partial f^i}{\partial t} + \dfrac{\partial }{\partial x} \left( (v [\rho] \cos\theta_i) f^i(t, \x) \right)=0
   \,\,\, \text{on } (t^k, t^{k+1}), \\ \label{eq:step1}
f^i(t^k, \x)=f^{i,k}.
\end{cases} 
\end{equation}
 Set $f^{i,k+\frac{1}{3}}=f^i(t^{k+1}, \x)$.
 
\bigskip
\item[-] {\bf Step 2}:  Find $f^i$, for $i = 1, \dots, N_d$, such that \\
\begin{equation}
\begin{cases}
 \dfrac{\partial f^i}{\partial t} + \dfrac{\partial }{\partial y} \left( (v [\rho] \sin\theta_i) f^i(t, \x) \right)=0
   \,\,\, \text{on } (t^k, t^{k+1}),  \\ \label{eq:step2}
f^i(t^k, \x)=f^{i, k+\frac{1}{3}}. 
\end{cases}
\end{equation}
 Set $f^{i, k+\frac{2}{3}}=f^i(t^{k+1}, \x)$.

\bigskip
\item[-] {\bf Step 3}:  Find $f_i$, for $i = 1, \dots, N_d$, such that\\
\begin{equation}
\begin{cases}
 \dfrac{\partial f^i}{\partial t} = \mathcal{J}^i[f](t, \x)  \,\,\, \text{on } (t^k, t^{k+1}),  \\ \label{eq:step3}
f^i(t^k, \x)=f^{i, k+\frac{2}{3}}.
\end{cases}
\end{equation}
 Set $f^{i,k+1}=f^i(t^{k+1}, \x)$.
\end{itemize}
\bigskip

Notice that once $f^{i, k+1}$ is computed for $i = 1, \dots, N_d$, we use equation~\eqref{eq:rho} to get the density $\rho^{k+1}$ and
equation~\eqref{eq:v},\eqref{eq:coeff} to get the velocity magnitude at time $t^{k+1}$.

\subsection{Space and time discretization} \label{sec:discretization}

Let us assume for simplicity that 
the computational domain under consideration is a rectangle $[0, L] \times [0, H]$, for given $L$ and $H$. 
We mesh the computational domain by choosing $\Delta x$ and $\Delta y$ to partition interval $[0, L]$
and $ [0, H]$, respectively. Let $N_x = L/\Delta x$ and $N_y = H/\Delta y$.
We define the discrete mesh points $\x_{pq} = (x_p, \, y_q)$ by
\begin{align}
x_p &=p \Delta x, \quad p= 0, 1, \dots, N_x,    \cl
y_q  &=q \Delta y, \quad q= 0, 1, \dots, N_y.    \el
\end{align}
It will also be useful to define 
\begin{align}
x_{p+1/2}=x_{p}+\Delta x/2=\Big(p+\frac{1}{2}\Big)\Delta x, \cl
y_{q+1/2}=y_{q}+\Delta y/2=\Big(q+\frac{1}{2}\Big)\Delta y. \el
\end{align}

In order to simplify notation, let us set $\phi = f^i$, $\theta = \theta_i$, $t_0 = t^k$, $t_f = t^{k+1}$.
Let $M$ be a positive integer ($\geq 3$, in practice). We associate with $M$ a time discretization step $\tau = (t_f - t_0)/M$
and set $t^m = t_0 + m \tau$. Next, we proceed with the space and time discretization of each subproblem
in Section~\ref{sec:Lie_applied}. \\

\noindent{\bf Step 1}

\noindent Let $\phi_0 = f^{i,k}$.
Problem \eqref{eq:step1} can be rewritten as
\begin{equation}
\begin{cases}
 \dfrac{\partial \phi}{\partial t} + \dfrac{\partial }{\partial x} \left( (v [\rho] \cos\theta) \phi(t, \x) \right)=0
   \,\,\, \text{on } (t_0, t_f), \\ \label{eq:step1_bis}
\phi(t_0, \x)= \phi_0.
\end{cases} 
\end{equation}
The finite difference method we use produces an approximation $\Phi_{p,q}^{m} \in \mathbb{R}$ 
of the cell average 
\[
\Phi_{p,q}^{m} \approx \dfrac{1}{\Delta x \, \Delta y} \int_{y_{q-1/2}}^{y_{q+1/2}}  \int_{x_{p-1/2}}^{x_{p+1/2}}
\phi(t^m, x, y) dx\, dy, 
\]
where $m=1, \dots, M$, $1 \leq p \leq N_x-1$ and $1 \leq q \leq N_y-1$. 
Given an initial condition $\phi_0$,  function $\phi^m$ will be approximated by $\Phi^{m}$ with
\[
\Phi^m \bigg|_{[x_{p-1/2}, \, x_{p+1/2}] \times [y_{q-1/2}, \, y_{q+1/2}]} = \Phi_{p,q}^{m}
\]

The Lax-Friedrichs method for problem \eqref{eq:step1_bis} can be written in conservative form as follows:
\[\Phi_{p,q}^{m+1}=\Phi_{p,q}^{m}- \dfrac{\tau}{\Delta x}\Big(  \mathcal{F}(\Phi_{p,q}^{m}, \Phi_{p+1,q}^{m}) 
- \mathcal{F}(\Phi_{p-1,q}^{m}, \Phi_{p,q}^{m}) \Big)\]
where
\[ \mathcal{F}(\Phi_{p,q}^{m}, \Phi_{p+1,q}^{m}) =\dfrac{\Delta x}{2\tau}(\Phi_{p,q}^{m}- \Phi_{p+1,q}^{m}) + \dfrac{1}{2} \Big( (v [\rho^m_{p,q}] \cos\theta) \Phi_{p,q}^{m}+(v [\rho_{p+1,q}^m] \cos\theta) \Phi_{p+1,q}^{m} \Big).  \]
\bigskip

\noindent{\bf Step 2}

\noindent Let  $\phi_0 = f^{i,k+\frac{1}{3}}$.
Problem \eqref{eq:step2} can be rewritten as
\begin{equation}
\begin{cases}
 \dfrac{\partial \phi}{\partial t} + \dfrac{\partial }{\partial y} \left( (v [\rho] \sin\theta) \phi(t, \x) \right)=0
   \,\,\, \text{on } (t_0, t_f), \cl
\phi(t_0, \x)= \phi_0 \cl
\end{cases}
\end{equation}
Similarly to step 1, we use the conservative Lax-Friedrichs scheme:
\[\Phi_{p,q}^{m+1}=\Phi_{p,q}^{m}- \dfrac{\tau}{\Delta y}\Big(  \mathcal{F}(\Phi_{p,q}^{m}, \Phi_{p,q+1}^{m}) 
- \mathcal{F}(\Phi_{p,q-1}^{m} \Phi_{p,q}^{m}) \Big)\]
where
\[ \mathcal{F}(\Phi_{p,q}^{m}, \Phi_{p,q+1}^{m}) =\dfrac{\Delta y}{2\tau}(\Phi_{p,q}^{m}-\Phi_{p,q+1}^{m}) + \dfrac{1}{2} \Big( (v [\rho^m_{p,q}] \sin\theta) \Phi_{p,q}^{m}+(v [\rho_{p,q+1}^m] \sin\theta) \Phi_{p,q+1}^{m} \Big).  \]
\bigskip

\noindent{\bf Step 3}

\noindent Let $\mathcal{J} = \mathcal{J}^i $ and $\phi_0 = f^{i,k+\frac{2}{3}}$.
Problem \eqref{eq:step3} can be rewritten as
\begin{equation}
\begin{cases}
 \dfrac{\partial \phi}{\partial t} = \mathcal{J}[f](t, \x)  \,\,\, \text{on } (t_0, t_{f}), \cl
\phi(t_0, \x)= \phi_0. \el
\end{cases}
\end{equation}
For the approximation of the above problem, we use the Forward Euler scheme:
\[ \Phi_{p,q}^{m+1}= \Phi_{p,q}^{m} + \tau \Big (\mathcal{J}^m[F^{m}] \Big ), \] 
where $F^m$ is the approximation of the reduced distribution function \eqref{eq:f} at time $t^m$.
\medskip

For stability, the subtime step $\tau$ is chosen to satisfy the Courant-Friedrichs-Lewy (CFL) 
condition (see, e.g., \cite{leveque1992numerical}):
\[\max \Bigg\{\cfrac{\tau}{\Delta x}, \,\, \cfrac{\tau}{\Delta y} \Bigg\}\leq 1. \]

\section{Numerical results} \label{sec:results}

We present a series of numerical results to showcase the features of our model. 
We start from the simulation of evacuation from a room with one exit in
Section~\ref{sec:one_exit}. We use this first test to validate our implementation
of the model presented in Section~\ref{sec:mathematicl_model} against the results
in \cite{Agnelli2015}.
In order to further validate our software, in Section~\ref{sec:two_exits} we compare our numerical results
with the experimental data reported in \cite{Kemloh}, where the authors study
the evacuation from a room with two exits.  Our successfully validated code is then used to
study evacuation from a room with obstacles in Section~\ref{sec:res_obstacle} and lane formation in a corridor 
in Section~\ref{sec:lane}.

For all the simulations, we consider eight different velocity directions $N_d=8$ in the discrete set:
\[ I_{\theta}=\left \{ \theta_{i}= \frac{i-1}{8} 2\pi : i = 1, \dots, 8 \right \}. \]

\subsection{Evacuation from a room with one exit} \label{sec:one_exit}

\begin{figure}[h!]
\centering
\begin{overpic}[width=0.32\textwidth,grid=false,tics=10]{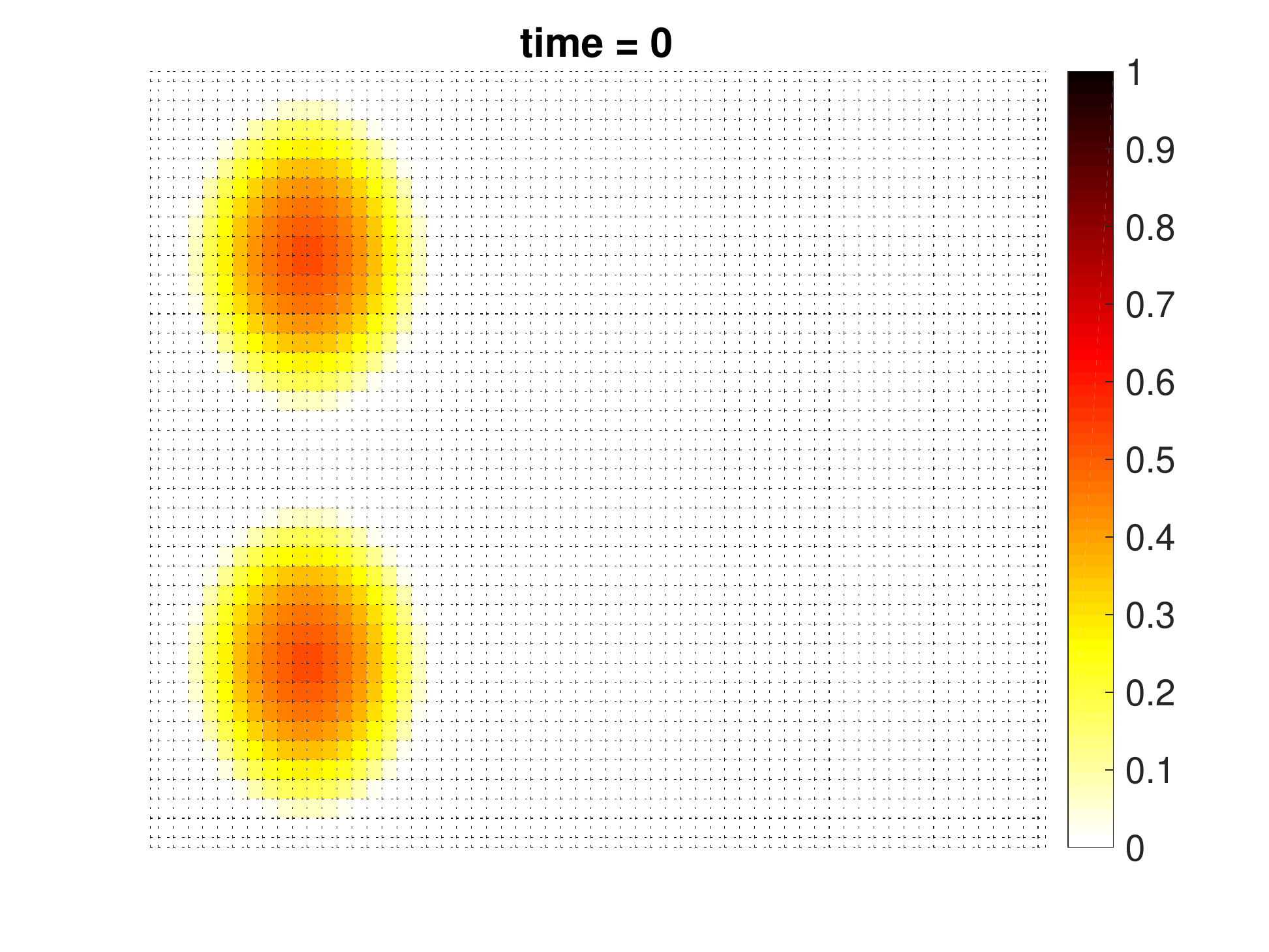}
\linethickness{1pt}
\put(24.25, 54.9){\vector(0,-1){10}}
\put(24.25, 22.7){\vector(0,1){10}}
\put(11.7, 69.5){\line(1,0){71}}
\put(11.7, 7.8){\line(0,1){61.9}}
\put(11.7, 8){\line(1,0){71}}
\put(61.2, 7.8){\line(0,1){23}}
\put(61.2, 46.7){\line(0,1){23}}
\end{overpic} 
\begin{overpic}[width=0.32\textwidth,grid=false,tics=10]{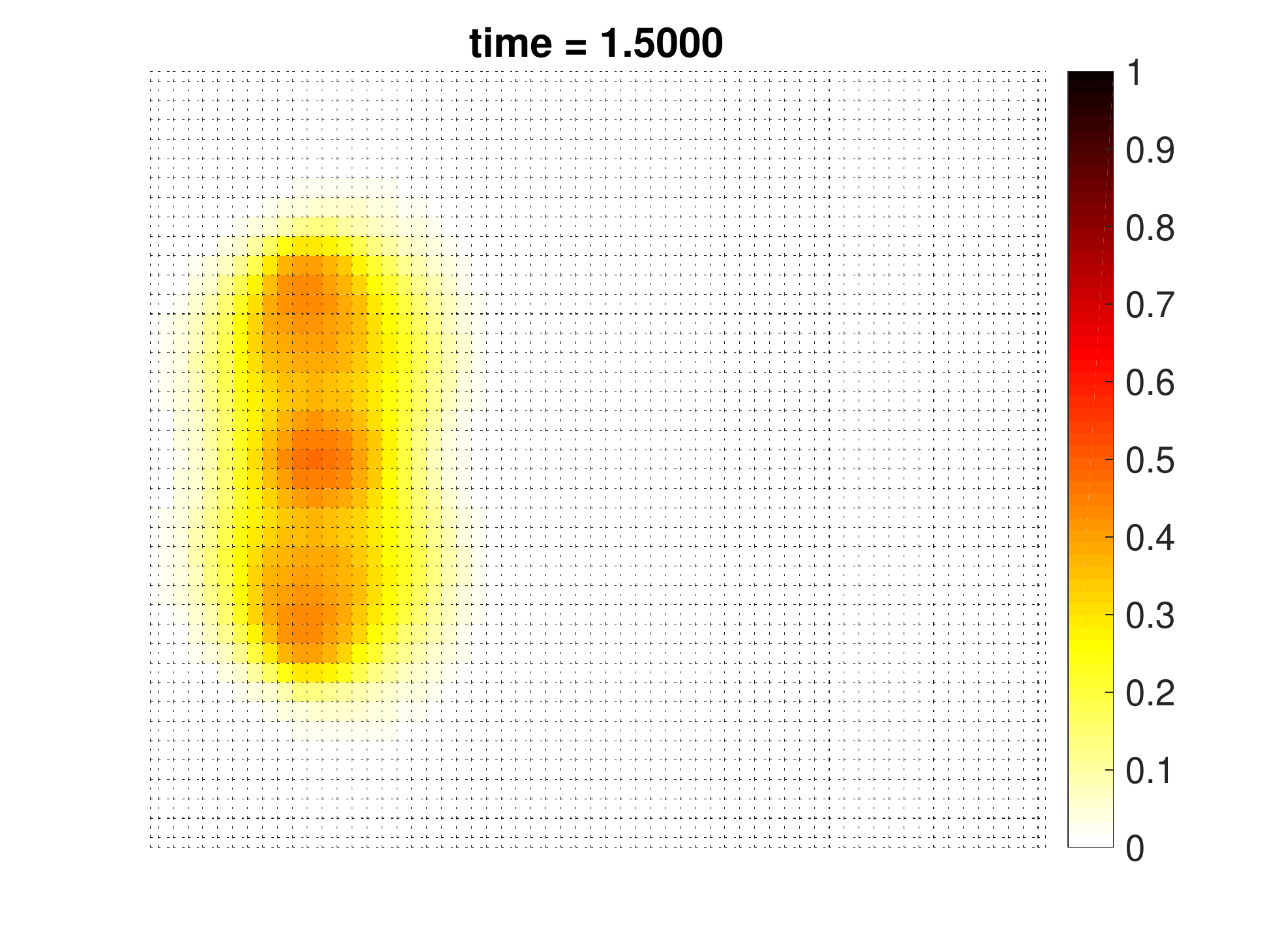}
\linethickness{1pt}
\put(11.7, 69.5){\line(1,0){71}}
\put(11.7, 7.8){\line(0,1){61.9}}
\put(11.7, 8){\line(1,0){71}}
\put(61.2, 7.8){\line(0,1){23}}
\put(61.2, 46.7){\line(0,1){23}}
\end{overpic}
\begin{overpic}[width=0.32\textwidth,grid=false,tics=10]{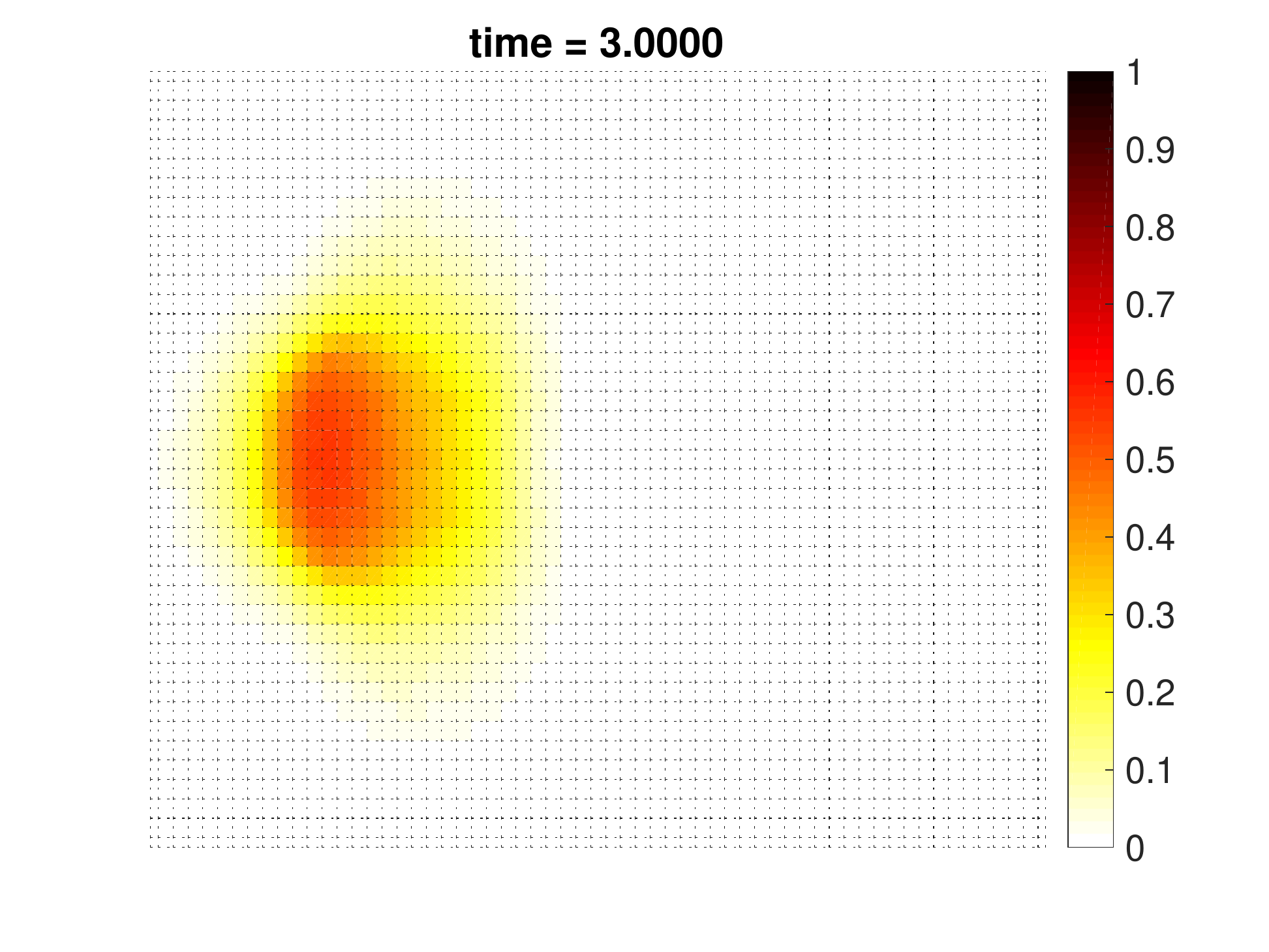}
\linethickness{1pt}
\put(11.7, 69.5){\line(1,0){71}}
\put(11.7, 7.8){\line(0,1){61.9}}
\put(11.7, 8){\line(1,0){71}}
\put(61.2, 7.8){\line(0,1){23}}
\put(61.2, 46.7){\line(0,1){23}}
\end{overpic}
\begin{overpic}[width=0.32\textwidth,grid=false,tics=10]{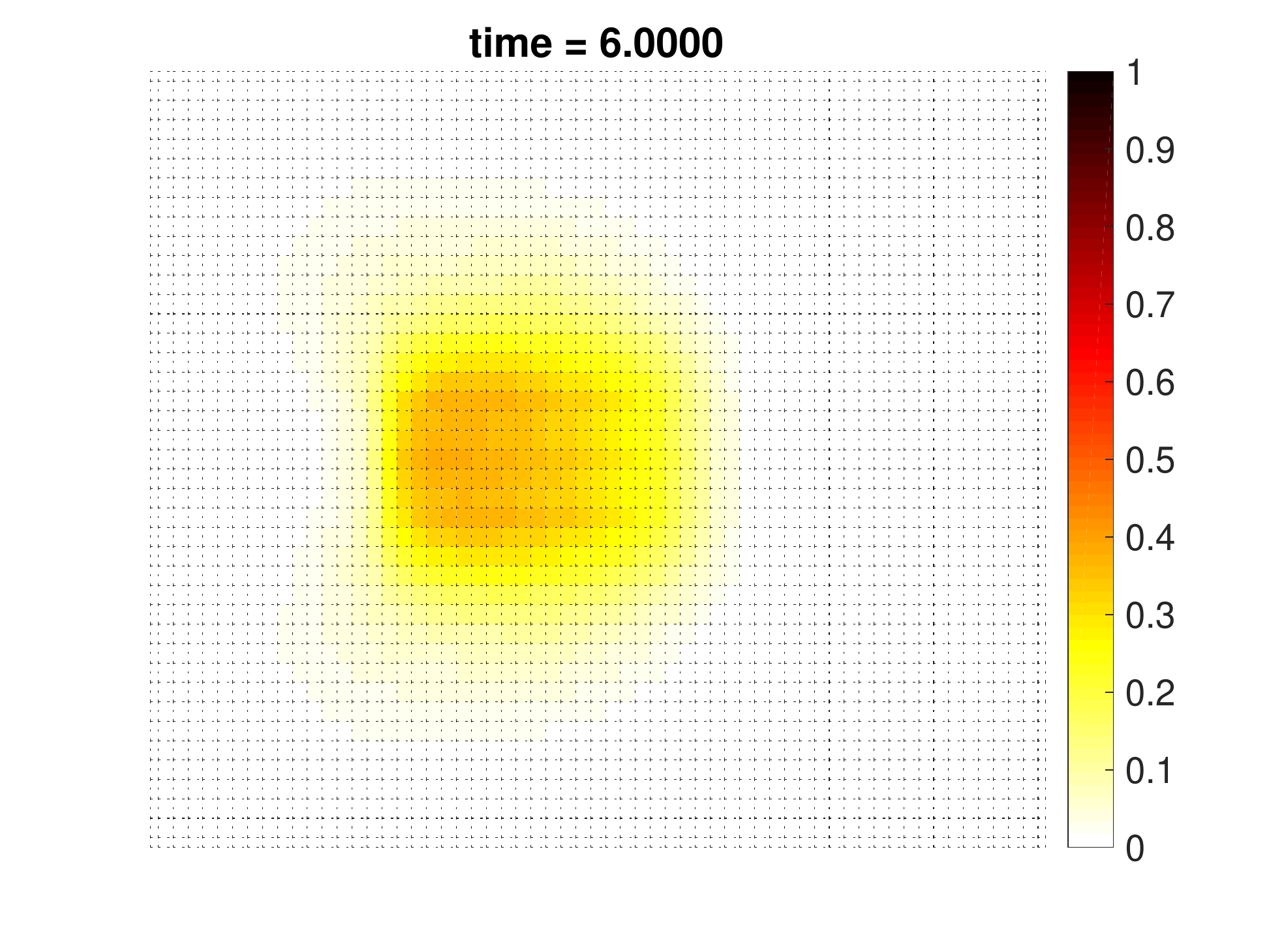}
\linethickness{1pt}
\put(11.7, 69.5){\line(1,0){71}}
\put(11.7, 7.8){\line(0,1){61.9}}
\put(11.7, 8){\line(1,0){71}}
\put(61.2, 7.8){\line(0,1){23}}
\put(61.2, 46.7){\line(0,1){23}}
\end{overpic}
\begin{overpic}[width=0.32\textwidth,grid=false,tics=10]{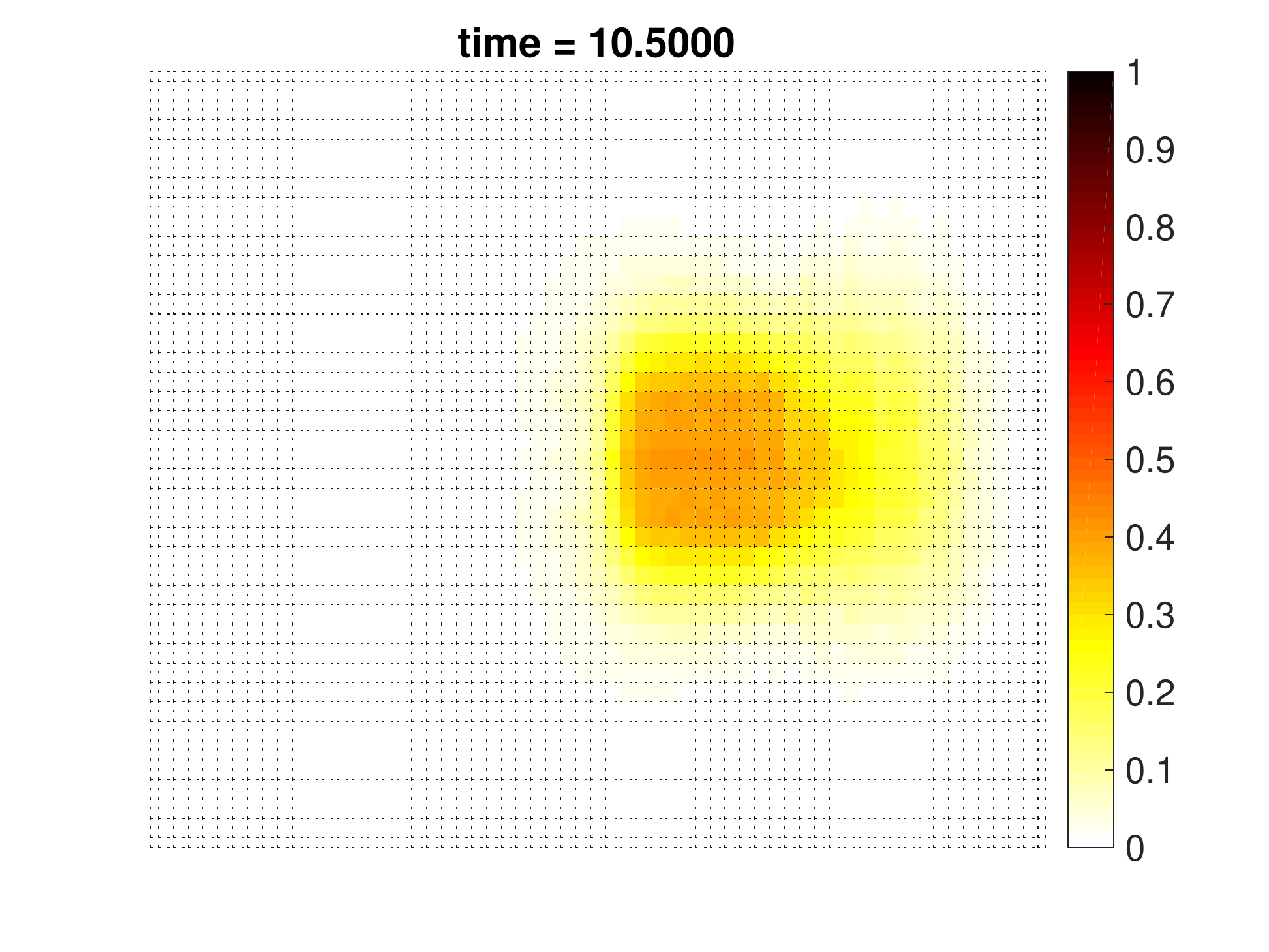}
\linethickness{1pt}
\put(11.7, 69.5){\line(1,0){71}}
\put(11.7, 7.8){\line(0,1){61.9}}
\put(11.7, 8){\line(1,0){71}}
\put(61.2, 7.8){\line(0,1){23}}
\put(61.2, 46.7){\line(0,1){23}}
\end{overpic}
\begin{overpic}[width=0.32\textwidth,grid=false,tics=10]{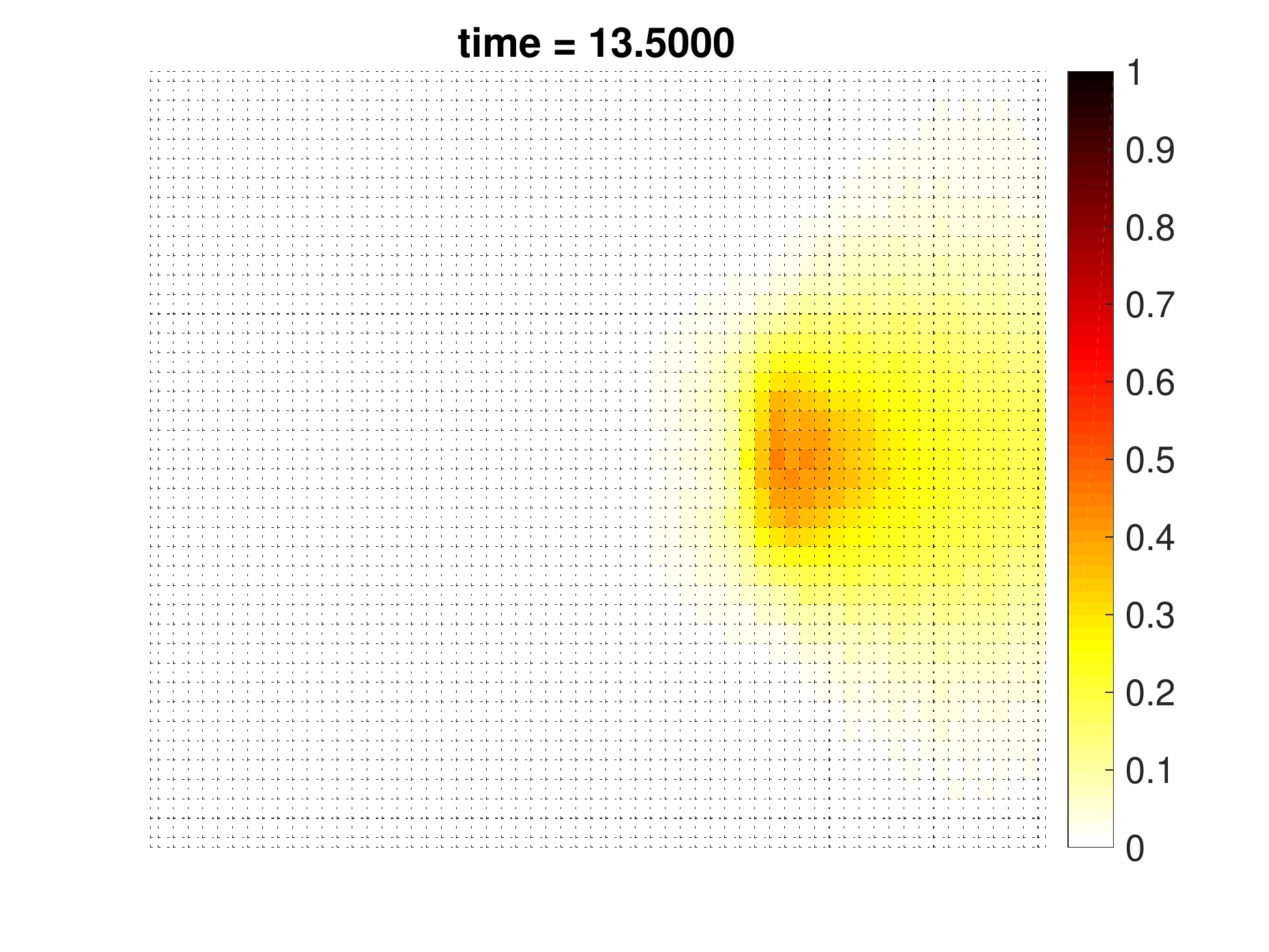}
\linethickness{1pt}
\put(11.7, 69.5){\line(1,0){71}}
\put(11.7, 7.8){\line(0,1){61.9}}
\put(11.7, 8){\line(1,0){71}}
\put(61.2, 7.8){\line(0,1){23}}
\put(61.2, 46.7){\line(0,1){23}}
\end{overpic}
\caption{Evacuation process of 46 pedestrians grouped into two clusters with opposite initial directions $\theta_3$ and $\theta_7$ using the medium mesh and $\Delta t_{medium}$ for time $t = 0, 1.5, 3, 6, 10.5,13.5$ s. The color refers to density.} 
\label{MM}
\end{figure}
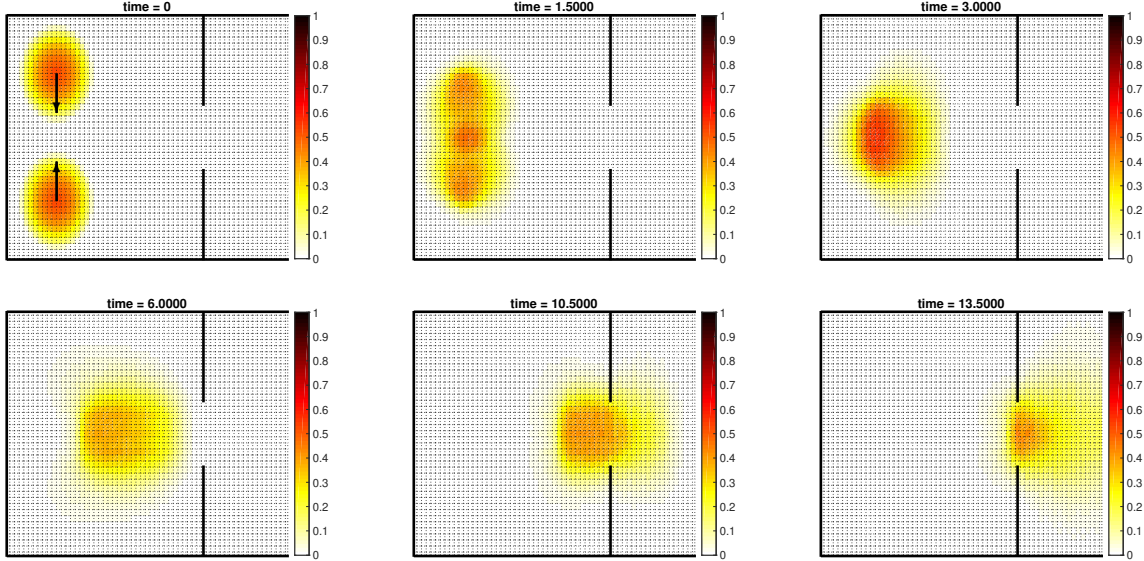

This first test case is taken from \cite{Agnelli2015}.
The computational domain encloses a square room with side $10$ m with
an exit door located in the middle of the right side. The exit size is $2.6$ m. 
The computational domain is larger than the room itself to follow the motion of the pedestrian
also once they have left the room. We aim at simulating the evacuation of 46 people located inside
the room and initially distributed into two equal-area circular clusters. See Figure~\ref{MM}, top
left panel. The two groups are initially moving against the each other with opposite initial directions  
$\theta_3$ and $\theta_7$. Following \cite{Agnelli2015}, simulations are performed with $\varepsilon=0.4$.

In order to work with dimensionless quantities as described in Section~\ref{sec:mathematicl_model}, 
we define the following reference quantities: $D=10\sqrt2$ m, $V_M= 2$ m/s, $T = D/V_M= 5 \sqrt2$ s, and
$\rho_M = 7$ people/m$^2$. However, once the results are computed we convert them back to
dimensional quantities. 

{In order to understand what level of refinement is needed for the mesh,} we consider three different meshes: 
\begin{itemize}
\item[-] \emph{coarse mesh} with $\Delta x = \Delta y=0.5$ m;
\item[-] \emph{medium mesh} with $\Delta x = \Delta y=0.25$ m;
\item[-] \emph{fine mesh} with $\Delta x = \Delta y=0.125$ m.
\end{itemize}
 Similarly, we consider three different time steps: a large time step $\Delta t_{large} = 1.5$ s,  a medium time step $\Delta t_{medium} = 0.75$ s, and a small time step $\Delta t_{small} = 0.375$ s. The value of $M$ for the Lie splitting scheme is set to 3.
Figure~\ref{MM} shows the density computed with medium mesh and $\Delta t_{medium}$ 
at times $t = 0$, 1.5, 3, 6, 10.5, 13.5 s.

\begin{figure}[h!]
\centering
\subfloat[Number of pedestrians vs time]{
\begin{overpic}[width=0.49\textwidth,grid=false]{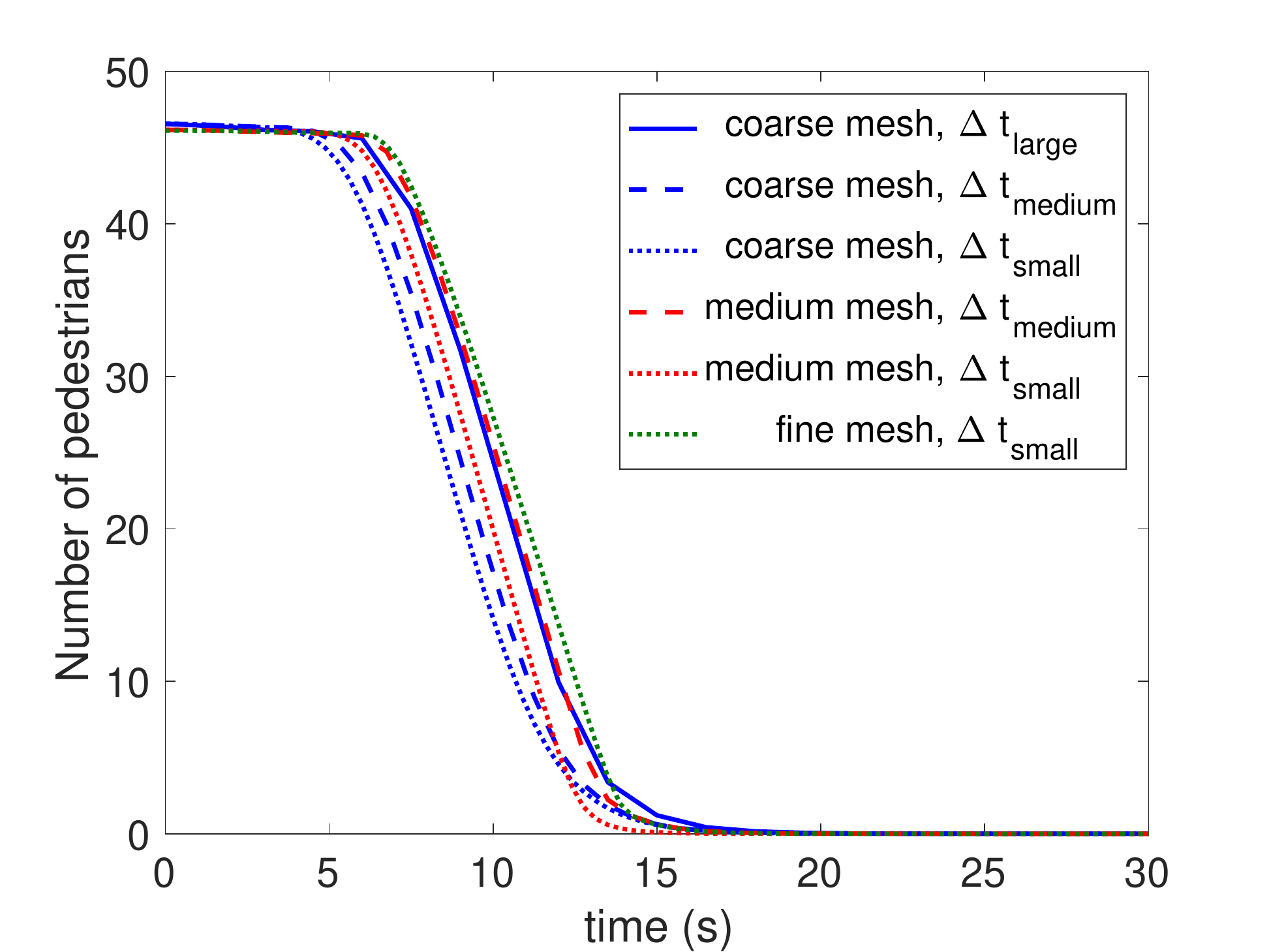}
\end{overpic}
}
\subfloat[Selected curves from (A)]{
\begin{overpic}[width=0.49\textwidth,grid=false]{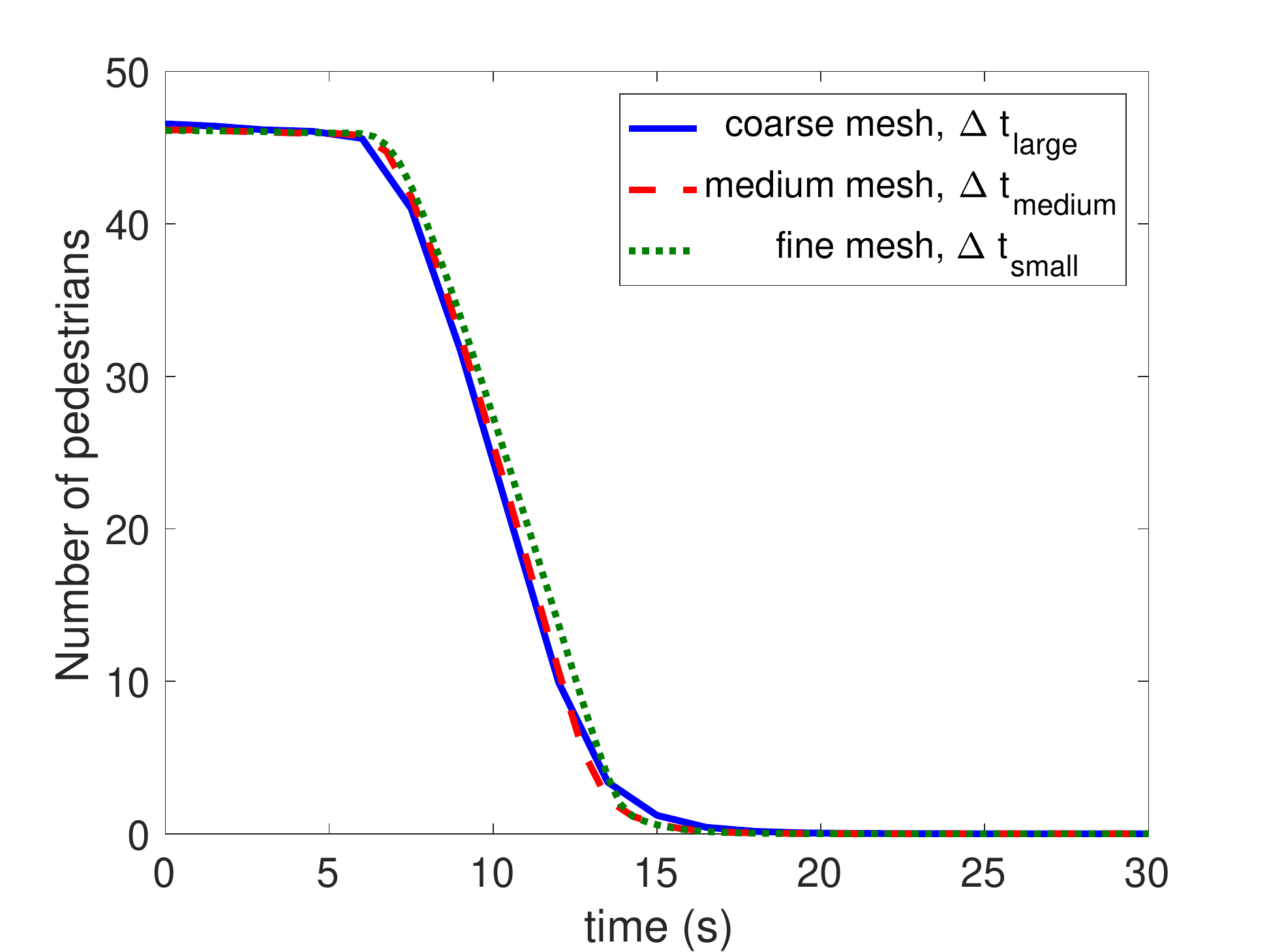}
\end{overpic}
}
\caption{(A) Number of pedestrians inside the room over time computed with six different 
combinations of mesh and time step. For ease of comparison, (B) shows only the curves in (A)
obtained with simultaneous refinements of mesh and time step.}
\label{Population_change}
\end{figure}

Figure~\ref{Population_change} (A) reports the number of pedestrians left in the room 
computed {with the following combinations of mesh and time step:
coarse mesh and $\Delta_{large}$, coarse mesh and $\Delta_{medium}$, coarse mesh and $\Delta_{small}$,
medium mesh and $\Delta_{medium}$, medium mesh and $\Delta_{small}$, and fine mesh and $\Delta_{small}$.
}
In all the cases the total evacuation time is around 18 s, which agrees well
with the results reported in \cite{Agnelli2015}. However, we see that the 
evacuation dynamics varies when the mesh and time step change.
In fact, from Figure~\ref{Population_change} (A) one can observe that as the time step gets 
smaller with a given mesh people walk slightly faster, 
while as the mesh gets finer with a given time step pedestrians walk slightly slower. 
When using operator splitting methods, it is advised to reduce the time
step as the mesh is refined (see, e.g. \cite{glowinski2003finite}). So, for ease of comparison in 
Figure~\ref{Population_change} (B) we report only the results computed with 
coarse mesh and $\Delta t_{large}$, medium mesh and $\Delta t_{medium}$, fine mesh and $\Delta t_{small}$.
We can see very good agreement for those three curves. 

\begin{figure}[h!]
\centering
\subfloat[Our results]{
\centering
\begin{overpic}[width=0.4\textwidth,grid=false]{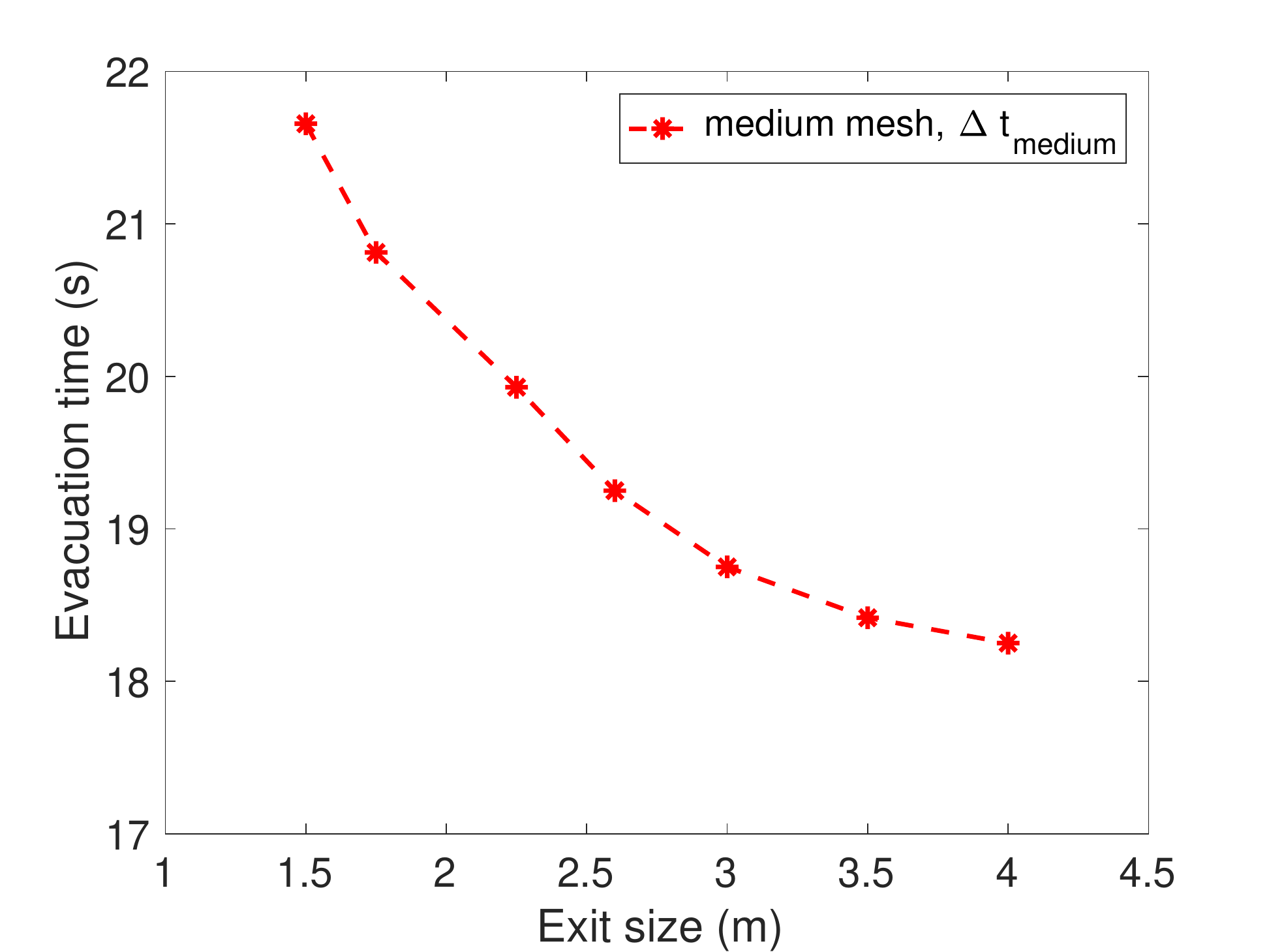}
\end{overpic}
}
\subfloat[{Results from \cite{Agnelli2015}}]{
\begin{overpic}[width=0.4\textwidth,grid=false]{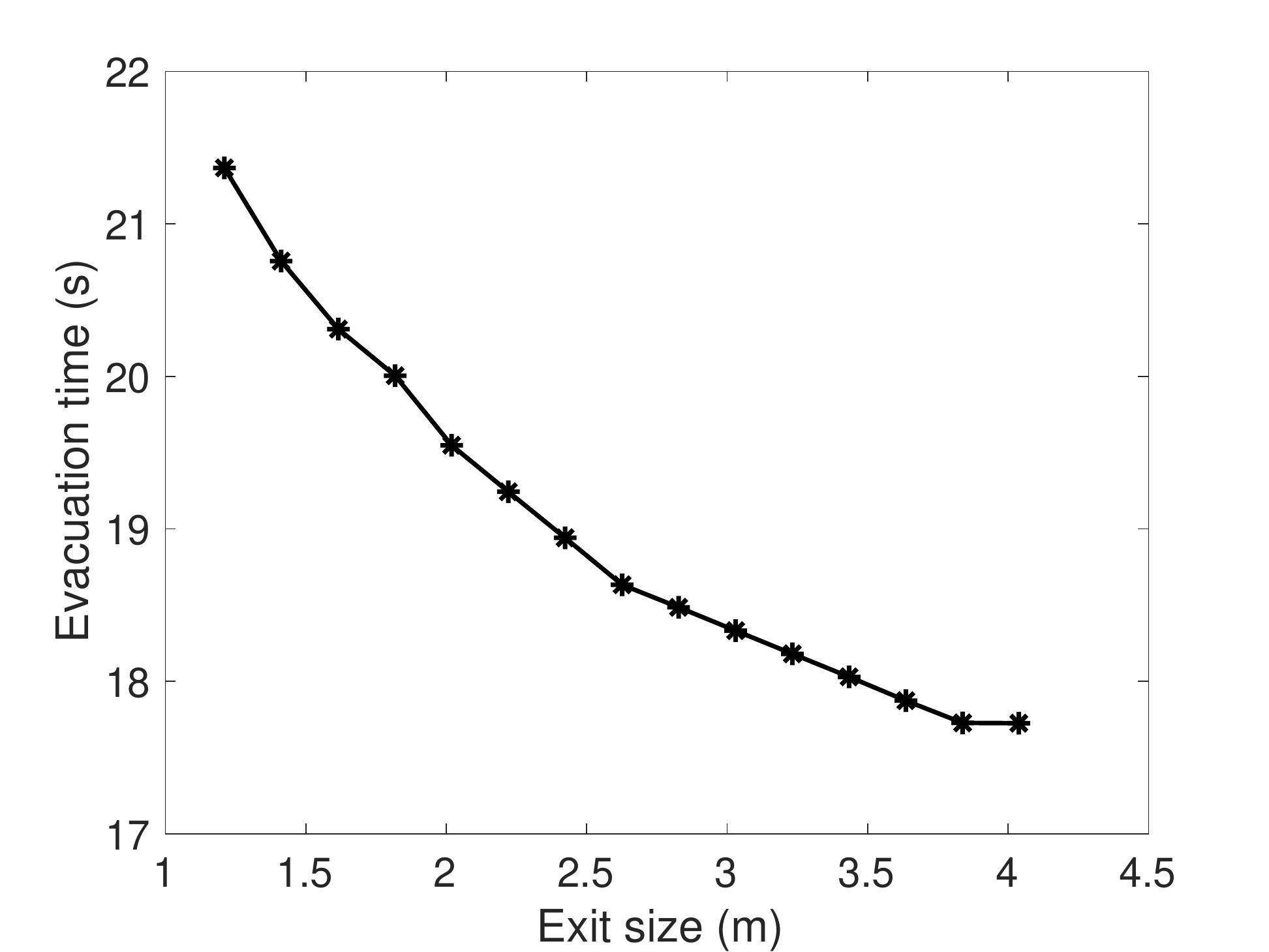}
\end{overpic}
}
\caption{Computed evacuation time from the room with one exit versus the exit size: {(A) our results and
(B) results from \cite{Agnelli2015}}.}
\label{Exit_sizes}
\end{figure}

Next, we consider the same room as in the previous test but we vary the exit size. 
We locate the exit symmetrically with respect to the room centerline and
let the exit size vary from 1.5 m to 4 m. 
We consider the medium meshes and $\Delta t_{medium}$ mentioned before
since Figure~\ref{Population_change} show that is an appropriate choice 
for the problem under consideration. 
All the other model and discretization parameters are set like in the previous test.
 Figure~\ref{Exit_sizes} shows the total evacuation time as a function
of the exit size. First, we notice that our results are in very good agreement with the results 
reported in \cite{Agnelli2015}. As expected, the total evacuation time decreases with the
exit size, but that once the exit is large enough for the crowd contained in the room
the evacuation time does not change significantly if the exit is further enlarged.

As last test for the room with one exit, we computed the total evacuation time
from a room with variable exit size for different values of $\varepsilon$.
The results did not show any visible difference from the ones displayed in 
 Figure~\ref{Exit_sizes} {(A)}, which corresponds to $\varepsilon = 0.4$, and thus
 are not reported here. This test case is probably too simple to see
 any difference in the crowd dynamics when (F3), i.e.~the tendency avoid the crowd, 
 is dominant over (F4), i.e.~the tendency to follow the crowd, and viceversa. 

\subsection{Evacuation from a room with two exits}\label{sec:two_exits}

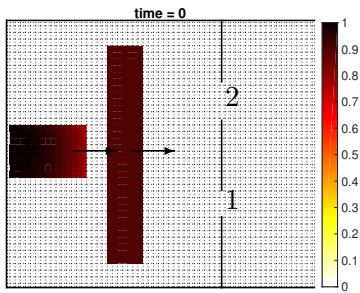
\begin{figure}[h!]
\centering
\begin{overpic}[width=0.35\textwidth,grid=false,tics=10]{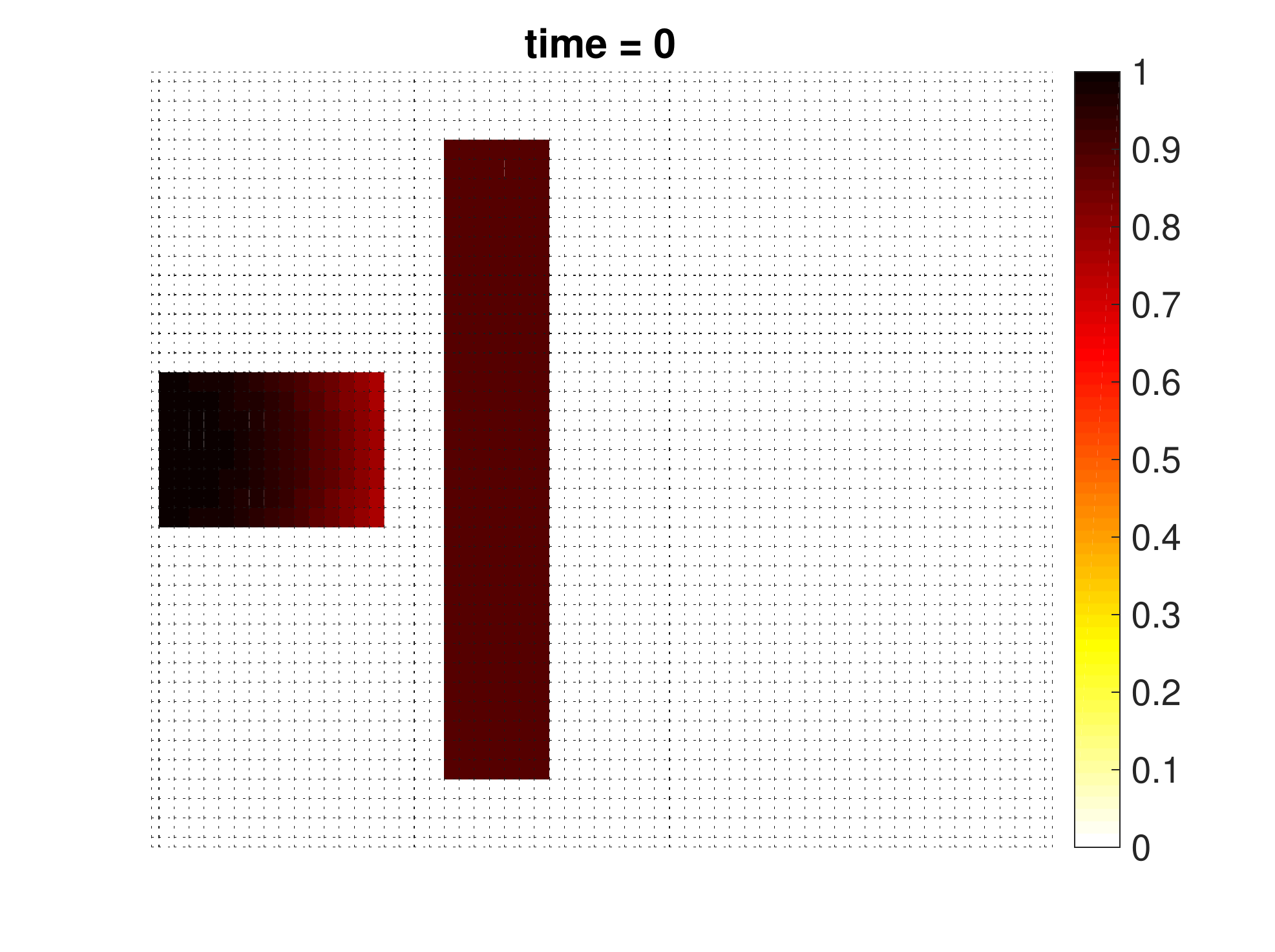}
\linethickness{0.6pt}
\put(27, 39.5){\vector(1,0){10}}
\put(40.5, 39.5){\vector(1,0){10}}
\put(11.7, 69.5){\line(1,0){71}}
\put(11.7, 7.8){\line(0,1){61.9}}
\put(11.7, 8){\line(1,0){71}}
\put(61.2, 7.8){\line(0,1){16.7}}
\put(61.2, 30.2){\line(0,1){16.5}}
\put(61.2, 55.2){\line(0,1){14.5}}
\put(62, 49.7){\text{2}}
\put(62, 26){\text{1}}
\end{overpic} 
\caption{Computational domain corresponding to the experimental set-up in \cite{Kemloh}
and initial density and direction (i.e.,~$\theta_1$) for the experiment with 138 pedestrians. 
}\label{fig:kemloh_138}
\end{figure}

After validating our software against the numerical results in 
\cite{Agnelli2015}, we proceed with the validation against
the experimental data reported in \cite{Kemloh}.
 The computational domain corresponds to the setting studied in \cite{Kemloh}:
a room of side 10 m with different sized exits placed on on the right side: 
exit 1 with size 0.7 m and exit 2 with size 1.1 m. 
The distance between the two exits is 3 m.  See Figure~$\ref{fig:kemloh_138}$.

The data in \cite{Kemloh} refer to ten trials: 2 trials with 18 pedestrians,
6 trials with 40 pedestrians, and 2 trials with 138 pedestrians. 
In the experiments with 18 and 40 pedestrians (and so in the simulation as well), 
the group is initially positioned in a square located symmetrically 
with respect to the horizontal axis of the room and towards the back of the room with 
linearly increasing people density from the front to the back. 
In the experiments with 138 pedestrians, pedestrians are 
initially distributed as follows:
a first group of 90 people is positioned in a rectangle in the middle of the room, 
while the remaining 48 people are positioned in a square behind the first group.
Thus, we adopt the same configuration in the simulations.
We impose constant density in the rectangle and 
prescribe a linearly increasing density from the front to the back of the square. 
See Figure~$\ref{fig:kemloh_138}$.
In all the cases, the pedestrians are given initial direction $\theta_1$. 
For all the simulations, we use the medium mesh and $\Delta t_{small}$
considered in Section~\ref{sec:one_exit}.

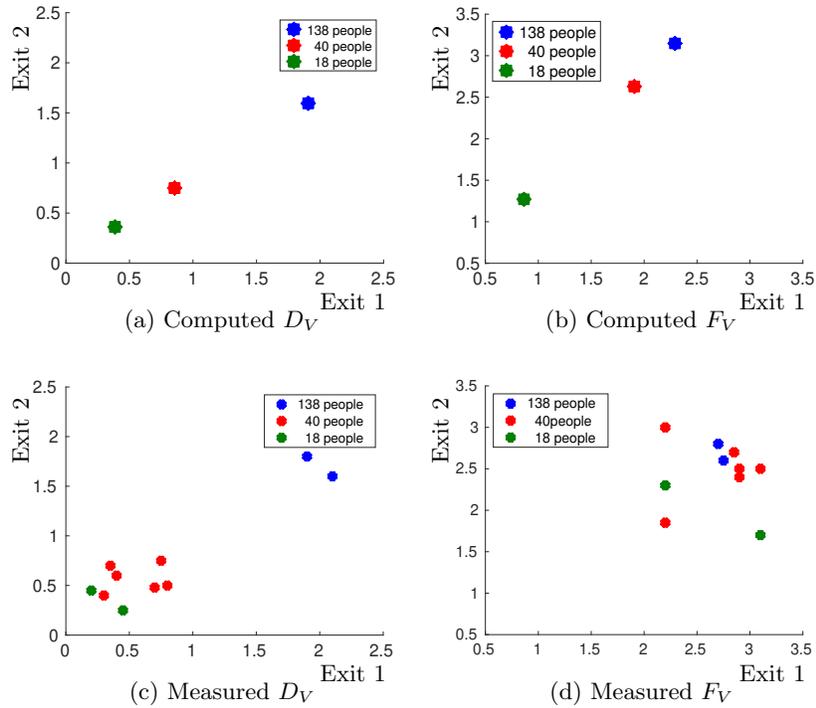
\begin{figure}[h!]
\centering
\subfloat[Computed $D_{V}$]{
\newline
\begin{overpic}[width=0.33\textwidth, grid=false,tics=10]{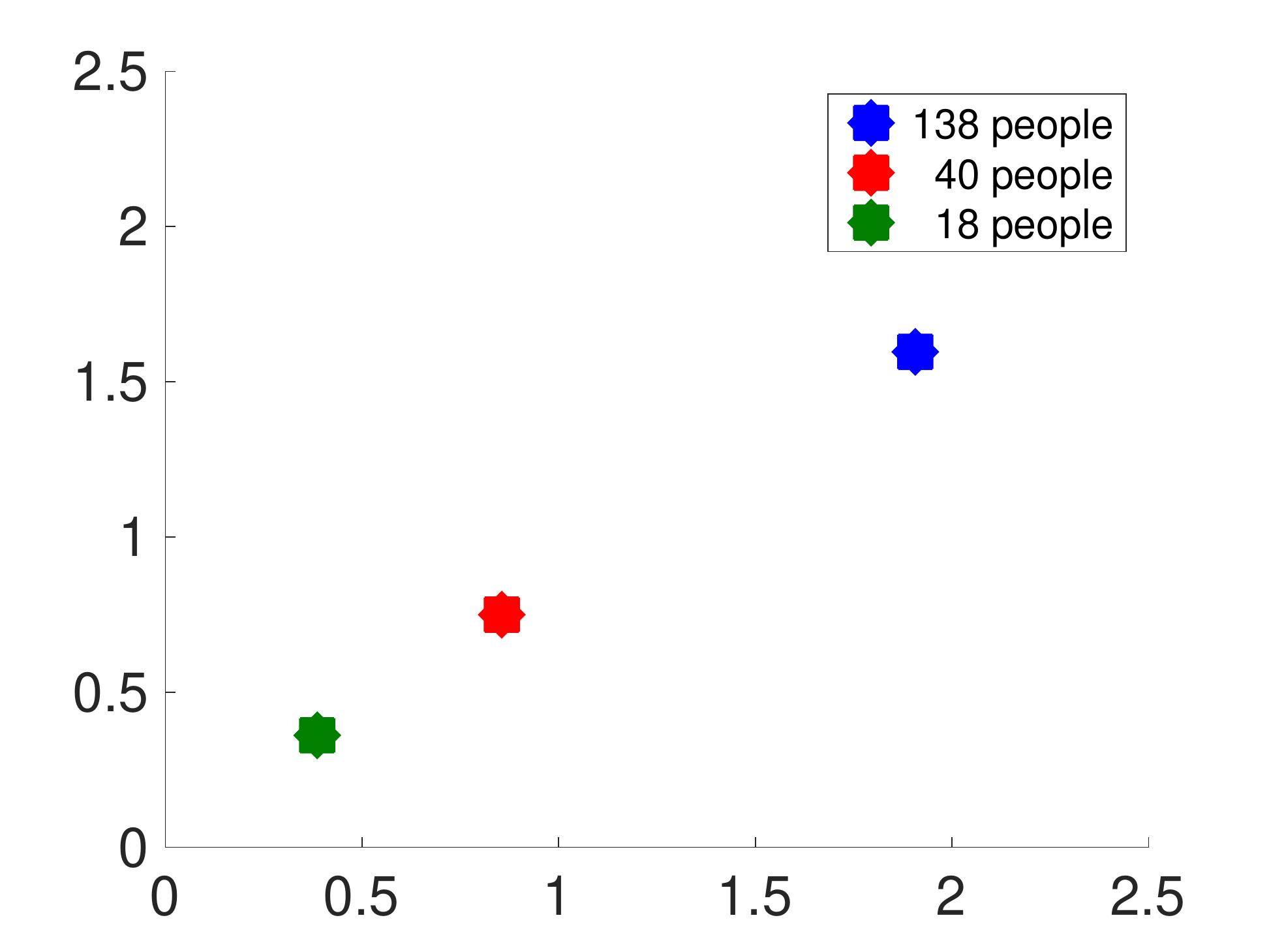}
\linethickness{1pt}
\put(75, -3){\color{black}\text{\footnotesize Exit 1}}
\put(0, 50){\color{black}\rotatebox{90}{\text{\footnotesize Exit 2}}}
\end{overpic}
}
\subfloat[Computed $F_{V}$]{
\newline
\begin{overpic}[width=0.33\textwidth,  grid=false,tics=10]{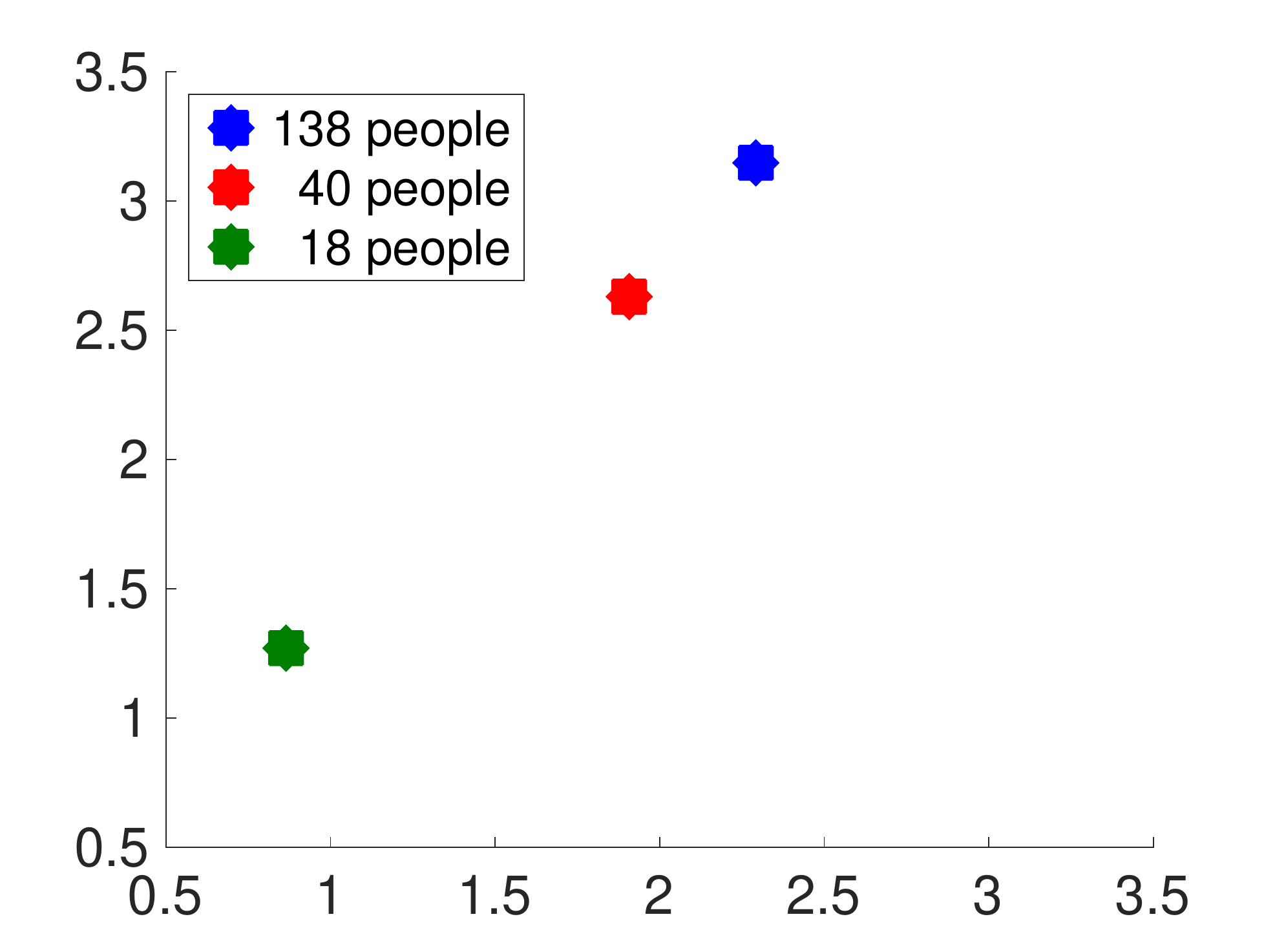}
\linethickness{1pt}
\put(75, -3){\color{black}\text{\footnotesize Exit 1}}
\put(0, 50){\color{black}\rotatebox{90}{\text{\footnotesize Exit 2}}}
\end{overpic}
} \\
\subfloat[Measured $D_{V}$]{
\newline
\begin{overpic}[width=0.33\textwidth, grid=false,tics=10]{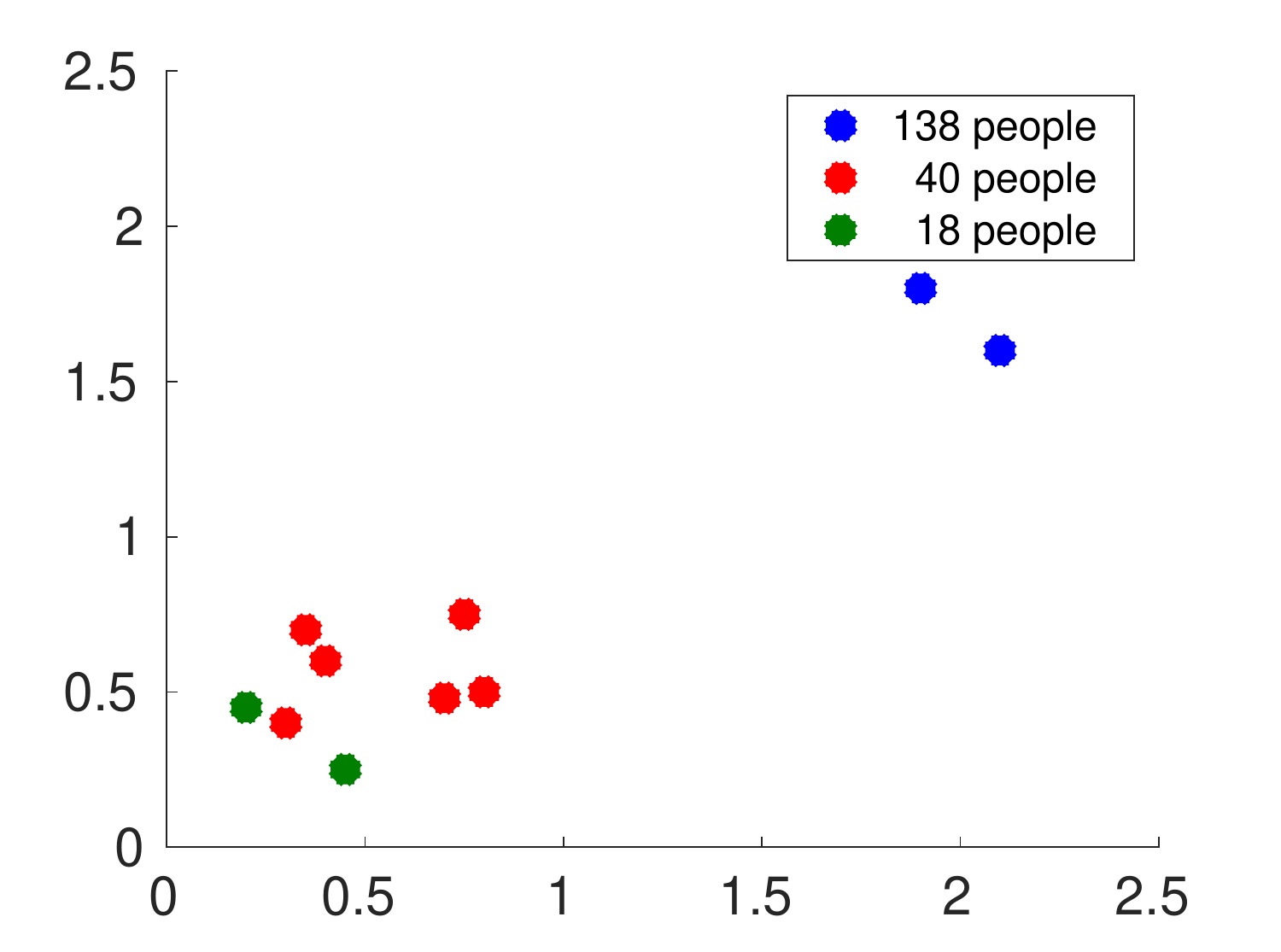}
\linethickness{1pt}
\put(75, -3){\color{black}\text{\footnotesize Exit 1}}
\put(0, 50){\color{black}\rotatebox{90}{\text{\footnotesize Exit 2}}}
\end{overpic}
}
\subfloat[Measured $F_{V}$]{
\newline
\begin{overpic}[width=0.33\textwidth, grid=false,tics=10]{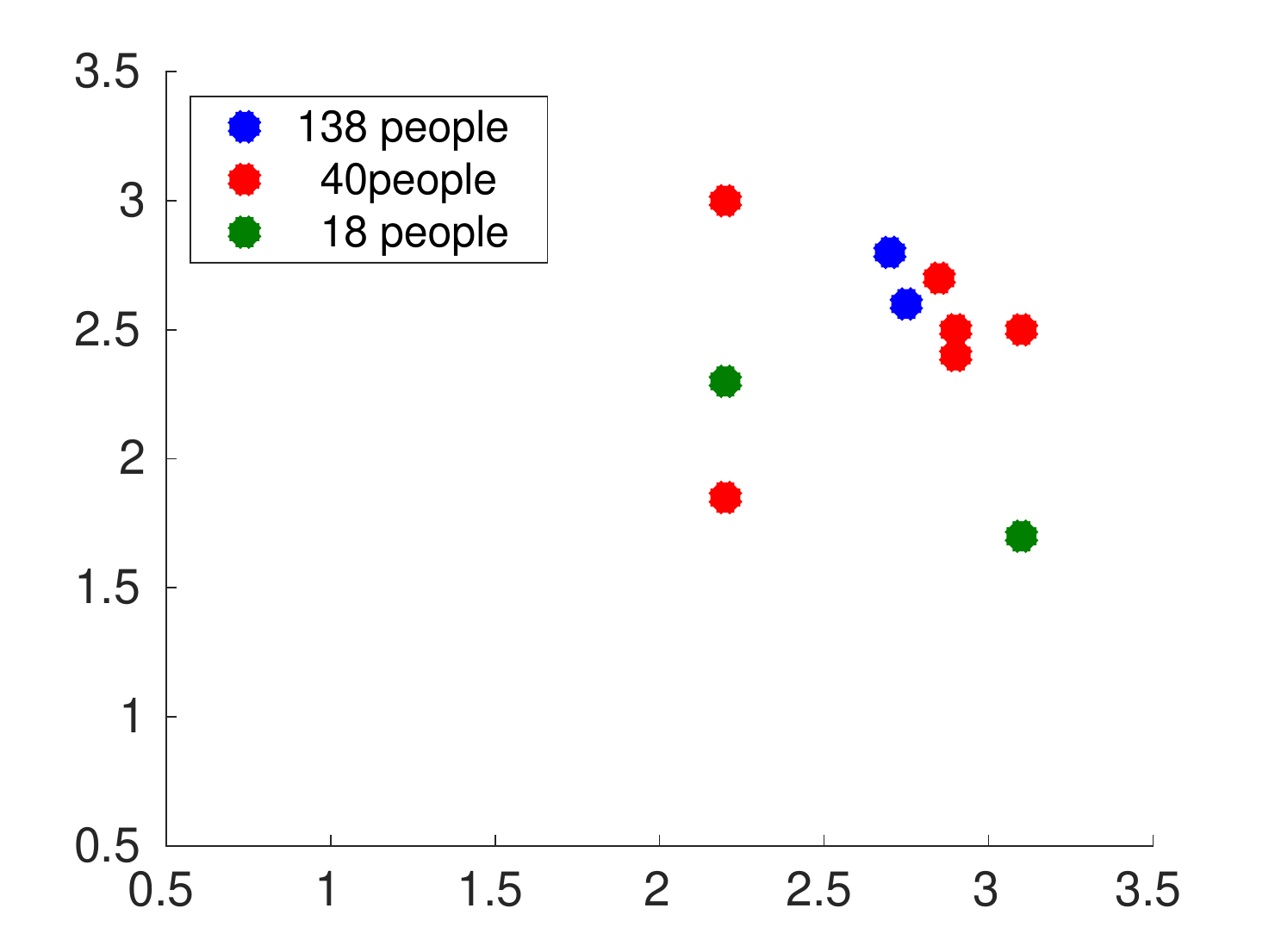}
\linethickness{1pt}
\put(75, -3){\color{black}\text{\footnotesize Exit 1}}
\put(0, 50){\color{black}\rotatebox{90}{\text{\footnotesize Exit 2}}}
\end{overpic}
}
\caption{Computed (A) mean density $D_{V}$ and (B) mean flow rate $F_{V}$ 
as defined in \eqref{Voronoimethod}, and 
measured (C) mean density and (D) mean flow rate from \cite{Kemloh}.
} 
\label{Voronoi}
\end{figure}

To compute the mean density $D_{V}$ and mean flow rate $F_{V}$ we use 
the Voronoi method \cite{Steffen2009}:
\begin{equation} \label{Voronoimethod}
D_{V}(t) = \cfrac{\sum_{\x \in \omega} \rho(t,\x)}{|\,\omega\,|},
\quad  F_{V}(t)=D_{V} V_{V} E,
\end{equation}
where $V_{V}$ is the mean velocity modulus over the area $\omega$ and $E$ is the exit width.
We choose  two 4 m$^2$ areas in front of the exits as $\omega$.
The mean density and mean flow rate computed from our simulations are shown 
in Figure~$\ref{Voronoi}$ (A) and (B). Figure~$\ref{Voronoi}$ (C) and (D)
report the measured mean density and mean flow rate form the experiments 
in \cite{Kemloh}. We see very good agreement between 
computed and measured quantities for the 40 and 138 people cases,
while there is no good agreement for the 18 people case. 
This is to be expected since the kinetic approach is not meant to simulate the dynamics 
of a small number of pedestrians. 
{With this test we concludes the validation of our software.}

\begin{figure}[h!]
\begin{center}
\subfloat[Different velocity moduli]{
\centering
\begin{overpic}[width=0.4\textwidth,grid=false,tics=10]{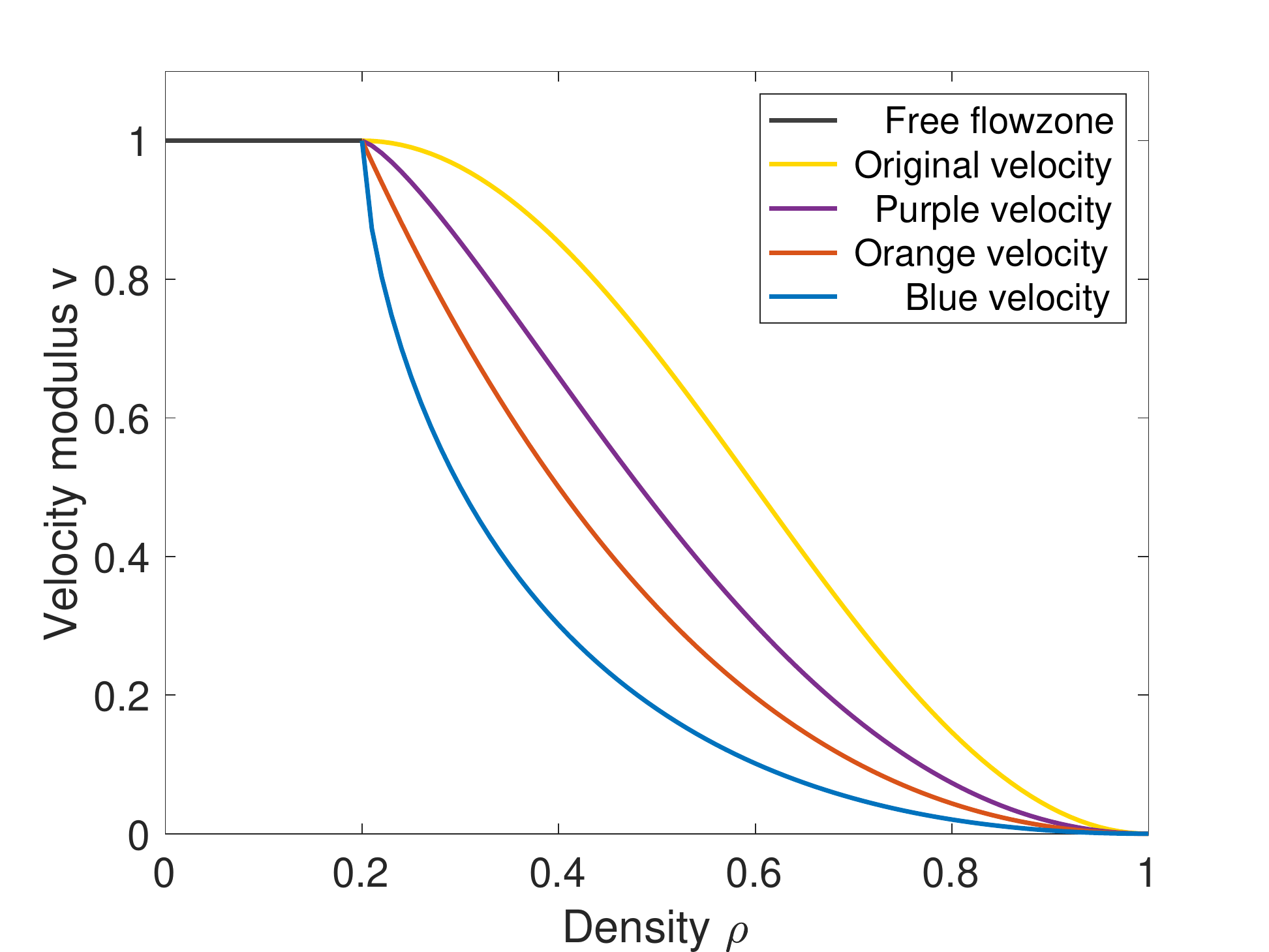}
\end{overpic}
}
\subfloat[Number of pedestrians]{
\begin{overpic}[width=0.4\textwidth,grid=false,tics=10]{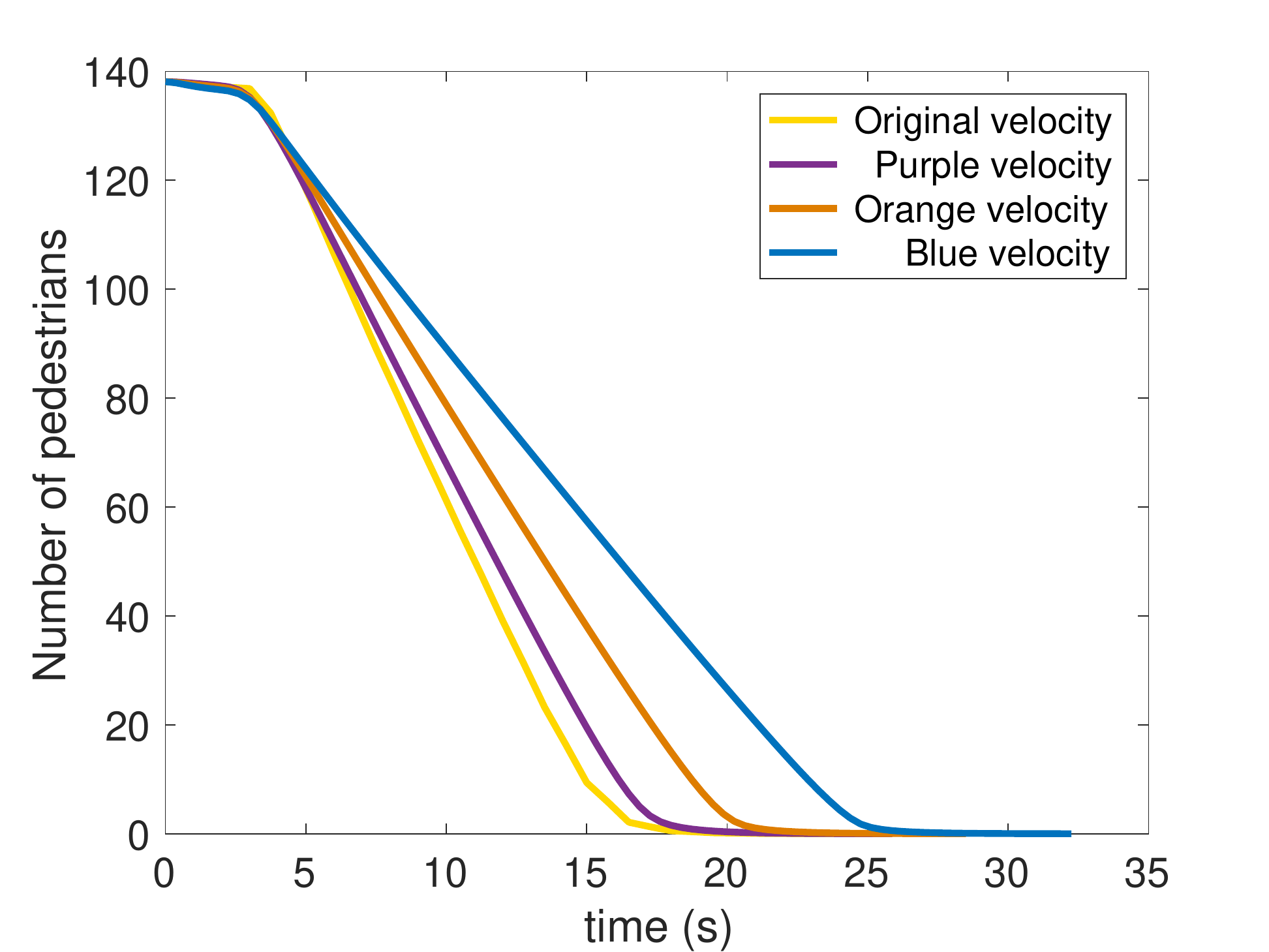}
\end{overpic}
}
\end{center}
\caption{(A) Different velocity moduli under consideration and (B) corresponding number of pedestrians
in the room versus time for the 138 pedestrian case.}
\label{different_velocity}
\end{figure}

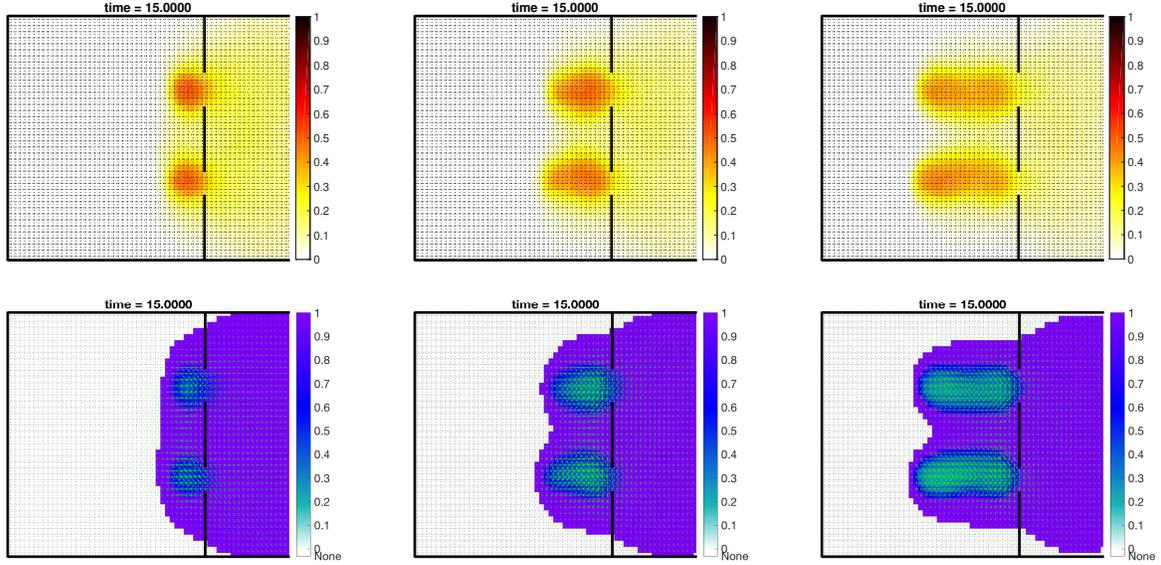
\begin{figure}[h!]
\centering
\begin{overpic}[width=0.32\textwidth,grid=false,tics=10]{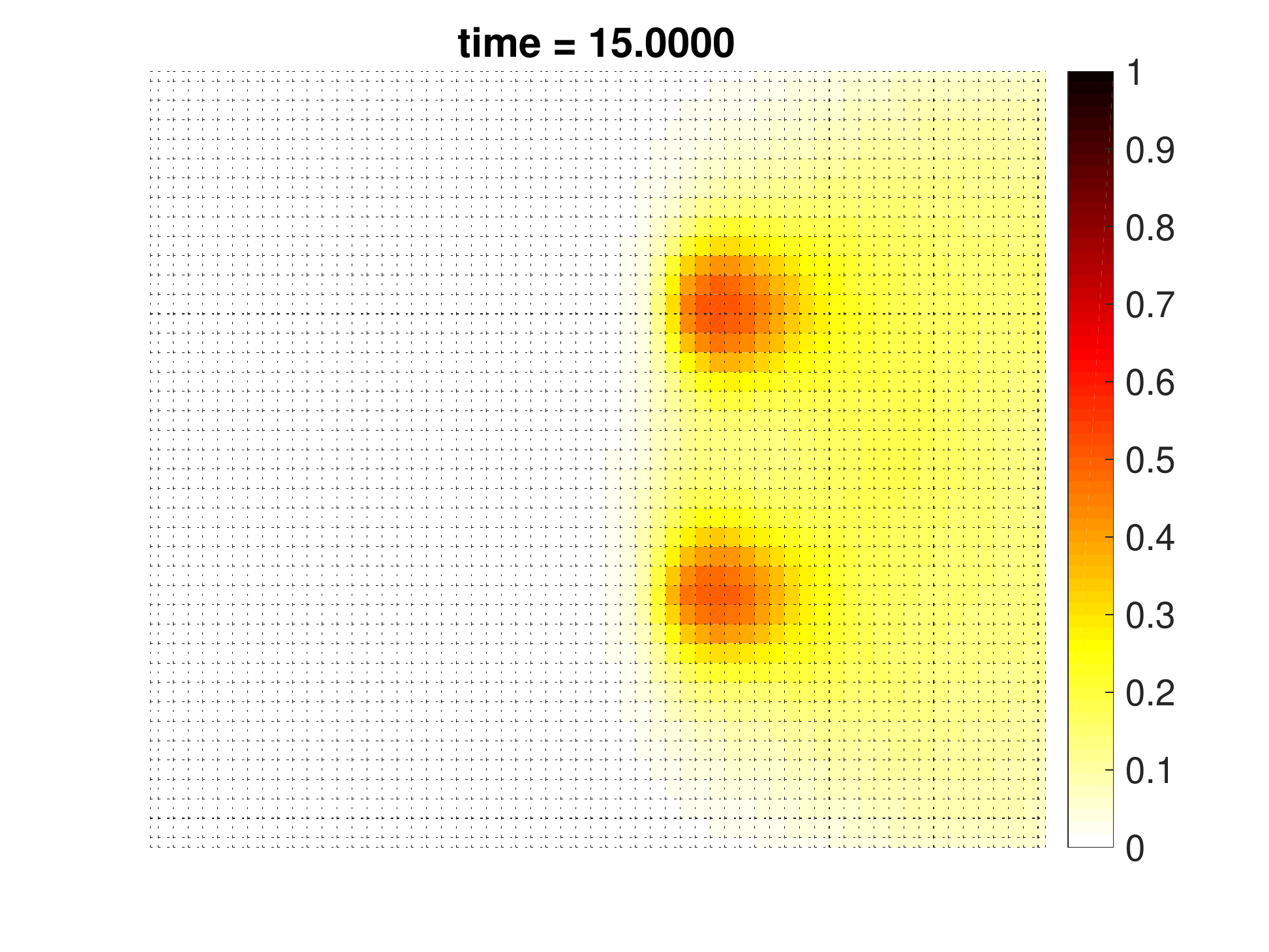}
\linethickness{1pt}
\put(11.7, 69.5){\line(1,0){71}}
\put(11.7, 7.8){\line(0,1){61.9}}
\put(11.7, 8){\line(1,0){71}}
\put(61.2, 7.8){\line(0,1){16.7}}
\put(61.2, 30.2){\line(0,1){16.5}}
\put(61.2, 55.2){\line(0,1){14.5}}
\end{overpic}
\begin{overpic}[width=0.32\textwidth,grid=false,tics=10]{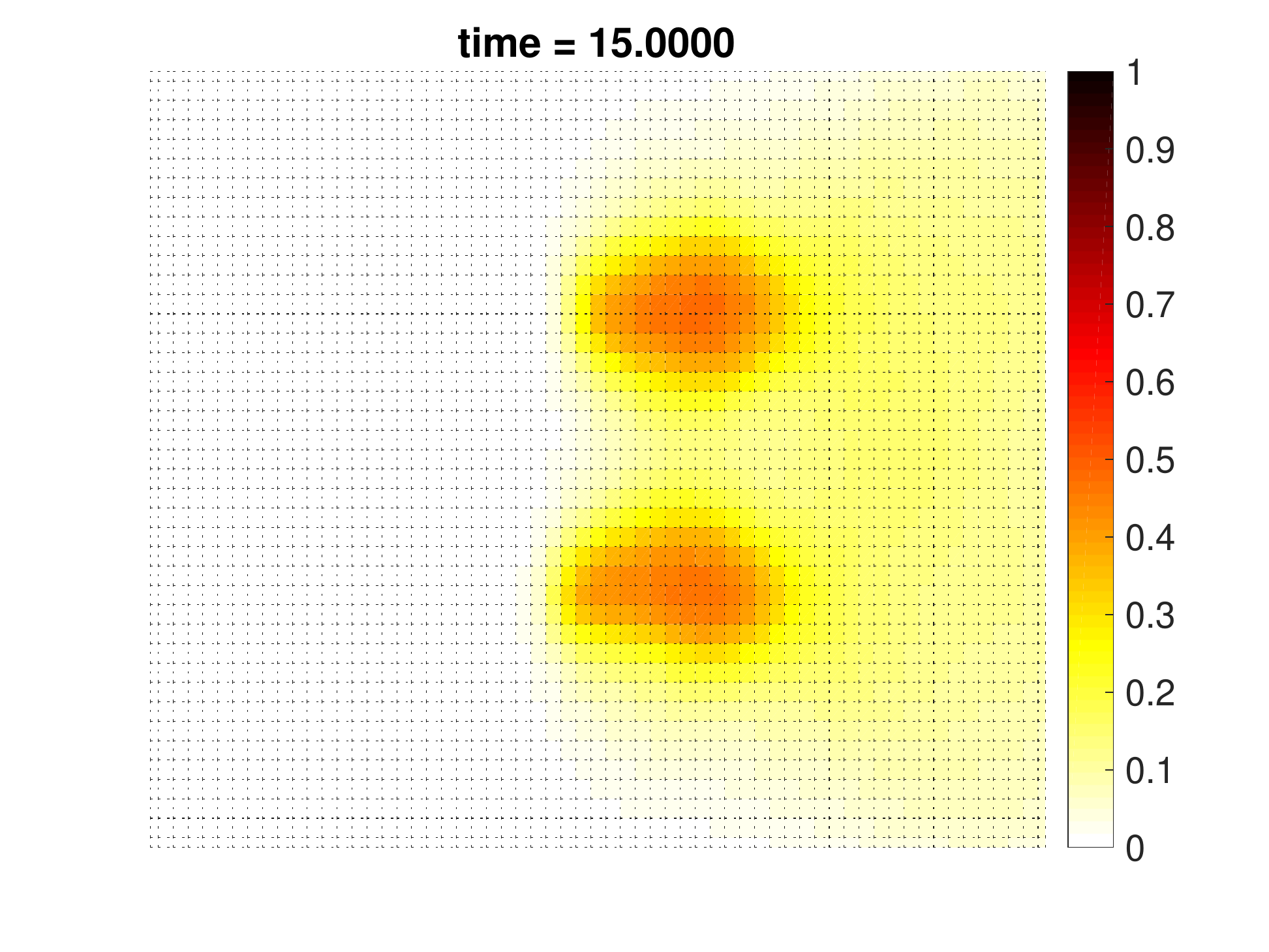}
\linethickness{1pt}
\put(11.7, 69.5){\line(1,0){71}}
\put(11.7, 7.8){\line(0,1){61.9}}
\put(11.7, 8){\line(1,0){71}}
\put(61.2, 7.8){\line(0,1){16.7}}
\put(61.2, 30.2){\line(0,1){16.5}}
\put(61.2, 55.2){\line(0,1){14.5}}
\end{overpic}
\begin{overpic}[width=0.32\textwidth,grid=false,tics=10]{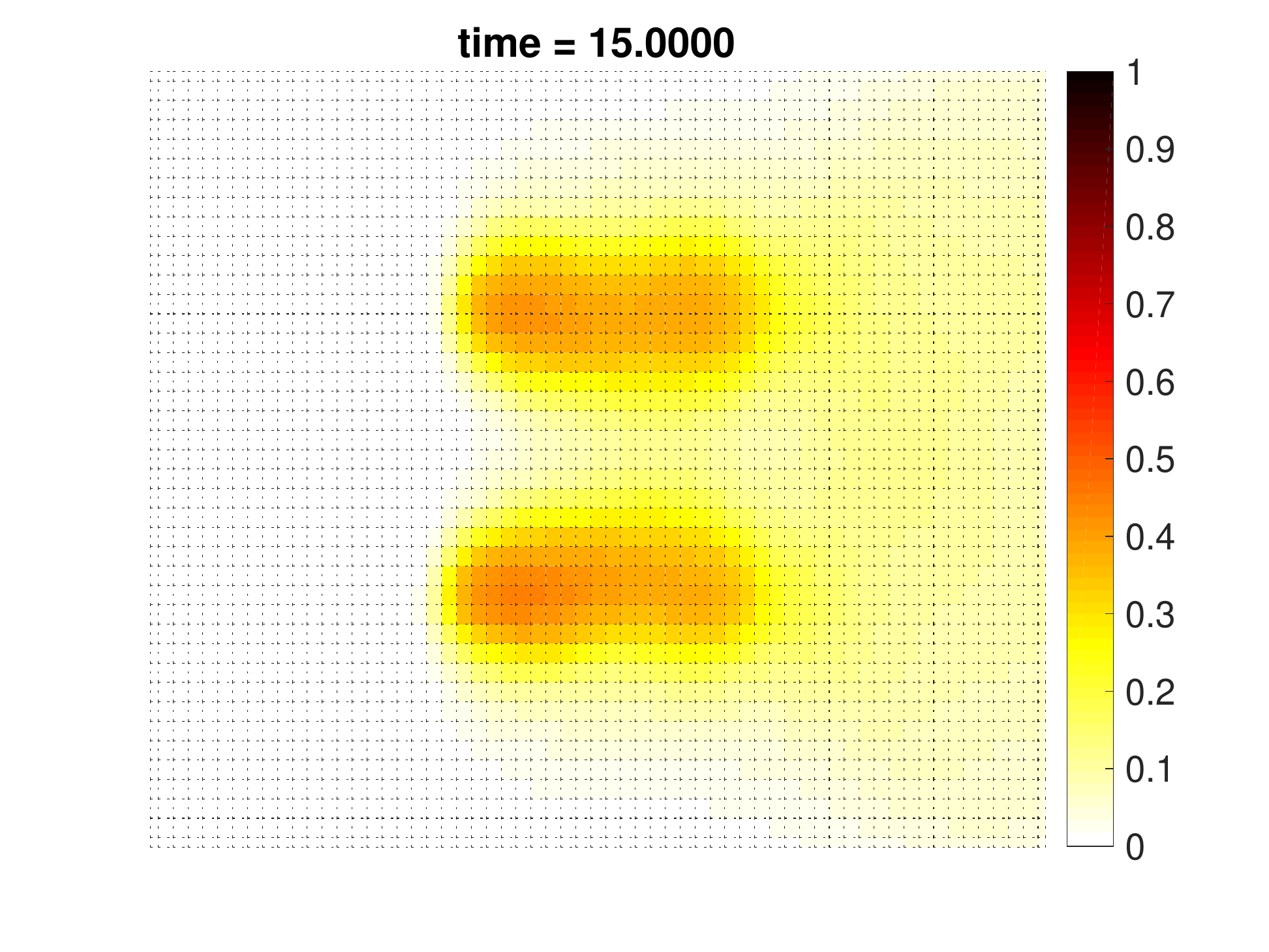}
\linethickness{1pt}
\put(11.7, 69.5){\line(1,0){71}}
\put(11.7, 7.8){\line(0,1){61.9}}
\put(11.7, 8){\line(1,0){71}}
\put(61.2, 7.8){\line(0,1){16.7}}
\put(61.2, 30.2){\line(0,1){16.5}}
\put(61.2, 55.2){\line(0,1){14.5}}
\end{overpic}
\begin{overpic}[width=0.32\textwidth,grid=false,tics=10]{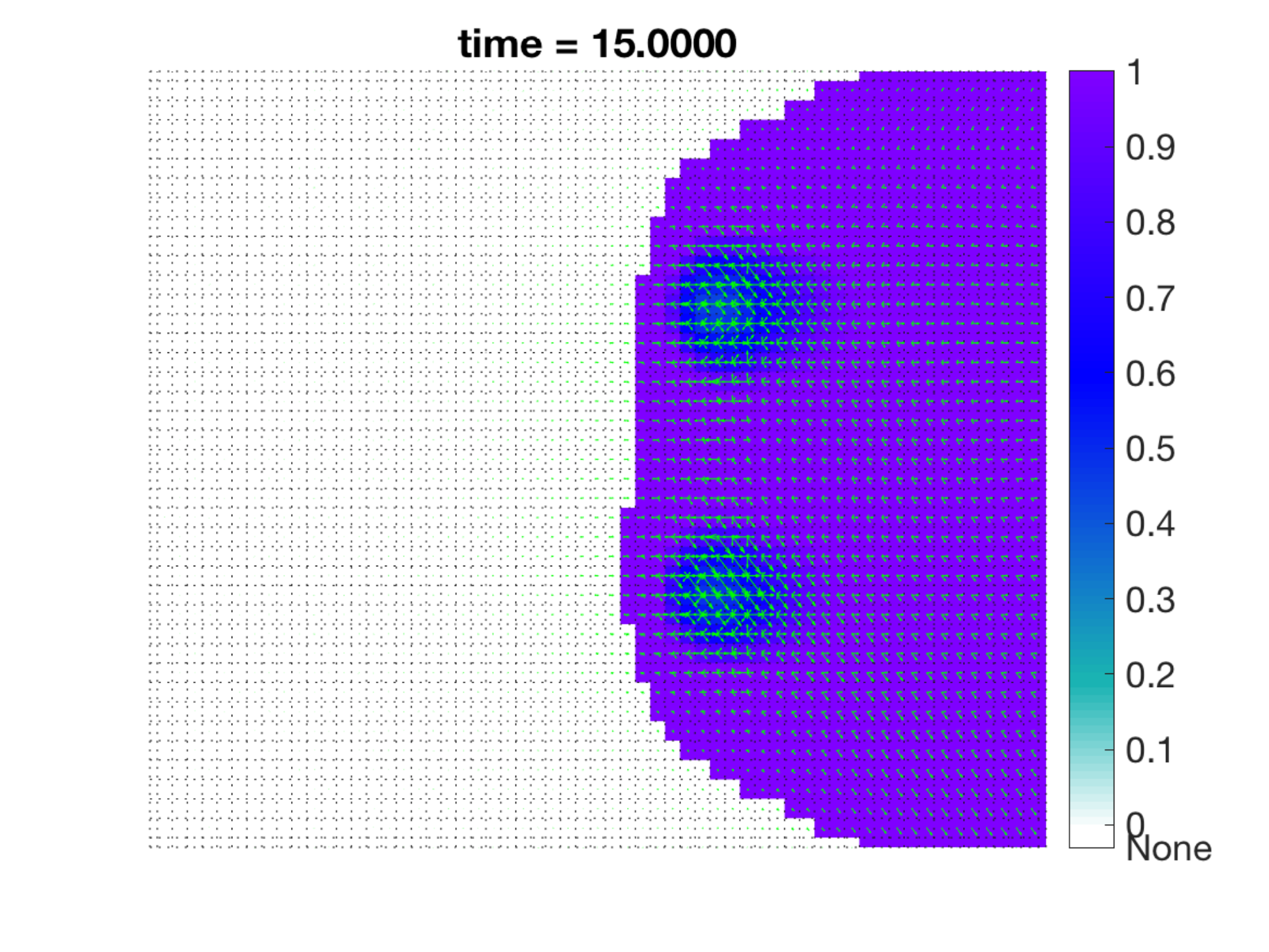}
\linethickness{1pt}
\put(11.7, 69.5){\line(1,0){71}}
\put(11.7, 7.8){\line(0,1){61.9}}
\put(11.7, 8){\line(1,0){71}}
\put(61.4, 7.8){\line(0,1){16.7}}
\put(61.4, 30.4){\line(0,1){16.5}}
\put(61.4, 55.2){\line(0,1){14.5}}
\end{overpic}
\begin{overpic}[width=0.32\textwidth,grid=false,tics=10]{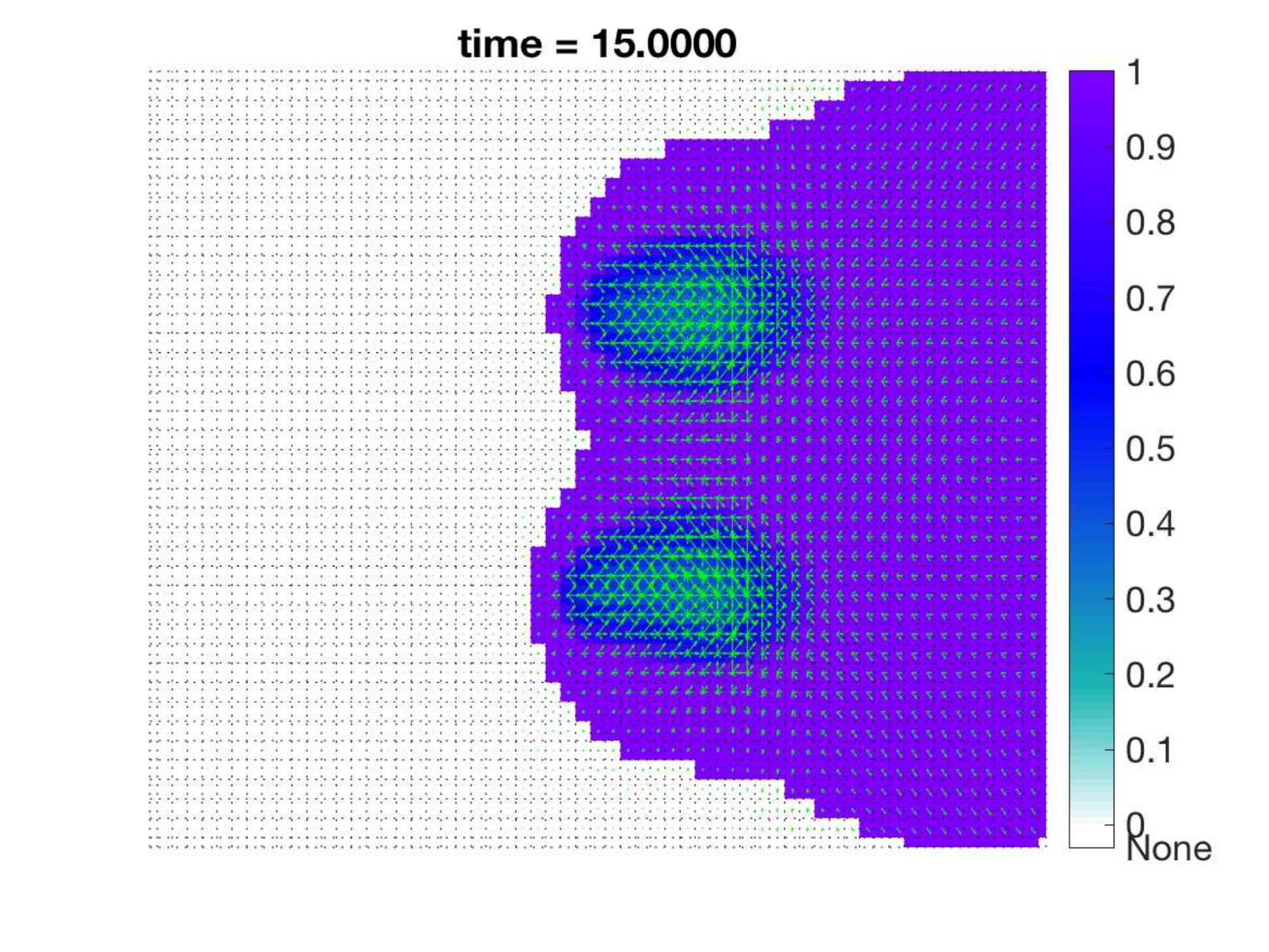}
\linethickness{1pt}
\put(11.7, 69.5){\line(1,0){71}}
\put(11.7, 7.8){\line(0,1){61.9}}
\put(11.7, 8){\line(1,0){71}}
\put(61.4, 7.8){\line(0,1){16.7}}
\put(61.4, 30.4){\line(0,1){16.5}}
\put(61.4, 55.2){\line(0,1){14.5}}
\end{overpic}
\begin{overpic}[width=0.32\textwidth,grid=false,tics=10]{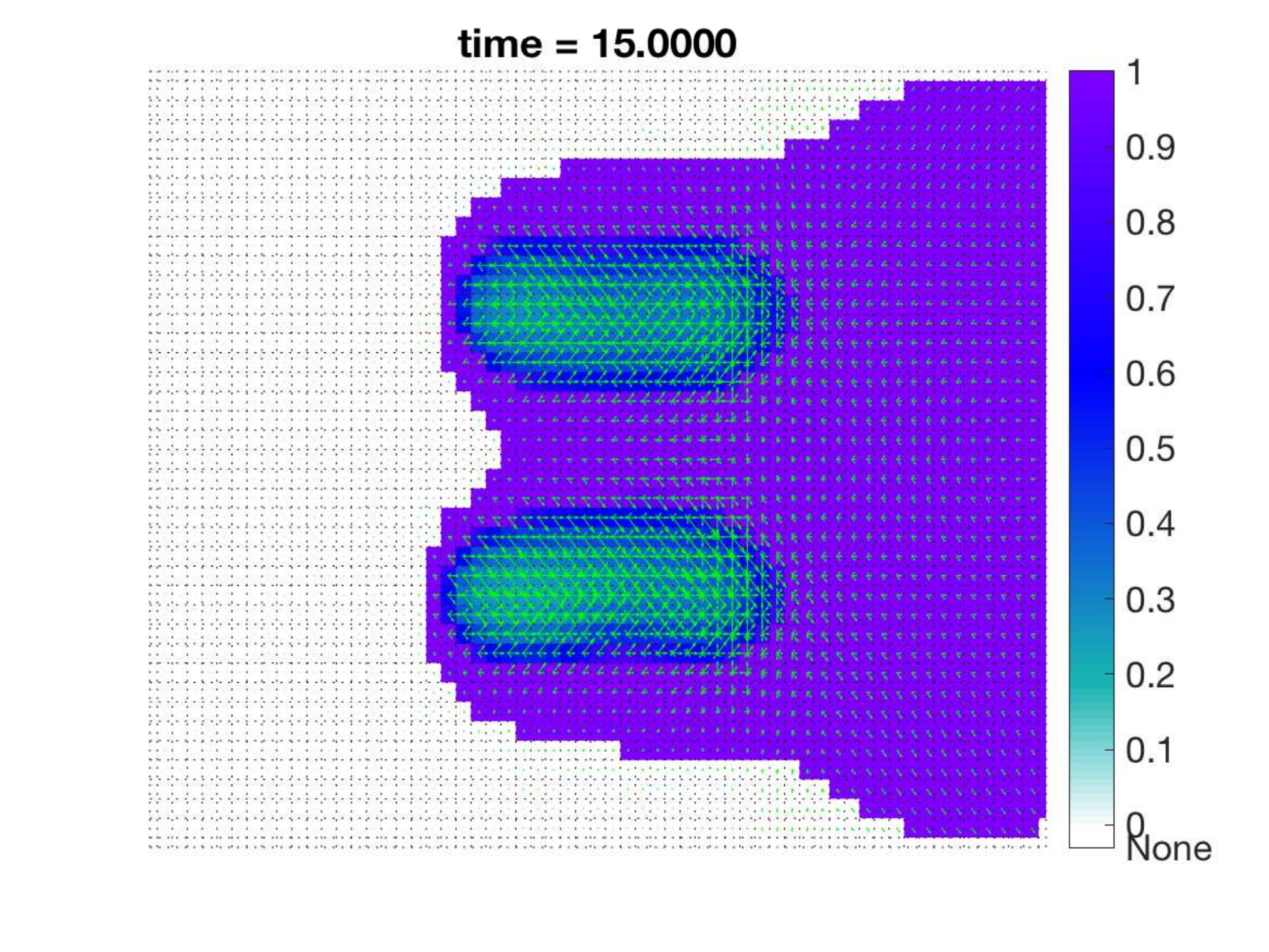}
\linethickness{1pt}
\put(11.7, 69.5){\line(1,0){71}}
\put(11.7, 7.8){\line(0,1){61.9}}
\put(11.7, 8){\line(1,0){71}}
\put(61.4, 7.8){\line(0,1){16.7}}
\put(61.4, 30.4){\line(0,1){16.5}}
\put(61.4, 55.2){\line(0,1){14.5}}
\end{overpic}
\caption{Density (top) and velocity magnitude with selected velocity vectors (bottom) 
for the evacuation process of 138 pedestrians with the purple (left), orange (middle), and blue (right) velocity moduli 
at time $t=15$ s.} 
\label{Diff_Velo_Moduli}
\end{figure}

So far, we have only considered velocity modulus $v$ as defined (\ref{eq:v}), which 
uses a cubic polynomial outside the free flow regime. 
For the two exit test, we take into consideration other possible choices:
\begin{align}
v_{purple}(\rho)&=(1+\cos((\rho-0.2)^{\frac{2}{3}}\pi/0.8^{\frac{2}{3}}))/2, \cl
v_{orange}(\rho)&=(1+\cos((\rho-0.2)^{\frac{1}{2}}\pi/0.8^{\frac{1}{2}}))/2, \cl
v_{blue}(\rho)&=(1+\cos((\rho-0.2)^{\frac{1}{3}}\pi/0.8^{\frac{1}{3}}))/2. \label{eq:v_blue}
\end{align}
See Figure~\ref{different_velocity} (A). We repeat the test with 138 pedestrians for all the above velocity moduli.
The corresponding number of pedestrians in the room versus times is shown in Figure~\ref{different_velocity} (B). 
For a given density $\overline{\rho}$ in the slowdown zone, we have
$v(\overline{\rho}) > v_{purple}(\overline{\rho})  > v_{orange}(\overline{\rho})  > v_{blue}(\overline{\rho})$.
Thus, it is not surprising the the total evacuation time increases 
as we pass from the original definition of velocity modulus in (\ref{eq:v}) to
$v_{blue}$ through  $v_{purple}$ and $v_{orange}$.
To further illustrate the difference in the evacuation process when the above 
velocity moduli are used, in Figure~\ref{Diff_Velo_Moduli} we display 
the density and velocity modulus (with selected velocity vectors) at time $t  = 15$ 
when using the three velocity moduli. 

\begin{figure}[h!] 
\begin{center}
\subfloat[Exit 2 width and width ratios]{
\begin{overpic}[height=0.35\textwidth, grid=false,tics=10]{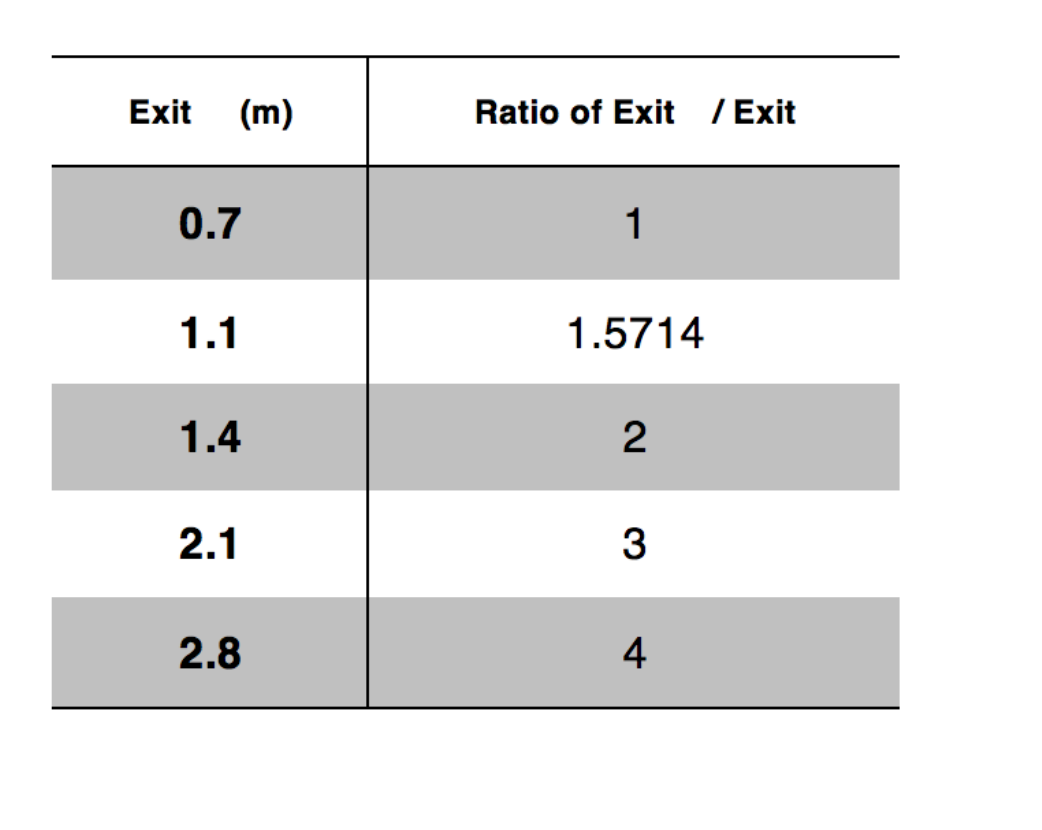}
\put(20, 67.7){\text{\small {2}}}
\put(65.5, 67.7){\text{\small {2}}}
\put(77.4, 67.7){\text{\small {1}}}
\end{overpic}
}
\subfloat[Evacuation time vs width ratio]{
\begin{overpic}[height=0.35\textwidth,grid=false]{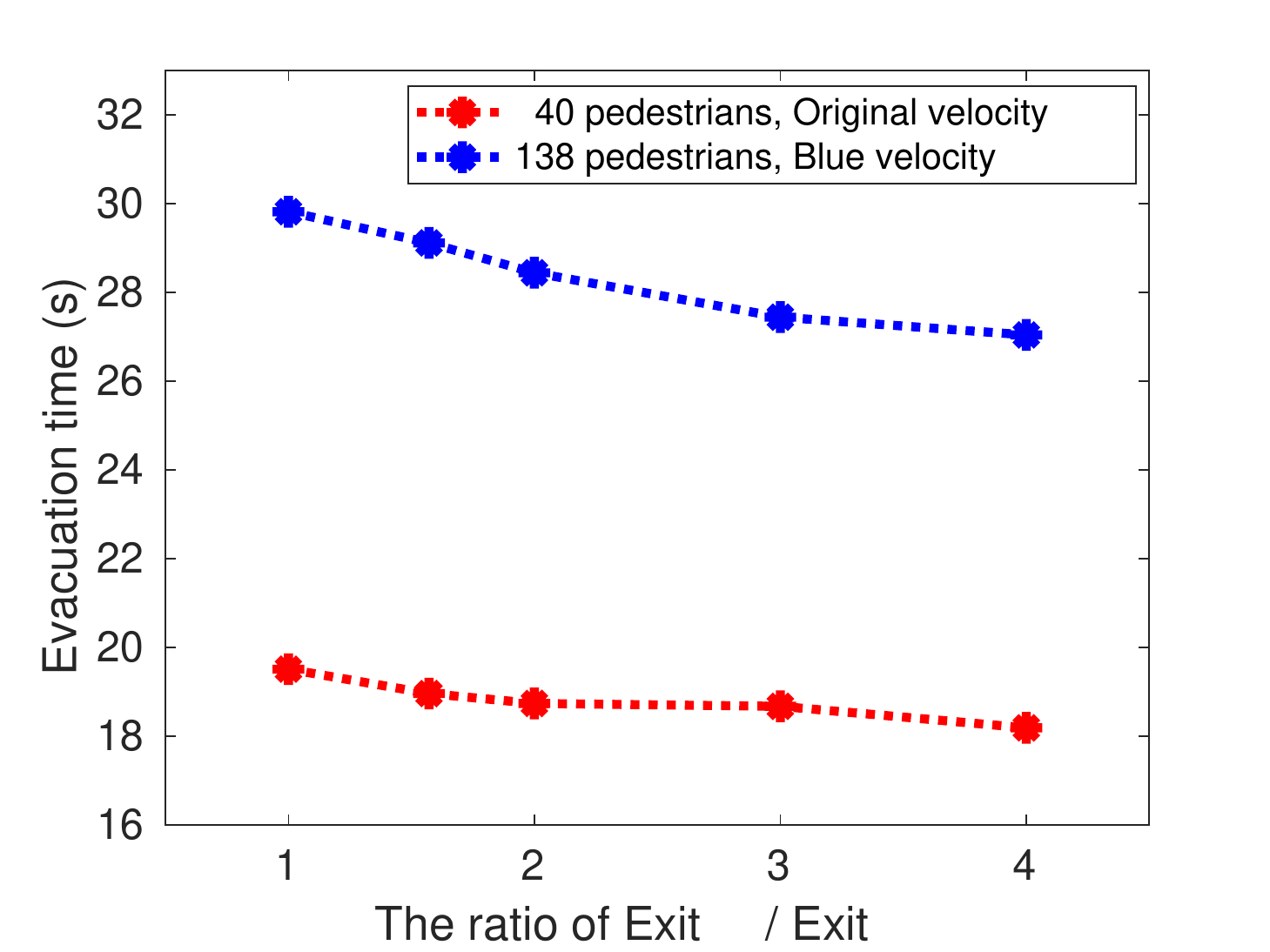}
\put(57, 1){\text{\small {2}}}
\put(70, 1){\text{\small {1}}}
\end{overpic}
}
\end{center}
\caption{(A) Widths of exit 2 under consideration and corresponding  width ratios and (B) evacuation 
time versus width ratios for two different scenarios.}
\label{width_ratios}
\end{figure}
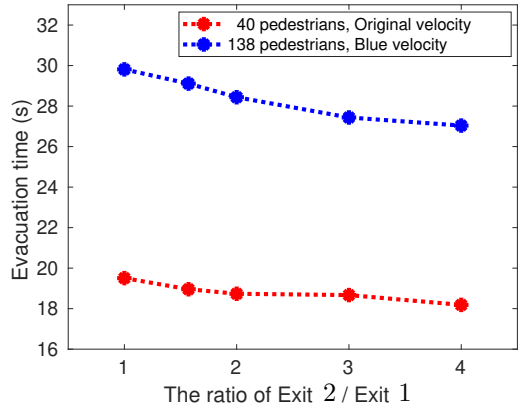

Next, we investigate the relationship between evacuation time and the width ratio of the two exits. 
We fix the size and position of exit {1}, and the position of the center of exit {2}, 
while the size of exit {2} varies. See Figure~\ref{width_ratios} (A) for the widths of exit 2 under consideration
and corresponding  width ratios.
We consider two scenarios: the group of 40 people with velocity modulus (\ref{eq:v}) 
and the group of 138 people with the blue velocity modulus \eqref{eq:v_blue}. 
In both cases, all other model and discretization parameters are 
set like for the results reported in Figure~\ref{Voronoi}.
Figure~\ref{width_ratios} (B) shows the total evacuation time versus the 
exit width ratio for both scenarios. As expected, when the ratio of exit ${2}$ width/exit ${1}$ width increases 
the total evacuation time decreases. 

\subsection{Evacuation from a room with obstacles}\label{sec:res_obstacle}

In the model introduced in Section~\ref{modelinginteractions},
parameter $\alpha$ represents the quality of the environment, which influences
also the maximal dimensionless velocity modulus (\ref{eq:v})
a pedestrian can reach.
In theory, parameter $\alpha=0$  forces pedestrians to stop, 
while the value $\alpha=1$ contributes to keep the maximal velocity modulus.
However, in practice this parameter alone is not sufficient to model obstacles
within the domain. To effectively model the obstacles, we use the strategy reported
in Section~\ref{sec:obstacles}.

\begin{figure}[h!]
\begin{center}
\subfloat[Configuration 1]{
\begin{overpic}[width=0.33\textwidth,grid=false,tics=10]{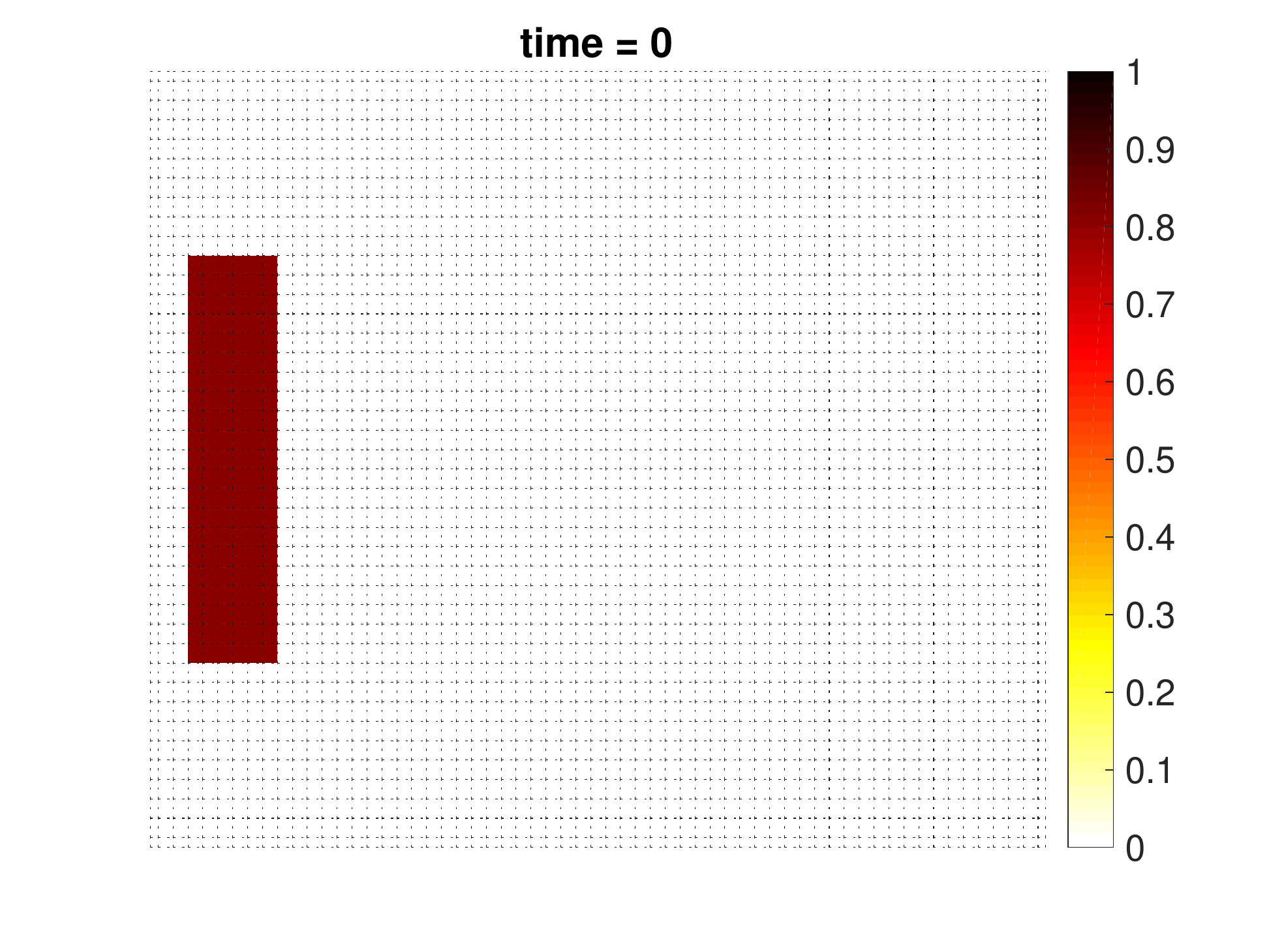}
\linethickness{1pt}
\put(20.3, 38.8){\vector(1,0){6}}
\put(11.7, 69.5){\line(1,0){71}}
\put(11.7, 7.8){\line(0,1){61.9}}
\put(11.7, 8){\line(1,0){71}}
\put(61.2, 7.8){\line(0,1){23}}
\put(61.2, 46.7){\line(0,1){23}}

\put(42.3, 46.5){\line(1,0){9.6}}
\put(42.3, 31.1){\line(1,0){9.6}}
\put(42.5, 31.1){\line(0,1){15.6}}
\put(51.8, 30.9){\line(0,1){15.8}}

\end{overpic}
}
\subfloat[Configuration 2]{
\begin{overpic}[width=0.33\textwidth,grid=false,tics=10]{./MM_alpha_no_time_0.pdf}
\linethickness{1pt}
\put(20.3, 38.8){\vector(1,0){6}}
\put(11.7, 69.5){\line(1,0){71}}
\put(11.7, 7.8){\line(0,1){61.9}}
\put(11.7, 8){\line(1,0){71}}
\put(61.2, 7.8){\line(0,1){23}}
\put(61.2, 46.7){\line(0,1){23}}

\put(42.3, 46.5){\line(1,0){9.6}}
\put(42.3, 61.9){\line(1,0){9.6}}
\put(42.5, 46.5){\line(0,1){15.6}}
\put(51.8, 46.3){\line(0,1){15.8}}

\put(42.3, 14.3){\line(1,0){9.6}}
\put(42.3, 29.7){\line(1,0){9.6}}
\put(42.5, 14.3){\line(0,1){15.6}}
\put(51.8, 14.1){\line(0,1){15.8}}
\end{overpic} 
}

\caption{Computational domain with effective area for (A) an obstacle placed in the middle of the room, towards the exit, and (B) two obstacles placed symmetrically with respect to the exit. }
\label{obstacle_position}
\end{center}
\end{figure}
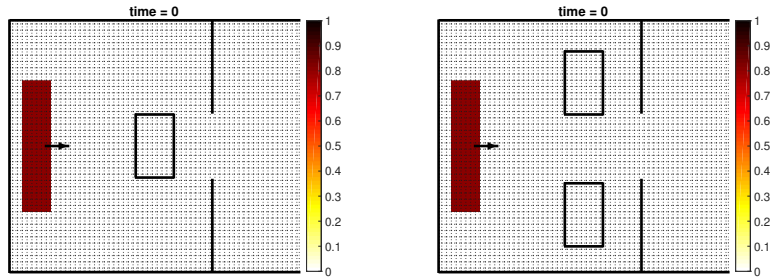

We consider a square room of side $10$ m with a $2.6$ m wide exit located on the right wall
as in Section~\ref{sec:one_exit}, with the following obstacle configuration:
\begin{enumerate}
\item One obstacle close to the exit, i.e. in the middle of the right wall, with effective area
depicted in Figure~\ref{obstacle_position} (A).
\item Two obstacles close to the right wall, place symmetrically with respect to the exit. See
Figure~\ref{obstacle_position} (B) for the effective area of both obstacles. 
\end{enumerate}
Pedestrians are initially distributed in a rectangular region with constant density $\rho=0.80$, 
as shown in Figure~\ref{obstacle_position}. The total number of of pedestrians is 44. 
The initial direction is  $\theta_1$.

\begin{figure}[h!]
\centering
\begin{overpic}[width=0.32\textwidth,grid=false,tics=10]{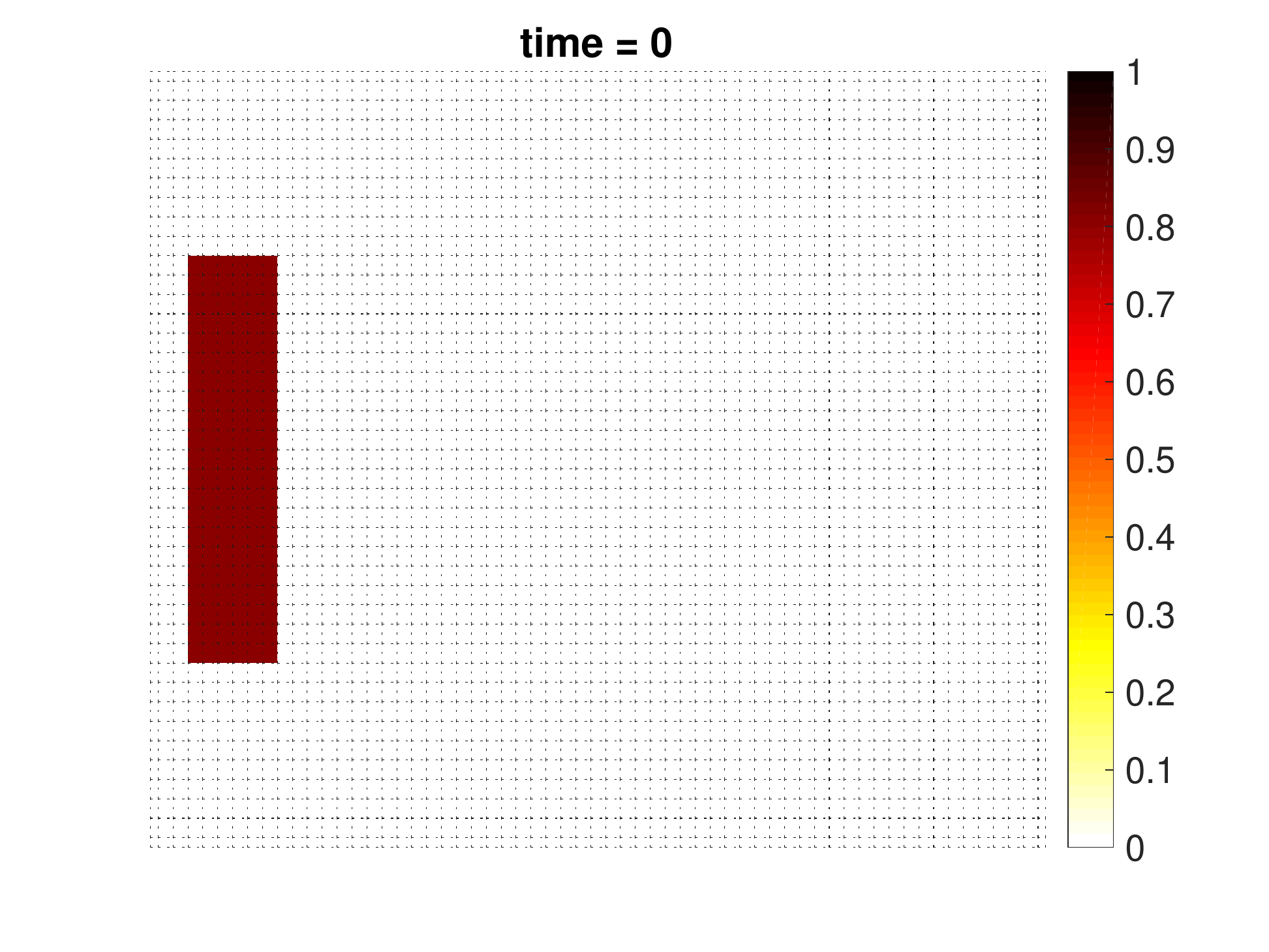}
\linethickness{1pt}
\put(20.3, 38.8){\vector(1,0){6}}
\put(11.7, 69.5){\line(1,0){71}}
\put(11.7, 7.8){\line(0,1){61.9}}
\put(11.7, 8){\line(1,0){71}}
\put(61.2, 7.8){\line(0,1){23}}
\put(61.2, 46.7){\line(0,1){23}}

\put(42.3, 46.5){\line(1,0){9.6}}
\put(42.3, 31.1){\line(1,0){9.6}}
\put(42.5, 31.1){\line(0,1){15.6}}
\put(51.8, 30.9){\line(0,1){15.8}}

\put(42.5, 42){\line(1,0){6}}
\put(42.5, 35.7){\line(1,0){6}}
\put(42.7, 35.9){\line(0,1){6}}
\put(48.6, 35.5){\line(0,1){6.7}}

\end{overpic} 
\begin{overpic}[width=0.32\textwidth,grid=false,tics=10]{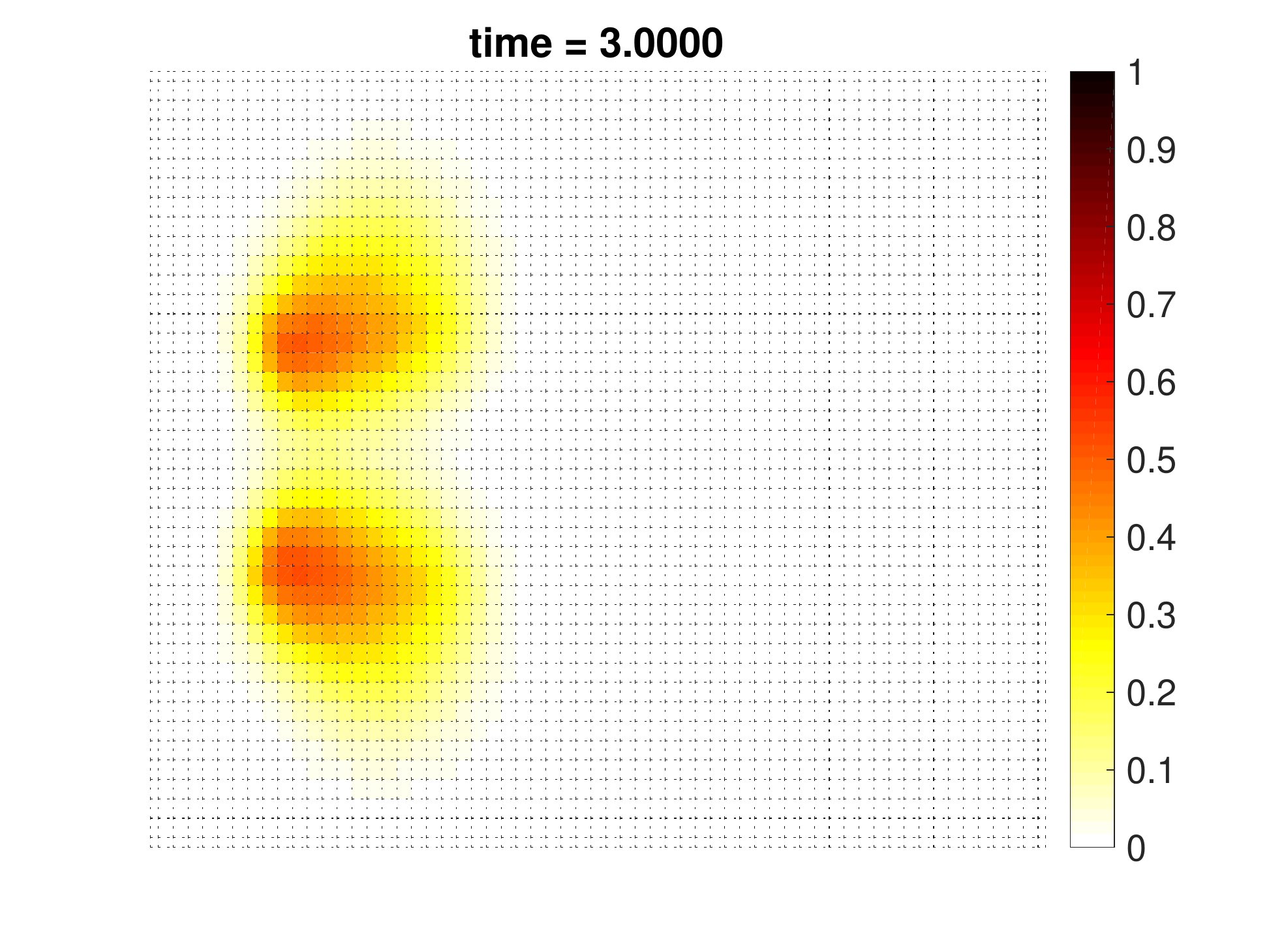}
\linethickness{1pt}
\put(11.7, 69.5){\line(1,0){71}}
\put(11.7, 7.8){\line(0,1){61.9}}
\put(11.7, 8){\line(1,0){71}}
\put(61.2, 7.8){\line(0,1){23}}
\put(61.2, 46.7){\line(0,1){23}}

\put(42.3, 46.5){\line(1,0){9.6}}
\put(42.3, 31.1){\line(1,0){9.6}}
\put(42.5, 31.1){\line(0,1){15.6}}
\put(51.8, 30.9){\line(0,1){15.8}}
\put(42.5, 42){\line(1,0){6}}
\put(42.5, 35.7){\line(1,0){6}}
\put(42.7, 35.9){\line(0,1){6}}
\put(48.6, 35.5){\line(0,1){6.7}}

\end{overpic} 
\begin{overpic}[width=0.32\textwidth,grid=false,tics=10]{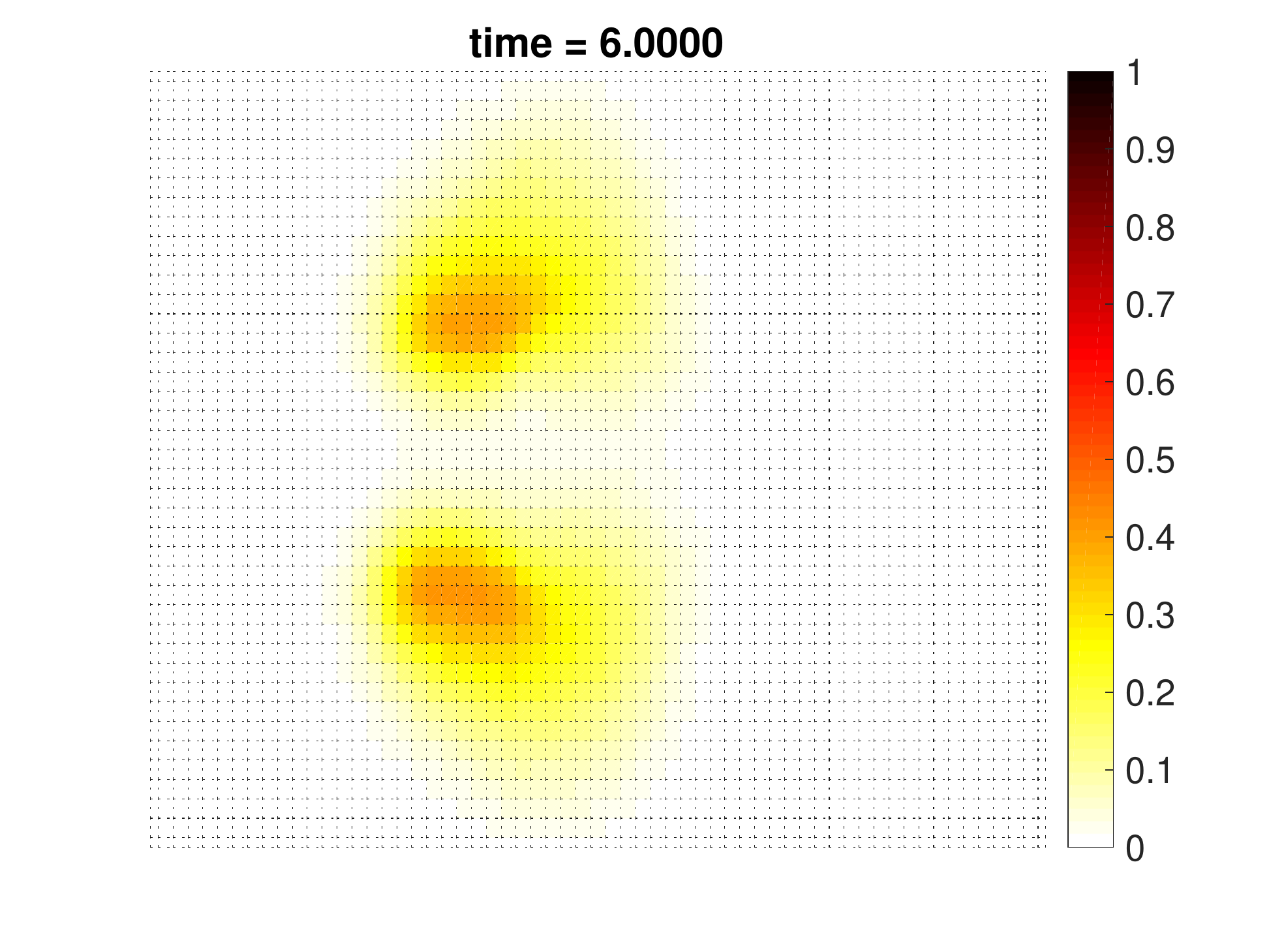}
\linethickness{1pt}
\put(11.7, 69.5){\line(1,0){71}}
\put(11.7, 7.8){\line(0,1){61.9}}
\put(11.7, 8){\line(1,0){71}}
\put(61.2, 7.8){\line(0,1){23}}
\put(61.2, 46.7){\line(0,1){23}}

\put(42.3, 46.5){\line(1,0){9.6}}
\put(42.3, 31.1){\line(1,0){9.6}}
\put(42.5, 31.1){\line(0,1){15.6}}
\put(51.8, 30.9){\line(0,1){15.8}}
\put(42.5, 42){\line(1,0){6}}
\put(42.5, 35.7){\line(1,0){6}}
\put(42.7, 35.9){\line(0,1){6}}
\put(48.6, 35.5){\line(0,1){6.7}}

\end{overpic}
\begin{overpic}[width=0.32\textwidth,grid=false,tics=10]{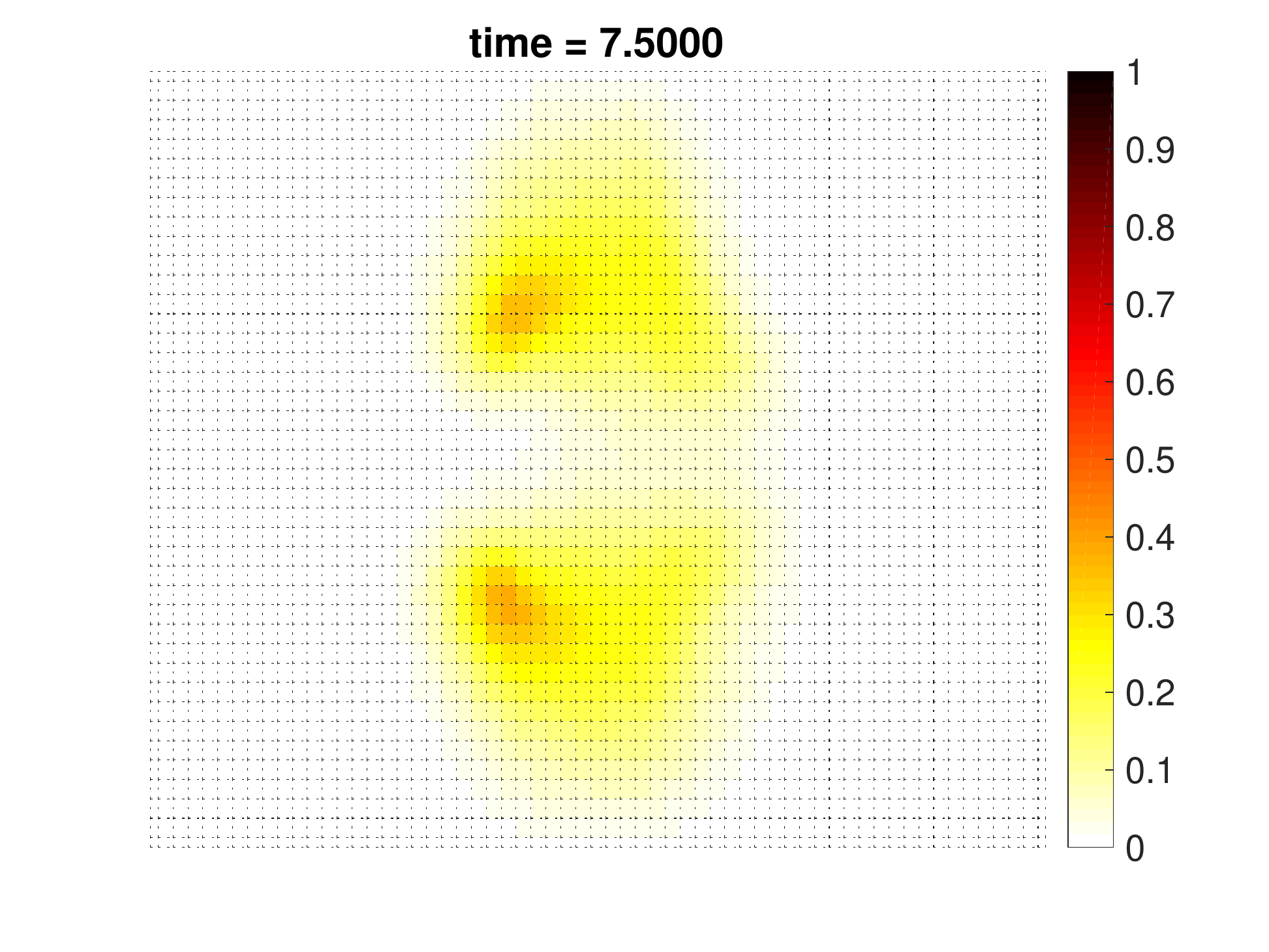}
\linethickness{1pt}
\put(11.7, 69.5){\line(1,0){71}}
\put(11.7, 7.8){\line(0,1){61.9}}
\put(11.7, 8){\line(1,0){71}}
\put(61.2, 7.8){\line(0,1){23}}
\put(61.2, 46.7){\line(0,1){23}}

\put(42.3, 46.5){\line(1,0){9.6}}
\put(42.3, 31.1){\line(1,0){9.6}}
\put(42.5, 31.1){\line(0,1){15.6}}
\put(51.8, 30.9){\line(0,1){15.8}}
\put(42.5, 42){\line(1,0){6}}
\put(42.5, 35.7){\line(1,0){6}}
\put(42.7, 35.9){\line(0,1){6}}
\put(48.6, 35.5){\line(0,1){6.7}}

\end{overpic}
\begin{overpic}[width=0.32\textwidth,grid=false,tics=10]{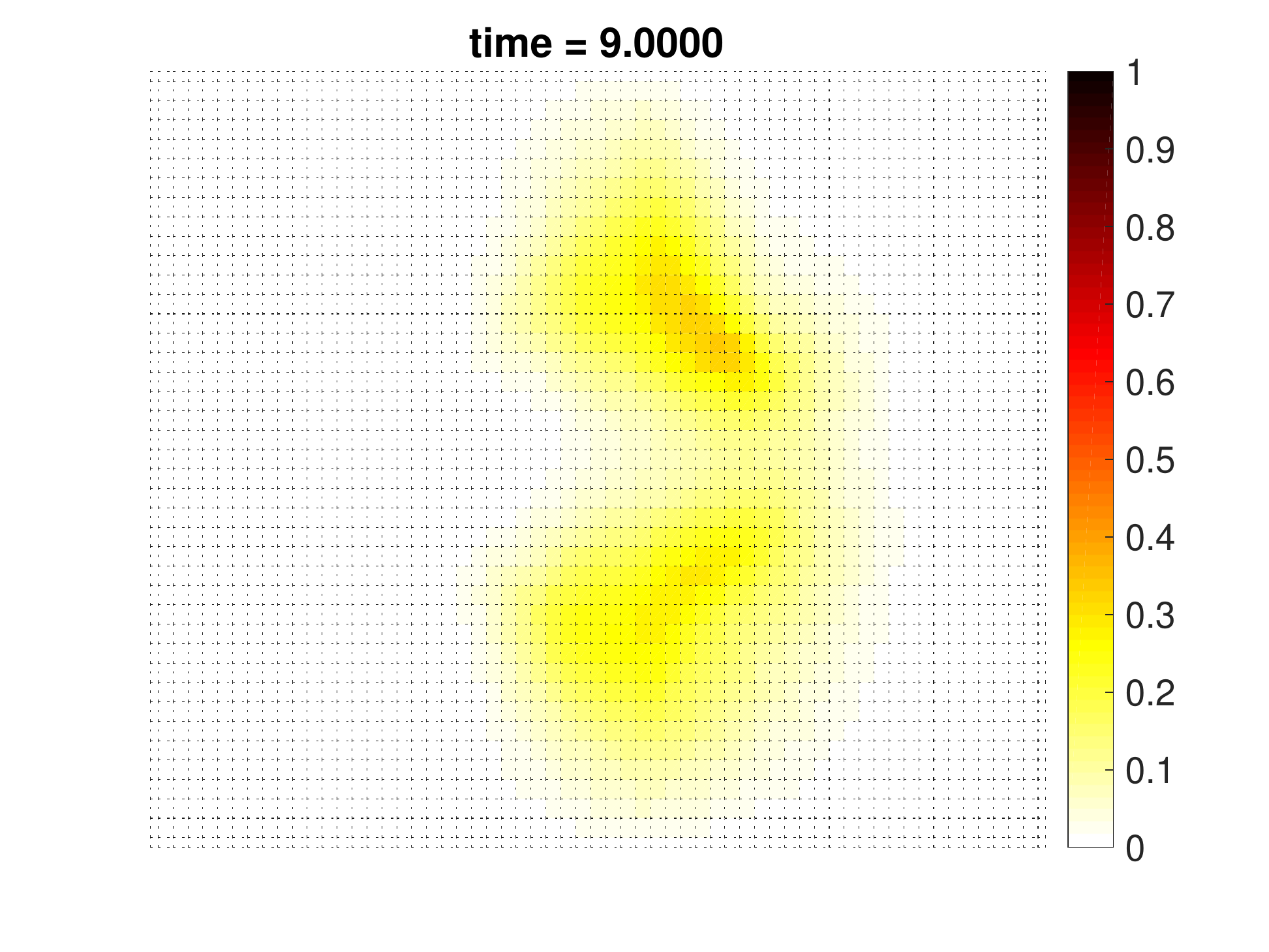}
\linethickness{1pt}
\put(11.7, 69.5){\line(1,0){71}}
\put(11.7, 7.8){\line(0,1){61.9}}
\put(11.7, 8){\line(1,0){71}}
\put(61.2, 7.8){\line(0,1){23}}
\put(61.2, 46.7){\line(0,1){23}}

\put(42.3, 46.5){\line(1,0){9.6}}
\put(42.3, 31.1){\line(1,0){9.6}}
\put(42.5, 31.1){\line(0,1){15.6}}
\put(51.8, 30.9){\line(0,1){15.8}}
\put(42.5, 42){\line(1,0){6}}
\put(42.5, 35.7){\line(1,0){6}}
\put(42.7, 35.9){\line(0,1){6}}
\put(48.6, 35.5){\line(0,1){6.7}}

\end{overpic}
\begin{overpic}[width=0.32\textwidth,grid=false,tics=10]{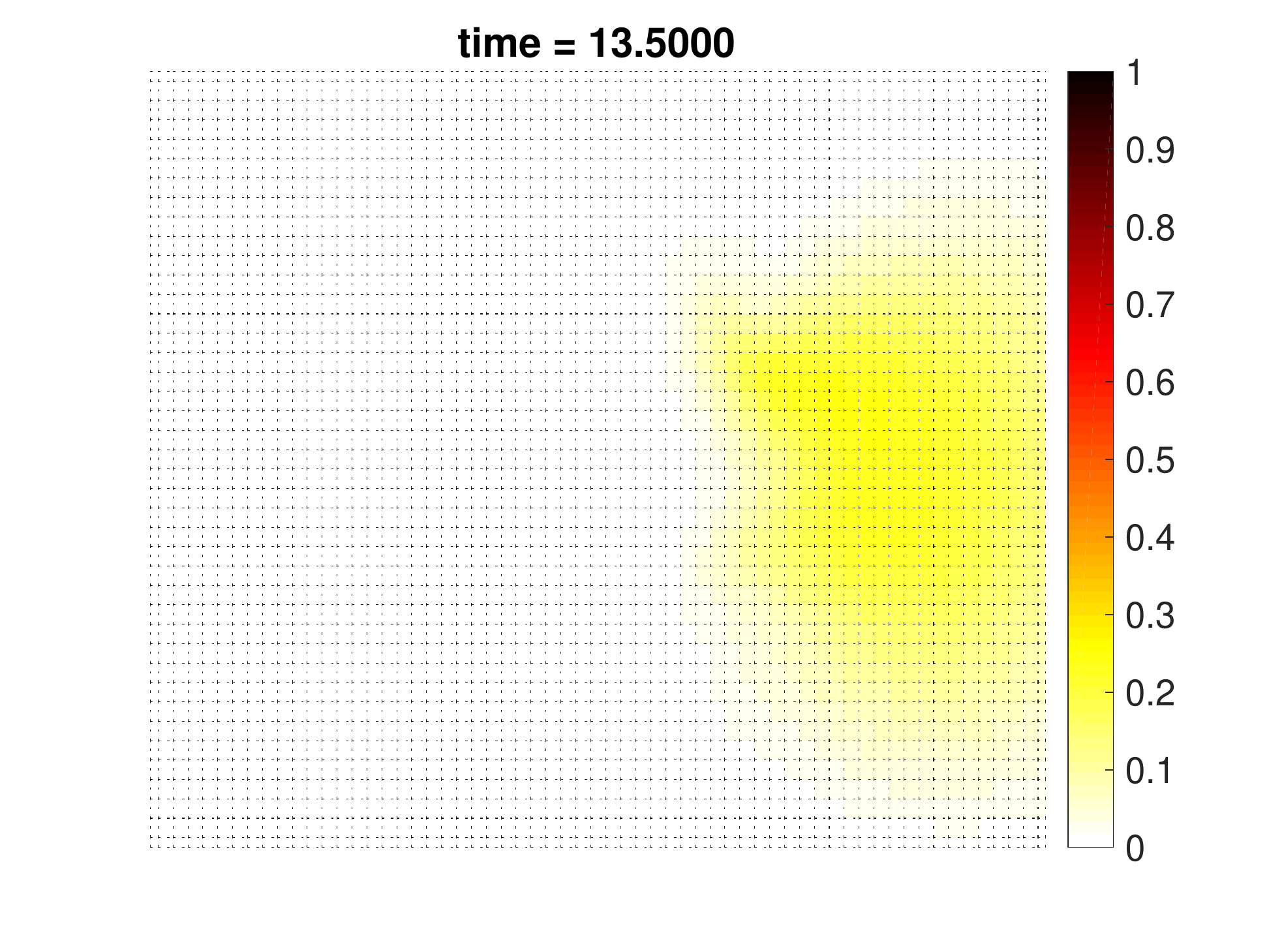}
\linethickness{1pt}
\put(11.7, 69.5){\line(1,0){71}}
\put(11.7, 7.8){\line(0,1){61.9}}
\put(11.7, 8){\line(1,0){71}}
\put(61.2, 7.8){\line(0,1){23}}
\put(61.2, 46.7){\line(0,1){23}}

\put(42.3, 46.5){\line(1,0){9.6}}
\put(42.3, 31.1){\line(1,0){9.6}}
\put(42.5, 31.1){\line(0,1){15.6}}
\put(51.8, 30.9){\line(0,1){15.8}}
\put(42.5, 42){\line(1,0){6}}
\put(42.5, 35.7){\line(1,0){6}}
\put(42.7, 35.9){\line(0,1){6}}
\put(48.6, 35.5){\line(0,1){6.7}}

\end{overpic}
\caption{Configuration 1 with $\alpha=1$ in the effective area: computed density for $t$ = 0, 3, 6, 7.5, 9, 13.5 s. The small square within the effective area represents the real obstacle.}
\label{middle_obstacle}
\end{figure}

\begin{figure}[h!]
\centering
\begin{overpic}[width=0.32\textwidth,grid=false,tics=10]{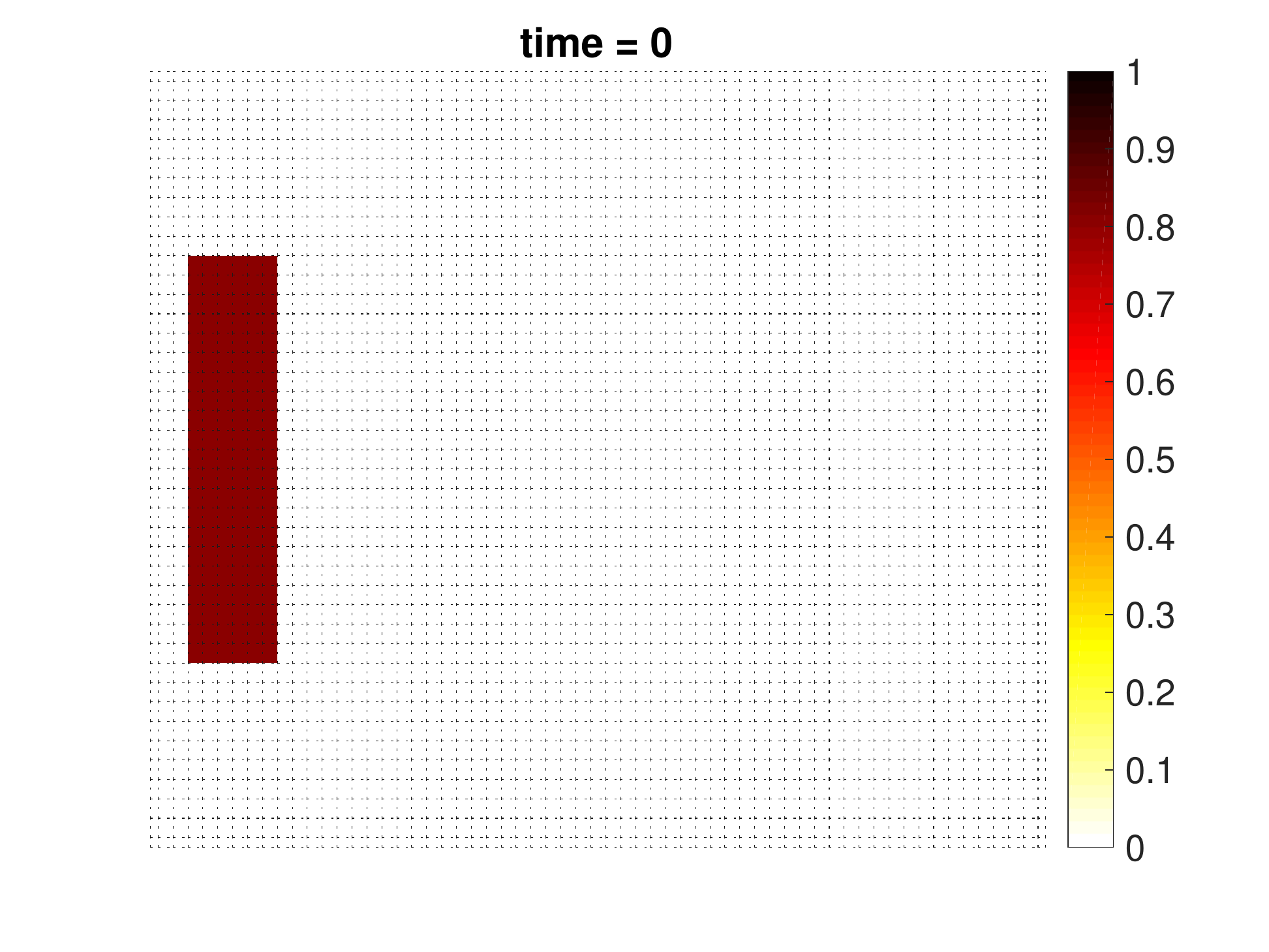}
\linethickness{1pt}
\put(20.3, 38.8){\vector(1,0){6}}
\put(11.7, 69.5){\line(1,0){71}}
\put(11.7, 7.8){\line(0,1){61.9}}
\put(11.7, 8){\line(1,0){71}}
\put(61.2, 7.8){\line(0,1){23}}
\put(61.2, 46.7){\line(0,1){23}}

\put(42.3, 46.5){\line(1,0){9.6}}
\put(42.3, 31.1){\line(1,0){9.6}}
\put(42.5, 31.1){\line(0,1){15.6}}
\put(51.8, 30.9){\line(0,1){15.8}}
\put(48, 34){\line(0,1){9.7}}
\put(51.4, 34){\line(0,1){9.3}}
\put(47.9, 34.2){\line(1,0){3.7}}
\put(47.9, 43.5){\line(1,0){3.7}}

\end{overpic} 
\begin{overpic}[width=0.32\textwidth,grid=false,tics=10]{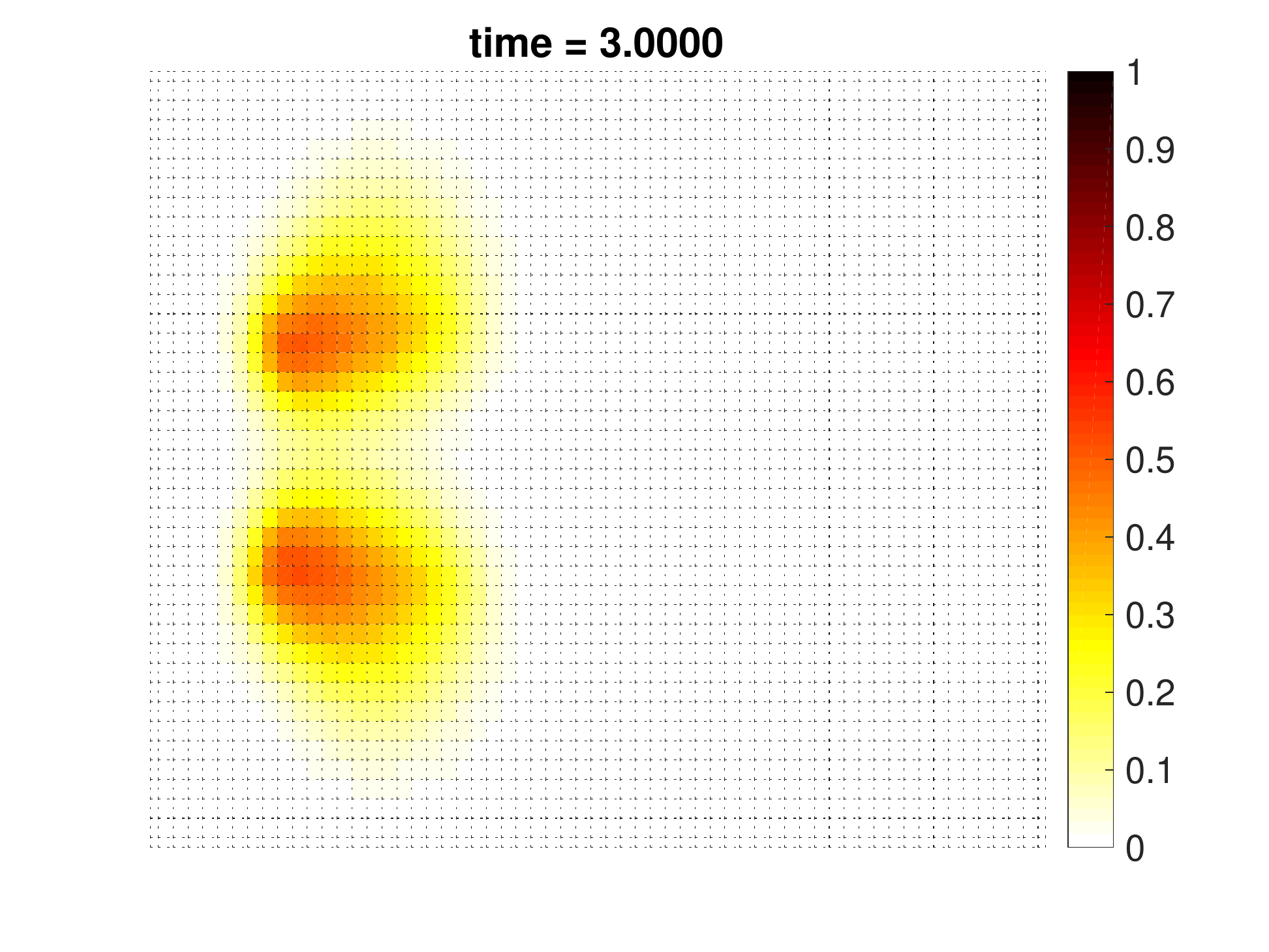}
\linethickness{1pt}
\put(11.7, 69.5){\line(1,0){71}}
\put(11.7, 7.8){\line(0,1){61.9}}
\put(11.7, 8){\line(1,0){71}}
\put(61.2, 7.8){\line(0,1){23}}
\put(61.2, 46.7){\line(0,1){23}}

\put(42.3, 46.5){\line(1,0){9.6}}
\put(42.3, 31.1){\line(1,0){9.6}}
\put(42.5, 31.1){\line(0,1){15.6}}
\put(51.8, 30.9){\line(0,1){15.8}}
\put(48, 34){\line(0,1){9.7}}
\put(51.4, 34){\line(0,1){9.3}}
\put(47.9, 34.2){\line(1,0){3.7}}
\put(47.9, 43.5){\line(1,0){3.7}}

\end{overpic} 
\begin{overpic}[width=0.32\textwidth,grid=false,tics=10]{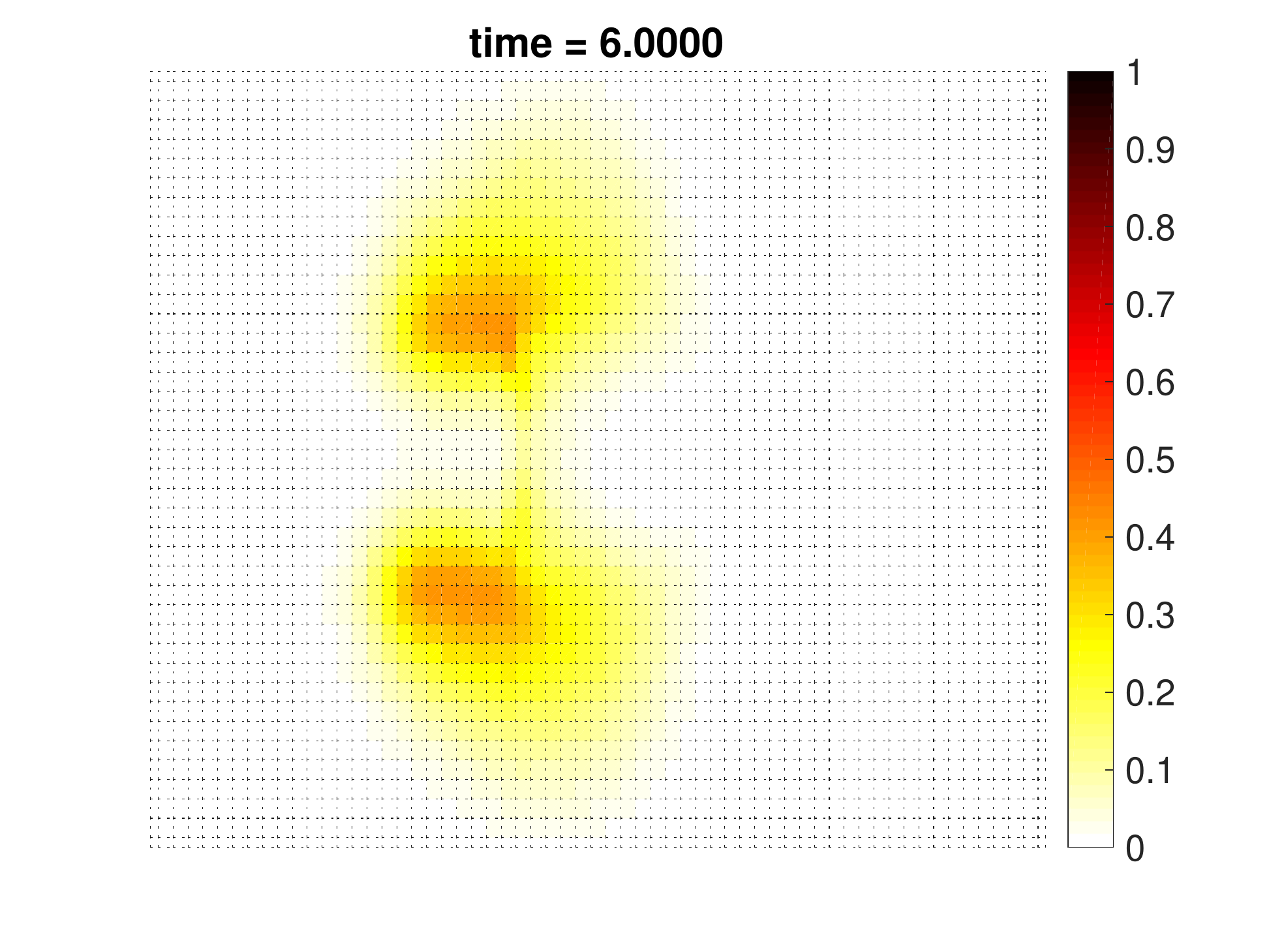}
\linethickness{1pt}
\put(11.7, 69.5){\line(1,0){71}}
\put(11.7, 7.8){\line(0,1){61.9}}
\put(11.7, 8){\line(1,0){71}}
\put(61.2, 7.8){\line(0,1){23}}
\put(61.2, 46.7){\line(0,1){23}}

\put(42.3, 46.5){\line(1,0){9.6}}
\put(42.3, 31.1){\line(1,0){9.6}}
\put(42.5, 31.1){\line(0,1){15.6}}
\put(51.8, 30.9){\line(0,1){15.8}}
\put(48, 34){\line(0,1){9.7}}
\put(51.4, 34){\line(0,1){9.3}}
\put(47.9, 34.2){\line(1,0){3.7}}
\put(47.9, 43.5){\line(1,0){3.7}}

\end{overpic}
\begin{overpic}[width=0.32\textwidth,grid=false,tics=10]{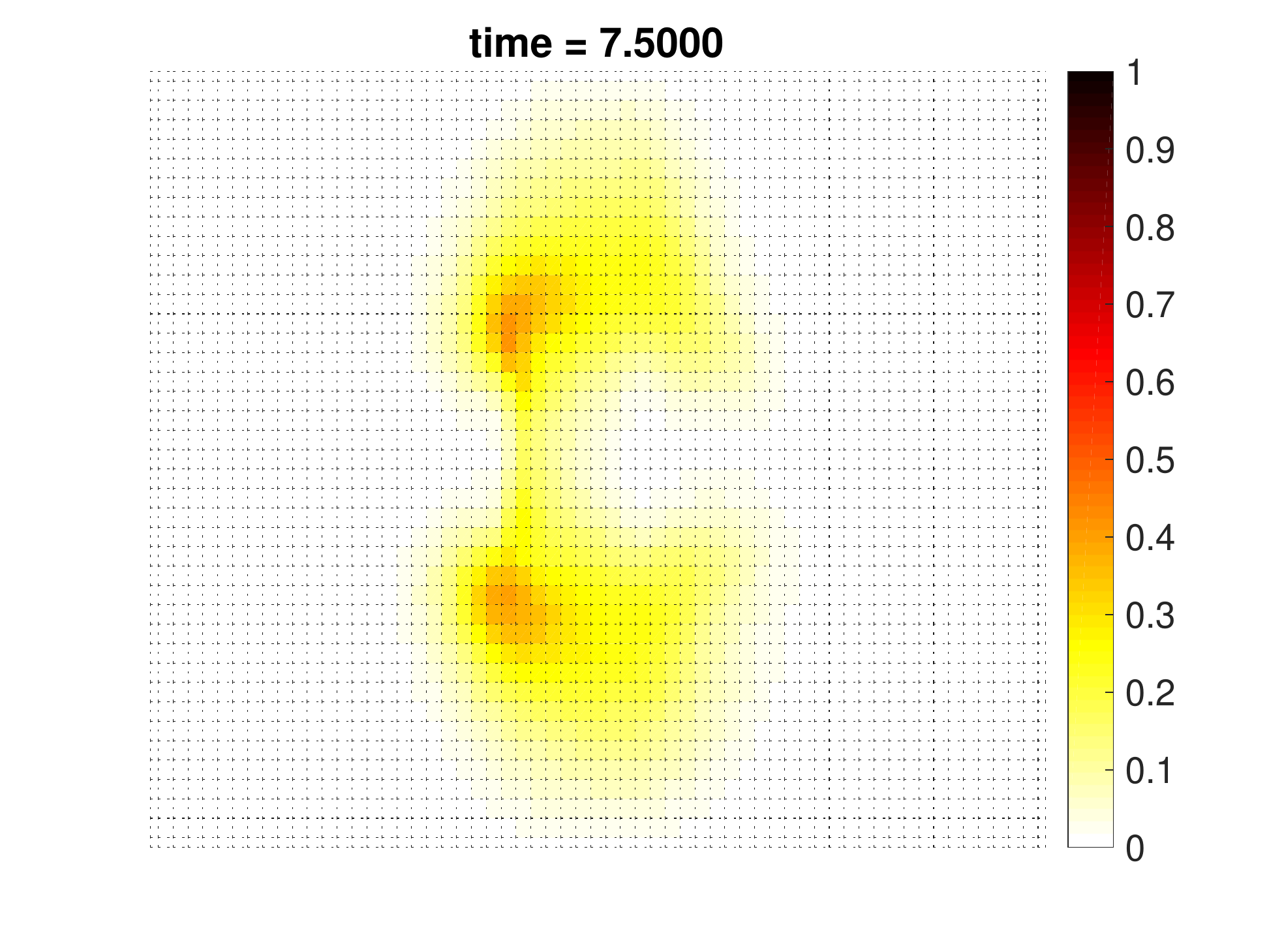}
\linethickness{1pt}
\put(11.7, 69.5){\line(1,0){71}}
\put(11.7, 7.8){\line(0,1){61.9}}
\put(11.7, 8){\line(1,0){71}}
\put(61.2, 7.8){\line(0,1){23}}
\put(61.2, 46.7){\line(0,1){23}}

\put(42.3, 46.5){\line(1,0){9.6}}
\put(42.3, 31.1){\line(1,0){9.6}}
\put(42.5, 31.1){\line(0,1){15.6}}
\put(51.8, 30.9){\line(0,1){15.8}}
\put(48, 34){\line(0,1){9.7}}
\put(51.4, 34){\line(0,1){9.3}}
\put(47.9, 34.2){\line(1,0){3.7}}
\put(47.9, 43.5){\line(1,0){3.7}}

\end{overpic}
\begin{overpic}[width=0.32\textwidth,grid=false,tics=10]{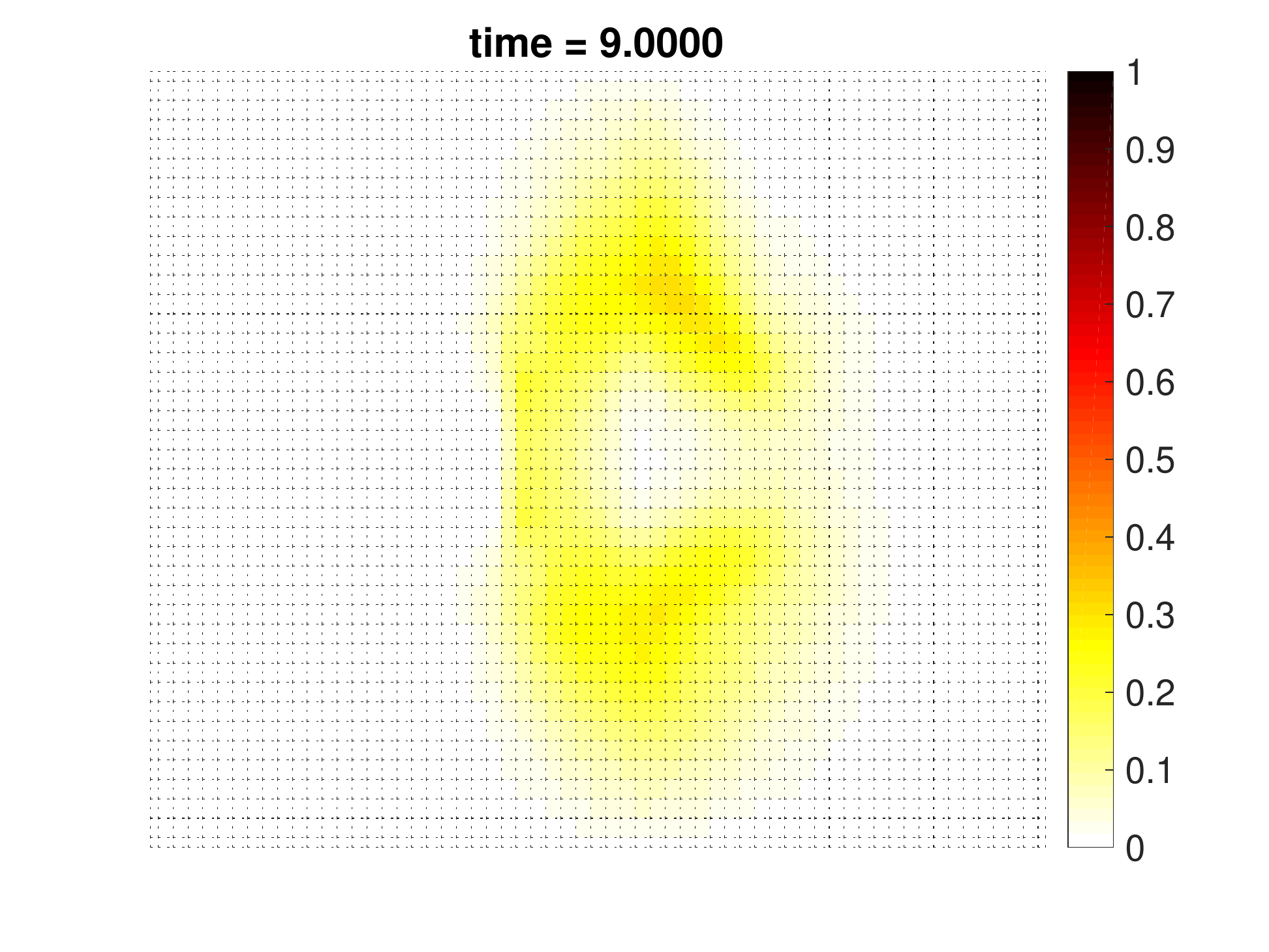}
\linethickness{1pt}
\put(11.7, 69.5){\line(1,0){71}}
\put(11.7, 7.8){\line(0,1){61.9}}
\put(11.7, 8){\line(1,0){71}}
\put(61.2, 7.8){\line(0,1){23}}
\put(61.2, 46.7){\line(0,1){23}}

\put(42.3, 46.5){\line(1,0){9.6}}
\put(42.3, 31.1){\line(1,0){9.6}}
\put(42.5, 31.1){\line(0,1){15.6}}
\put(51.8, 30.9){\line(0,1){15.8}}
\put(48, 34){\line(0,1){9.7}}
\put(51.4, 34){\line(0,1){9.3}}
\put(47.9, 34.2){\line(1,0){3.7}}
\put(47.9, 43.5){\line(1,0){3.7}}

\end{overpic}
\begin{overpic}[width=0.32\textwidth,grid=false,tics=10]{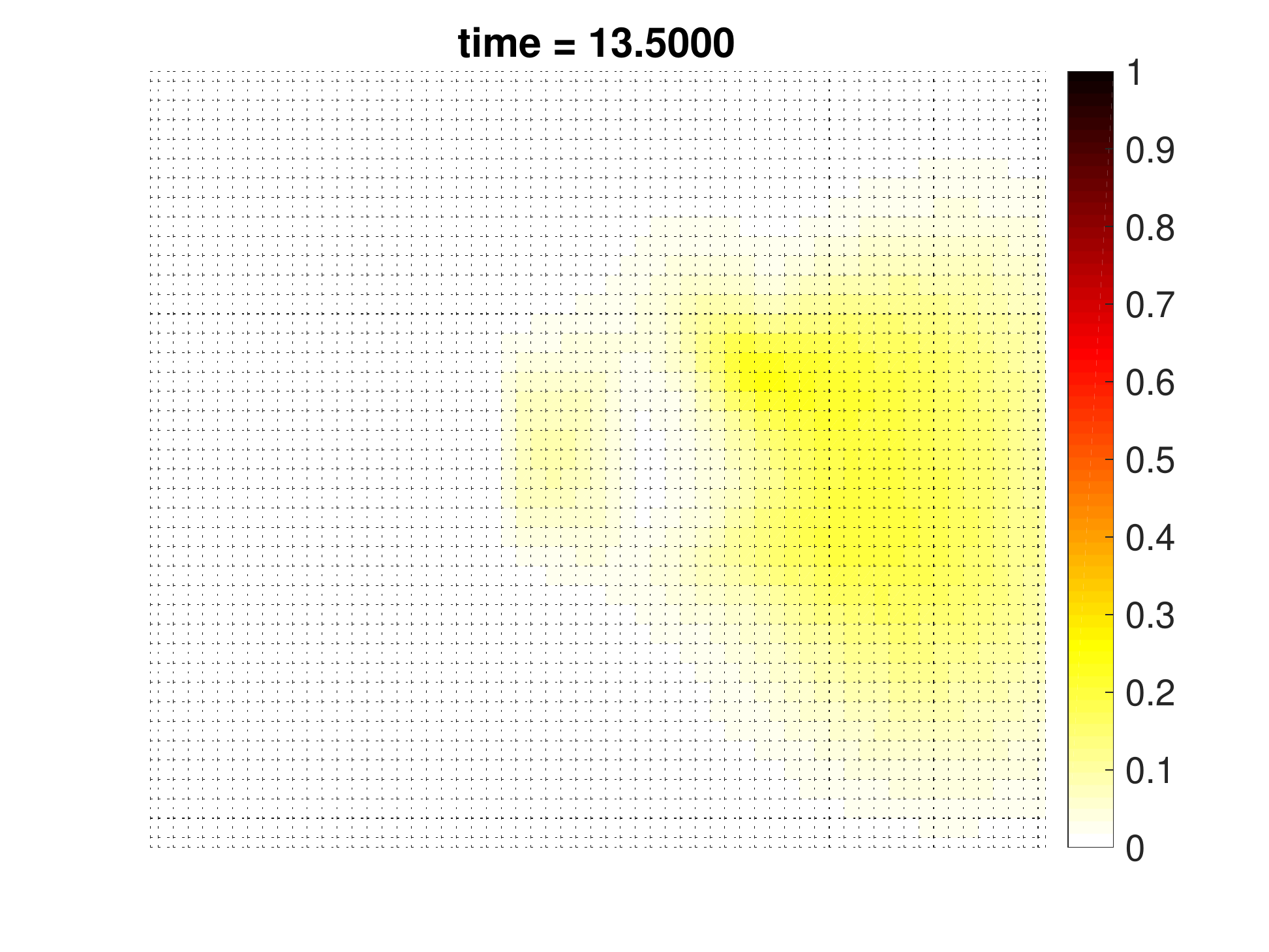}
\linethickness{1pt}
\put(11.7, 69.5){\line(1,0){71}}
\put(11.7, 7.8){\line(0,1){61.9}}
\put(11.7, 8){\line(1,0){71}}
\put(61.2, 7.8){\line(0,1){23}}
\put(61.2, 46.7){\line(0,1){23}}

\put(42.3, 46.5){\line(1,0){9.6}}
\put(42.3, 31.1){\line(1,0){9.6}}
\put(42.5, 31.1){\line(0,1){15.6}}
\put(51.8, 30.9){\line(0,1){15.8}}
\put(48, 34){\line(0,1){9.7}}
\put(51.4, 34){\line(0,1){9.3}}
\put(47.9, 34.2){\line(1,0){3.7}}
\put(47.9, 43.5){\line(1,0){3.7}}

\end{overpic}
\caption{ Configuration 1 with $\alpha=0$ in the effective area: computed density for $t$ = 0, 3, 6, 7.5, 9, 13.5 s. The small rectangle within the effective area represents the real obstacle.}
\label{middle_obstacle_0}
\end{figure}

\begin{figure}[h!] 
\centering
\begin{overpic}[width=0.32\textwidth,grid=false,tics=10]{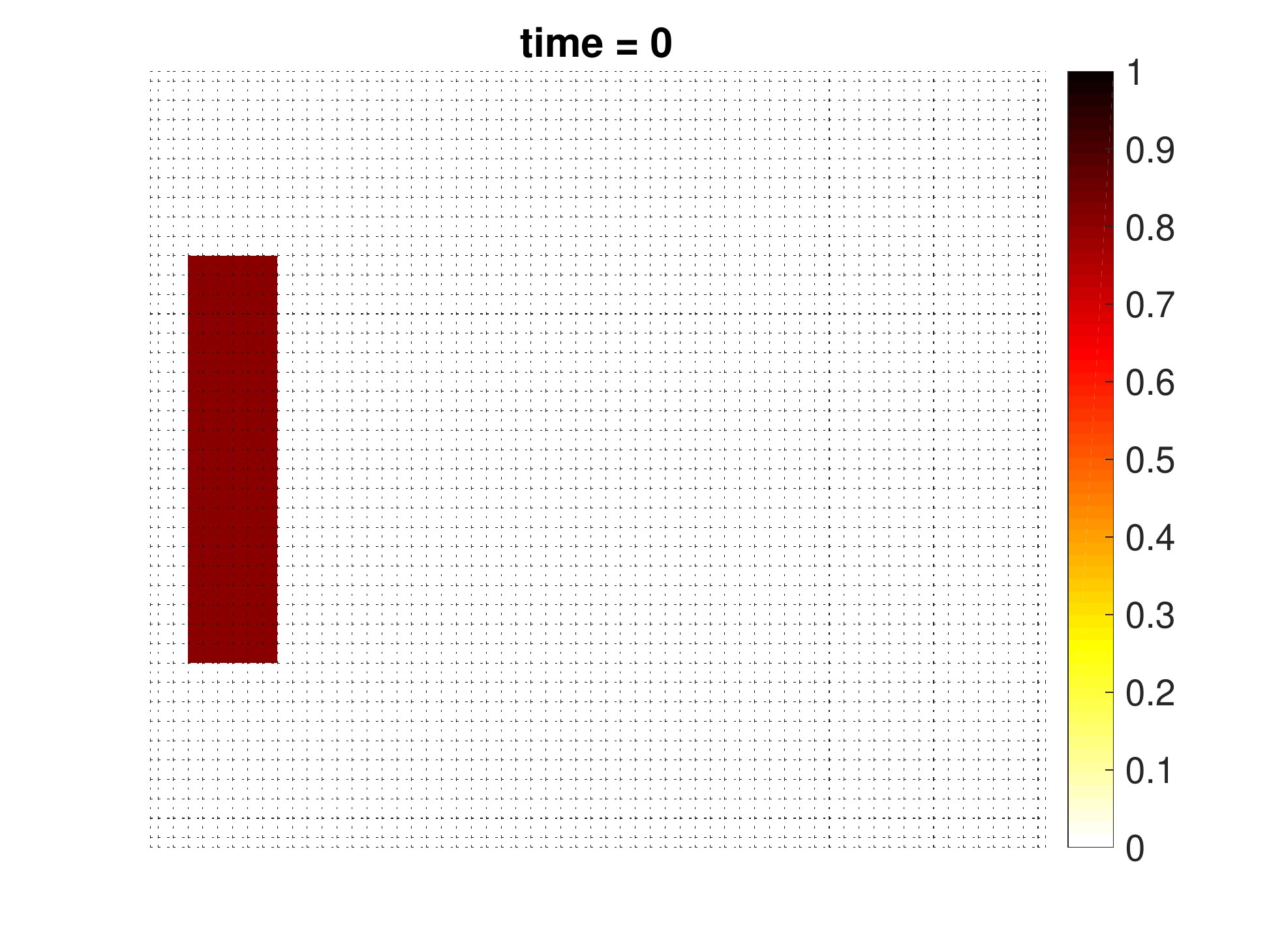}
\linethickness{1pt}
\put(20.3, 38.8){\vector(1,0){6}}
\put(11.7, 69.5){\line(1,0){71}}
\put(11.7, 7.8){\line(0,1){61.9}}
\put(11.7, 8){\line(1,0){71}}
\put(61.2, 7.8){\line(0,1){23}}
\put(61.2, 46.7){\line(0,1){23}}

\put(42.3, 46.5){\line(1,0){9.6}}
\put(42.3, 61.9){\line(1,0){9.6}}
\put(42.5, 46.5){\line(0,1){15.6}}
\put(51.8, 46.3){\line(0,1){15.8}}

\put(42.3, 14.3){\line(1,0){9.6}}
\put(42.3, 29.7){\line(1,0){9.6}}
\put(42.5, 14.3){\line(0,1){15.6}}
\put(51.8, 14.1){\line(0,1){15.8}}

\put(42.5, 51){\line(1,0){6}}
\put(42.5, 57.3){\line(1,0){6}}
\put(42.7, 51.2){\line(0,1){6}}
\put(48.4, 50.8){\line(0,1){6.7}}

\put(42.5, 18.8){\line(1,0){6}}
\put(42.5, 25.1){\line(1,0){6}}
\put(42.7, 19){\line(0,1){6}}
\put(48.4, 18.6){\line(0,1){6.7}}

\end{overpic} 
\begin{overpic}[width=0.32\textwidth,grid=false,tics=10]{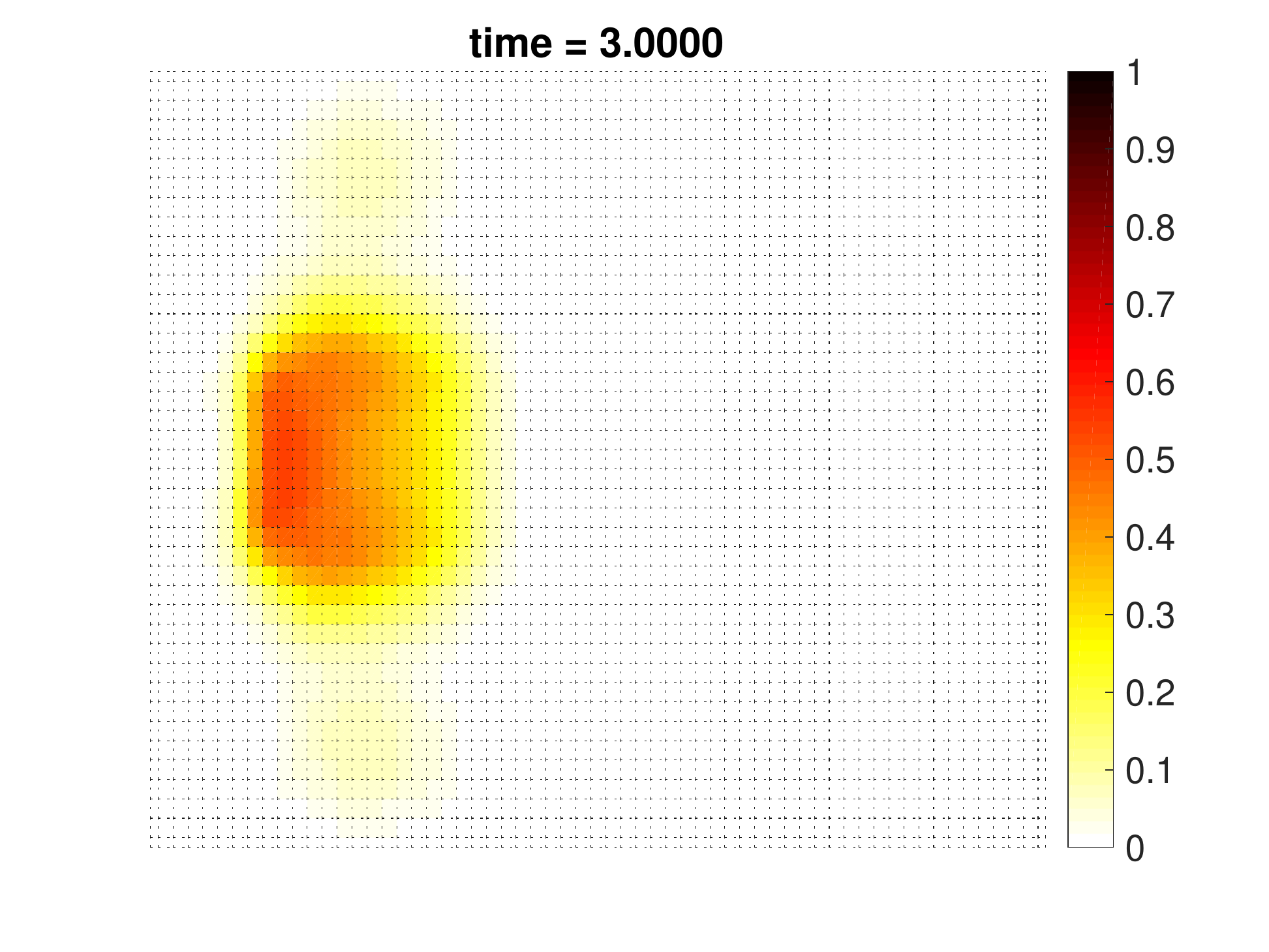}
\linethickness{1pt}
\put(11.7, 69.5){\line(1,0){71}}
\put(11.7, 7.8){\line(0,1){61.9}}
\put(11.7, 8){\line(1,0){71}}
\put(61.2, 7.8){\line(0,1){23}}
\put(61.2, 46.7){\line(0,1){23}}

\put(42.3, 46.5){\line(1,0){9.6}}
\put(42.3, 61.9){\line(1,0){9.6}}
\put(42.5, 46.5){\line(0,1){15.6}}
\put(51.8, 46.3){\line(0,1){15.8}}

\put(42.3, 14.3){\line(1,0){9.6}}
\put(42.3, 29.7){\line(1,0){9.6}}
\put(42.5, 14.3){\line(0,1){15.6}}
\put(51.8, 14.1){\line(0,1){15.8}}

\put(42.5, 51){\line(1,0){6}}
\put(42.5, 57.3){\line(1,0){6}}
\put(42.7, 51.2){\line(0,1){6}}
\put(48.4, 50.8){\line(0,1){6.7}}

\put(42.5, 18.8){\line(1,0){6}}
\put(42.5, 25.1){\line(1,0){6}}
\put(42.7, 19){\line(0,1){6}}
\put(48.4, 18.6){\line(0,1){6.7}}
\end{overpic} 
\begin{overpic}[width=0.32\textwidth,grid=false,tics=10]{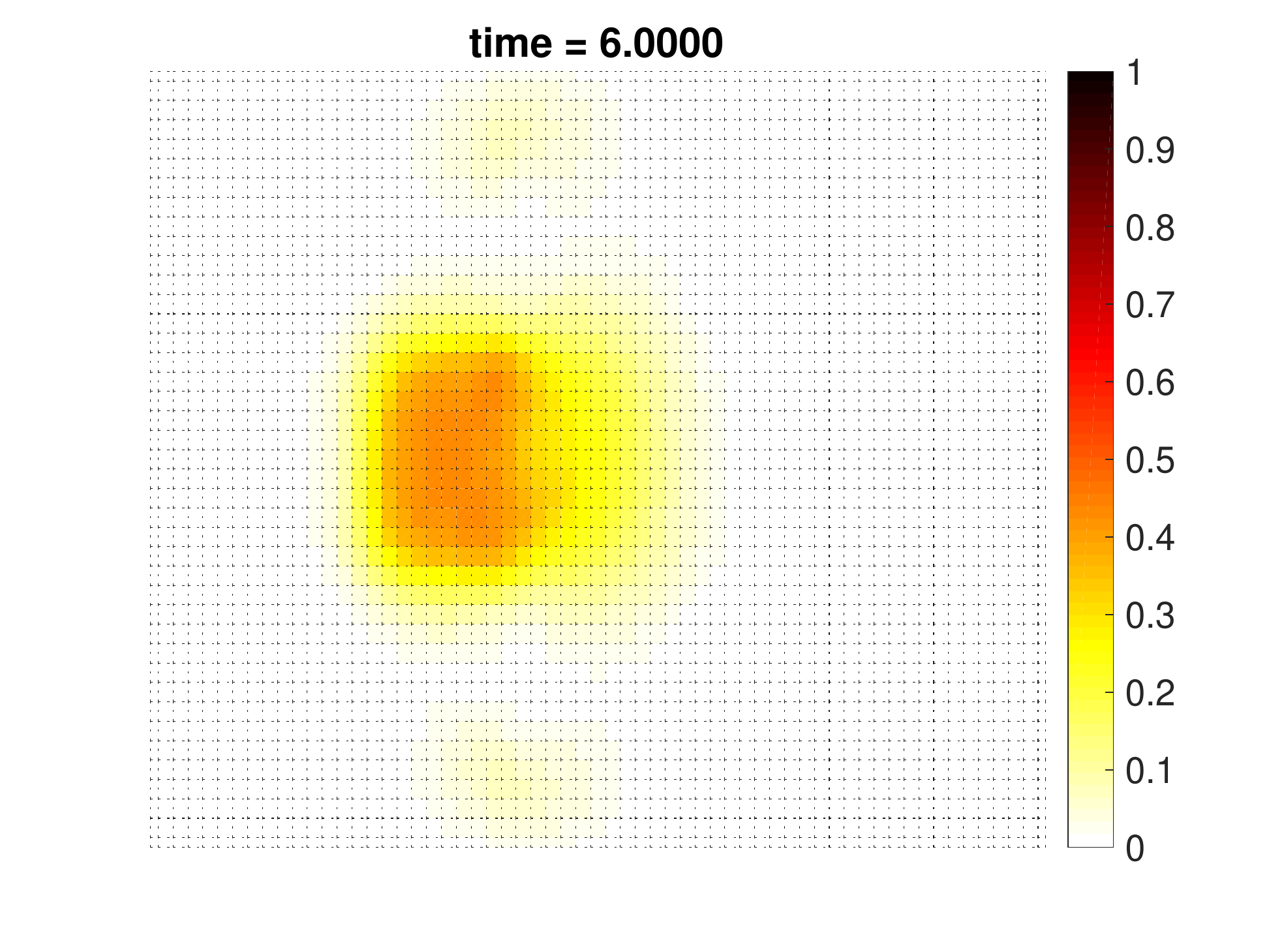}
\linethickness{1pt}
\put(11.7, 69.5){\line(1,0){71}}
\put(11.7, 7.8){\line(0,1){61.9}}
\put(11.7, 8){\line(1,0){71}}
\put(61.2, 7.8){\line(0,1){23}}
\put(61.2, 46.7){\line(0,1){23}}

\put(42.3, 46.5){\line(1,0){9.6}}
\put(42.3, 61.9){\line(1,0){9.6}}
\put(42.5, 46.5){\line(0,1){15.6}}
\put(51.8, 46.3){\line(0,1){15.8}}

\put(42.3, 14.3){\line(1,0){9.6}}
\put(42.3, 29.7){\line(1,0){9.6}}
\put(42.5, 14.3){\line(0,1){15.6}}
\put(51.8, 14.1){\line(0,1){15.8}}

\put(42.5, 51){\line(1,0){6}}
\put(42.5, 57.3){\line(1,0){6}}
\put(42.7, 51.2){\line(0,1){6}}
\put(48.4, 50.8){\line(0,1){6.7}}

\put(42.5, 18.8){\line(1,0){6}}
\put(42.5, 25.1){\line(1,0){6}}
\put(42.7, 19){\line(0,1){6}}
\put(48.4, 18.6){\line(0,1){6.7}}
\end{overpic}
\begin{overpic}[width=0.32\textwidth,grid=false,tics=10]{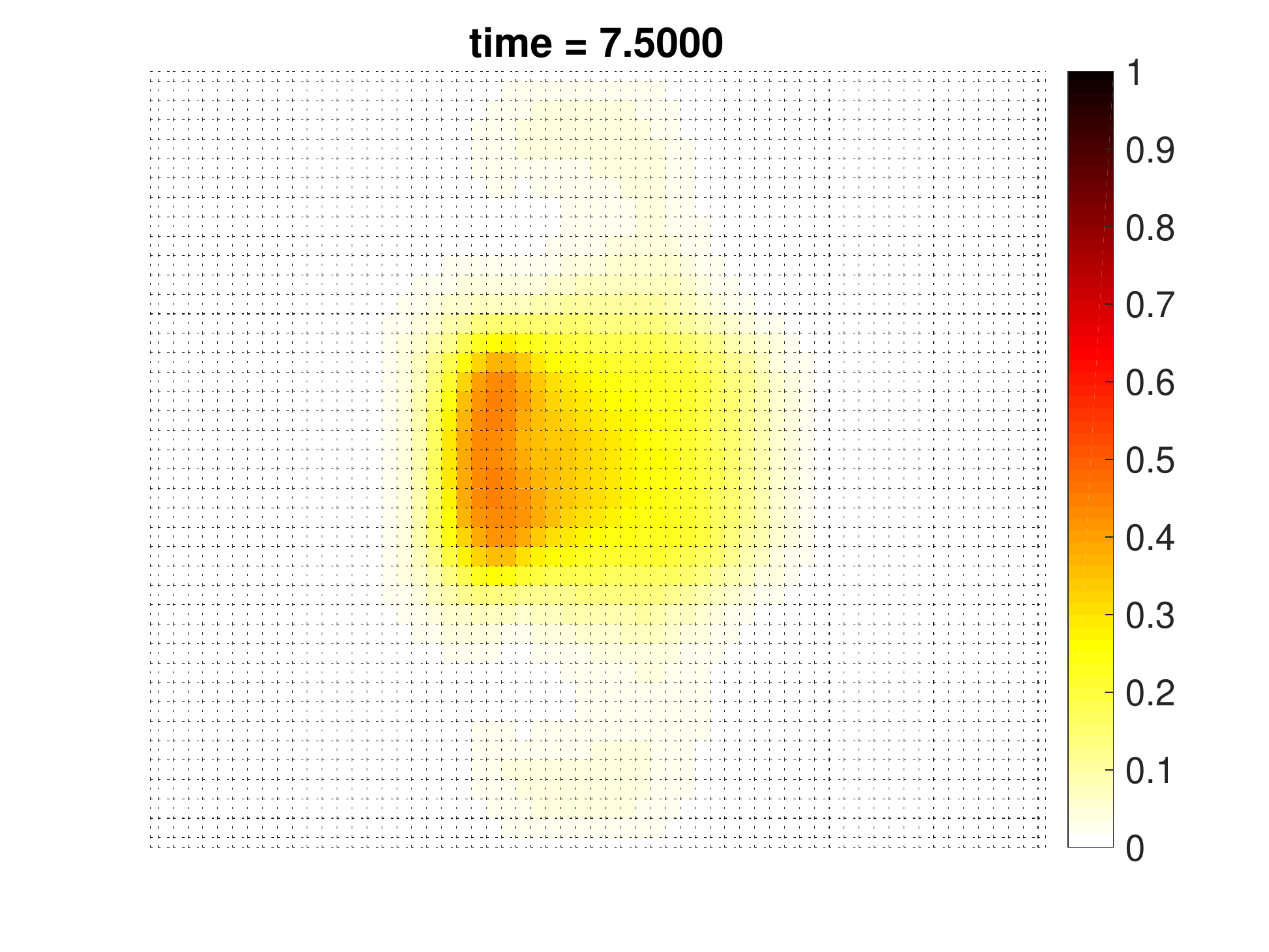}
\linethickness{1pt}
\put(11.7, 69.5){\line(1,0){71}}
\put(11.7, 7.8){\line(0,1){61.9}}
\put(11.7, 8){\line(1,0){71}}
\put(61.2, 7.8){\line(0,1){23}}
\put(61.2, 46.7){\line(0,1){23}}

\put(42.3, 46.5){\line(1,0){9.6}}
\put(42.3, 61.9){\line(1,0){9.6}}
\put(42.5, 46.5){\line(0,1){15.6}}
\put(51.8, 46.3){\line(0,1){15.8}}

\put(42.3, 14.3){\line(1,0){9.6}}
\put(42.3, 29.7){\line(1,0){9.6}}
\put(42.5, 14.3){\line(0,1){15.6}}
\put(51.8, 14.1){\line(0,1){15.8}}

\put(42.5, 51){\line(1,0){6}}
\put(42.5, 57.3){\line(1,0){6}}
\put(42.7, 51.2){\line(0,1){6}}
\put(48.4, 50.8){\line(0,1){6.7}}

\put(42.5, 18.8){\line(1,0){6}}
\put(42.5, 25.1){\line(1,0){6}}
\put(42.7, 19){\line(0,1){6}}
\put(48.4, 18.6){\line(0,1){6.7}}
\end{overpic}
\begin{overpic}[width=0.32\textwidth,grid=false,tics=10]{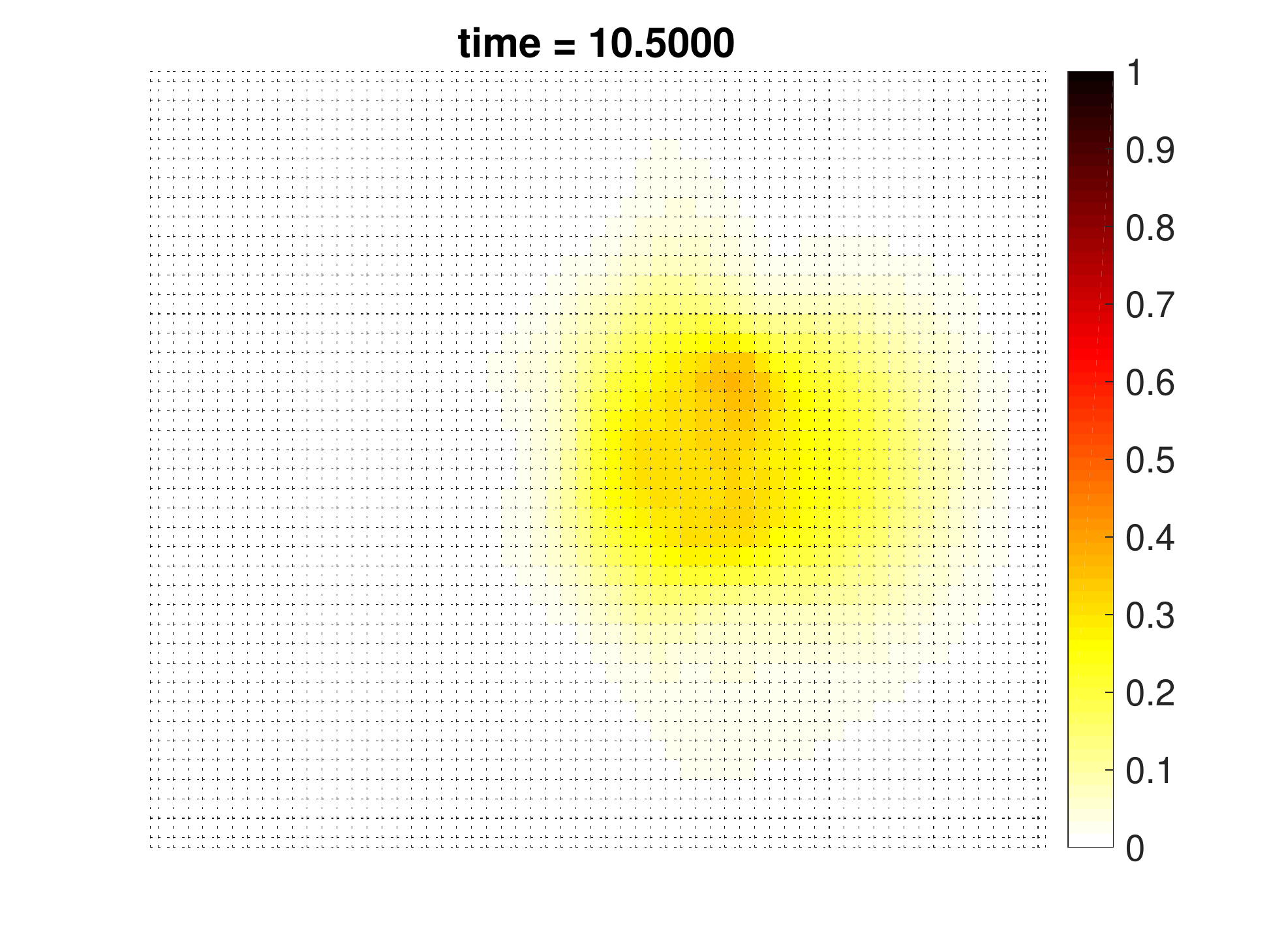}
\linethickness{1pt}
\put(11.7, 69.5){\line(1,0){71}}
\put(11.7, 7.8){\line(0,1){61.9}}
\put(11.7, 8){\line(1,0){71}}
\put(61.2, 7.8){\line(0,1){23}}
\put(61.2, 46.7){\line(0,1){23}}

\put(42.3, 46.5){\line(1,0){9.6}}
\put(42.3, 61.9){\line(1,0){9.6}}
\put(42.5, 46.5){\line(0,1){15.6}}
\put(51.8, 46.3){\line(0,1){15.8}}

\put(42.3, 14.3){\line(1,0){9.6}}
\put(42.3, 29.7){\line(1,0){9.6}}
\put(42.5, 14.3){\line(0,1){15.6}}
\put(51.8, 14.1){\line(0,1){15.8}}

\put(42.5, 51){\line(1,0){6}}
\put(42.5, 57.3){\line(1,0){6}}
\put(42.7, 51.2){\line(0,1){6}}
\put(48.4, 50.8){\line(0,1){6.7}}

\put(42.5, 18.8){\line(1,0){6}}
\put(42.5, 25.1){\line(1,0){6}}
\put(42.7, 19){\line(0,1){6}}
\put(48.4, 18.6){\line(0,1){6.7}}
\end{overpic}
\begin{overpic}[width=0.32\textwidth,grid=false,tics=10]{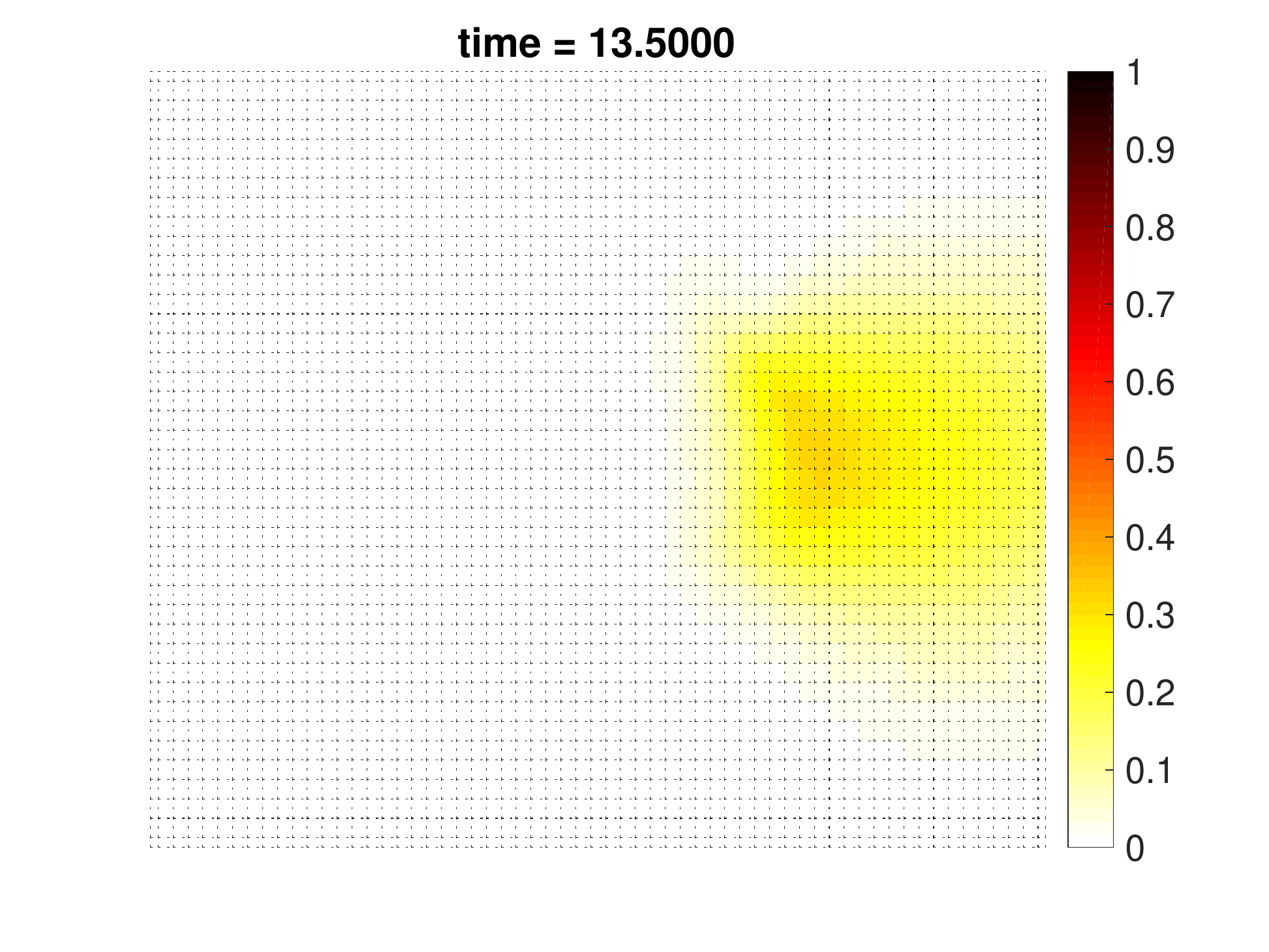}
\linethickness{1pt}
\put(11.7, 69.5){\line(1,0){71}}
\put(11.7, 7.8){\line(0,1){61.9}}
\put(11.7, 8){\line(1,0){71}}
\put(61.2, 7.8){\line(0,1){23}}
\put(61.2, 46.7){\line(0,1){23}}

\put(42.3, 46.5){\line(1,0){9.6}}
\put(42.3, 61.9){\line(1,0){9.6}}
\put(42.5, 46.5){\line(0,1){15.6}}
\put(51.8, 46.3){\line(0,1){15.8}}

\put(42.3, 14.3){\line(1,0){9.6}}
\put(42.3, 29.7){\line(1,0){9.6}}
\put(42.5, 14.3){\line(0,1){15.6}}
\put(51.8, 14.1){\line(0,1){15.8}}

\put(42.5, 51){\line(1,0){6}}
\put(42.5, 57.3){\line(1,0){6}}
\put(42.7, 51.2){\line(0,1){6}}
\put(48.4, 50.8){\line(0,1){6.7}}

\put(42.5, 18.8){\line(1,0){6}}
\put(42.5, 25.1){\line(1,0){6}}
\put(42.7, 19){\line(0,1){6}}
\put(48.4, 18.6){\line(0,1){6.7}}
\end{overpic}
\caption{Configuration 2 with $\alpha=1$ in the effective area: computed density for $t$ = 0, 3, 6, 7.5, 10.5, 13.5 s. The small square within the effective area represents the real obstacle.}
\label{two_obstacle}
\end{figure}

\begin{figure} 
\centering
\begin{overpic}[width=0.32\textwidth,grid=false,tics=10]{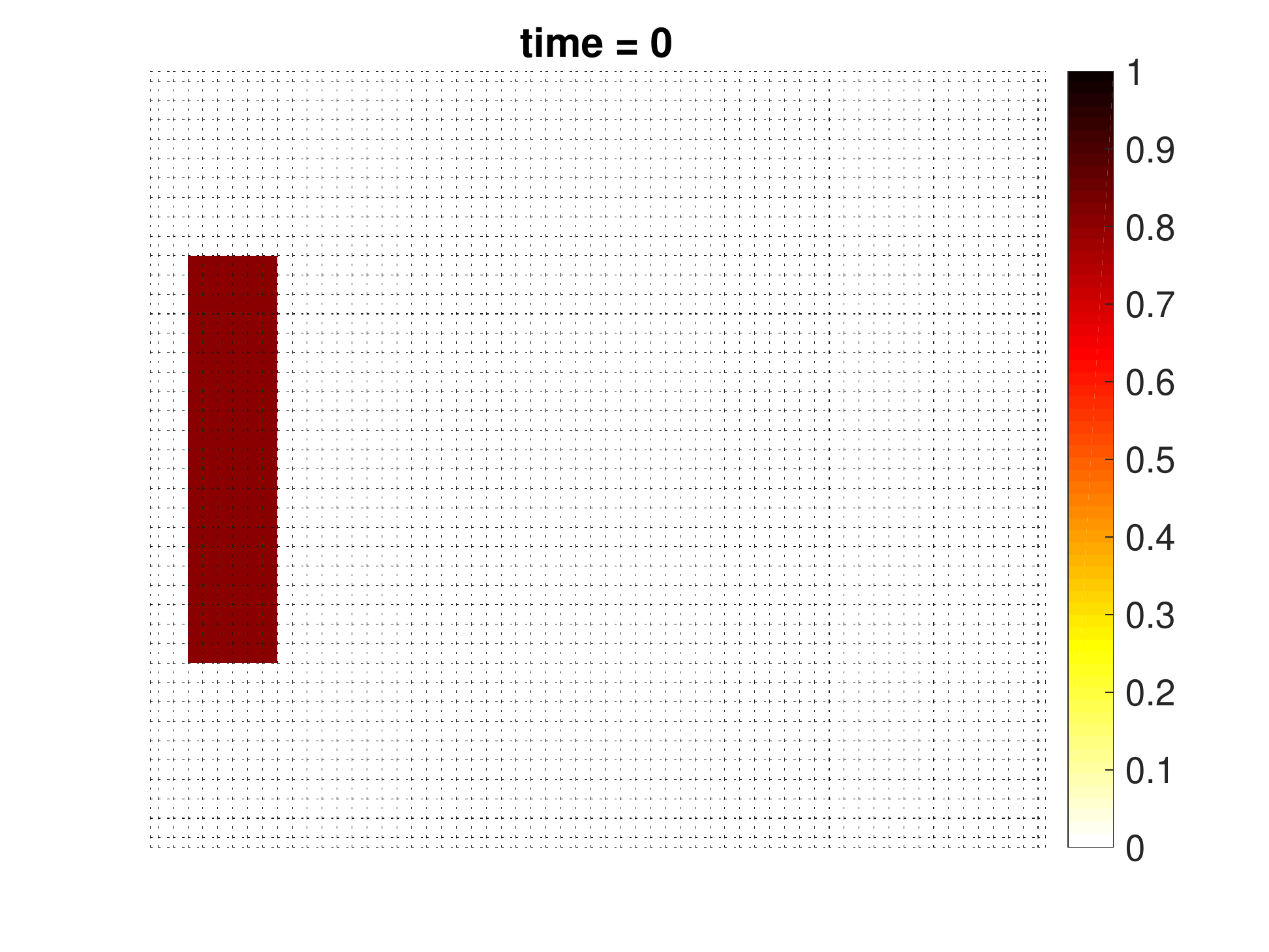}
\linethickness{1pt}
\put(20.3, 38.8){\vector(1,0){6}}
\put(11.7, 69.5){\line(1,0){71}}
\put(11.7, 7.8){\line(0,1){61.9}}
\put(11.7, 8){\line(1,0){71}}
\put(61.2, 7.8){\line(0,1){23}}
\put(61.2, 46.7){\line(0,1){23}}

\put(42.3, 46.5){\line(1,0){9.6}}
\put(42.3, 61.9){\line(1,0){9.6}}
\put(42.5, 46.5){\line(0,1){15.6}}
\put(51.8, 46.3){\line(0,1){15.8}}

\put(42.3, 14.3){\line(1,0){9.6}}
\put(42.3, 29.7){\line(1,0){9.6}}
\put(42.5, 14.3){\line(0,1){15.6}}
\put(51.8, 14.1){\line(0,1){15.8}}

\put(48, 49.3){\line(0,1){9.7}}
\put(51.4, 49.3){\line(0,1){9.3}}
\put(47.9, 49.5){\line(1,0){3.7}}
\put(47.9, 58.8){\line(1,0){3.7}}

\put(48, 17.1){\line(0,1){9.7}}
\put(51.4, 17.1){\line(0,1){9.3}}
\put(47.9, 17.2){\line(1,0){3.7}}
\put(47.9, 26.6){\line(1,0){3.7}}
\end{overpic} 
\begin{overpic}[width=0.32\textwidth,grid=false,tics=10]{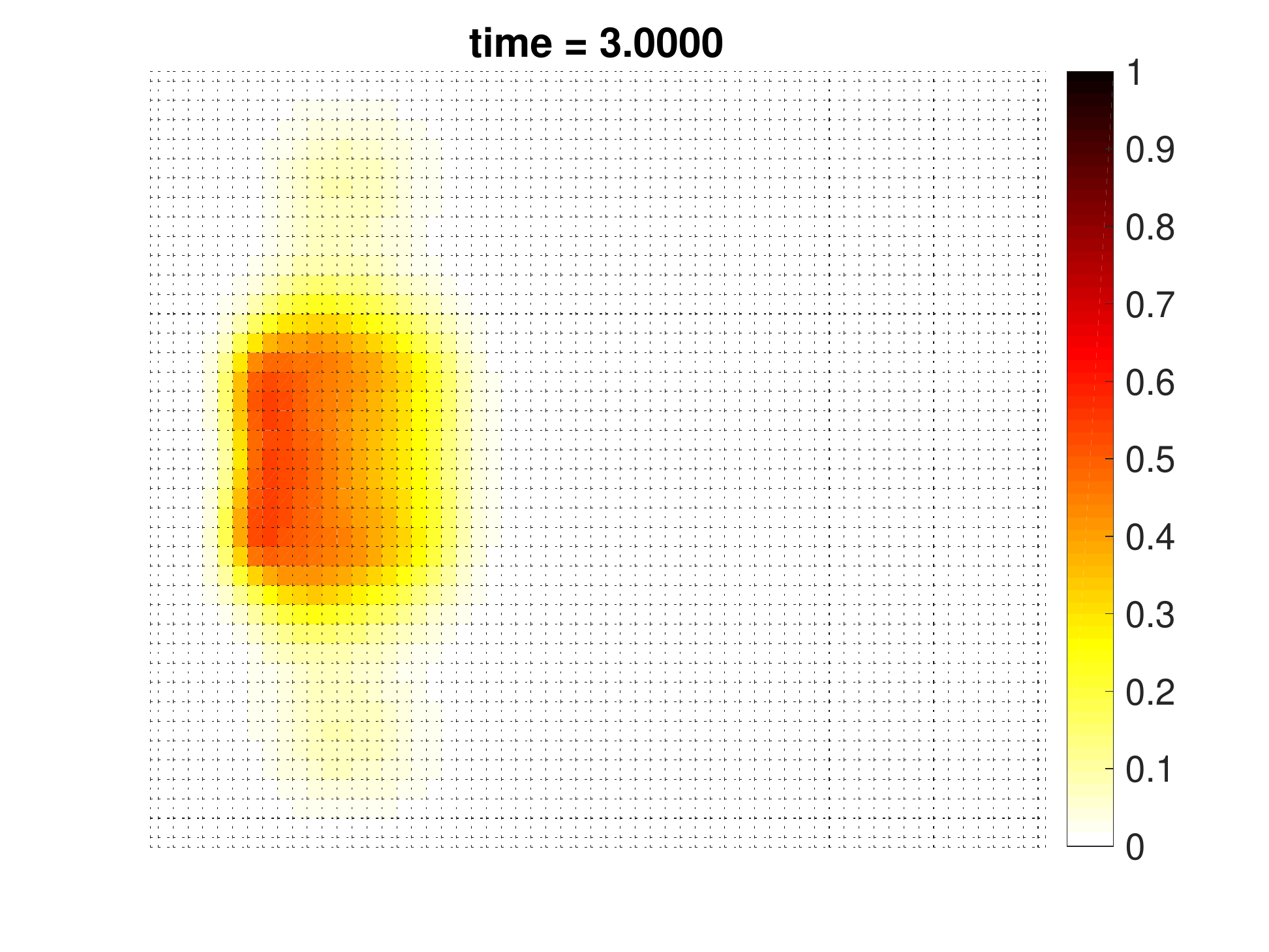}
\linethickness{1pt}
\put(11.7, 69.5){\line(1,0){71}}
\put(11.7, 7.8){\line(0,1){61.9}}
\put(11.7, 8){\line(1,0){71}}
\put(61.2, 7.8){\line(0,1){23}}
\put(61.2, 46.7){\line(0,1){23}}

\put(42.3, 46.5){\line(1,0){9.6}}
\put(42.3, 61.9){\line(1,0){9.6}}
\put(42.5, 46.5){\line(0,1){15.6}}
\put(51.8, 46.3){\line(0,1){15.8}}

\put(42.3, 14.3){\line(1,0){9.6}}
\put(42.3, 29.7){\line(1,0){9.6}}
\put(42.5, 14.3){\line(0,1){15.6}}
\put(51.8, 14.1){\line(0,1){15.8}}

\put(48, 49.3){\line(0,1){9.7}}
\put(51.4, 49.3){\line(0,1){9.3}}
\put(47.9, 49.5){\line(1,0){3.7}}
\put(47.9, 58.8){\line(1,0){3.7}}

\put(48, 17.1){\line(0,1){9.7}}
\put(51.4, 17.1){\line(0,1){9.3}}
\put(47.9, 17.2){\line(1,0){3.7}}
\put(47.9, 26.6){\line(1,0){3.7}}
\end{overpic} 
\begin{overpic}[width=0.32\textwidth,grid=false,tics=10]{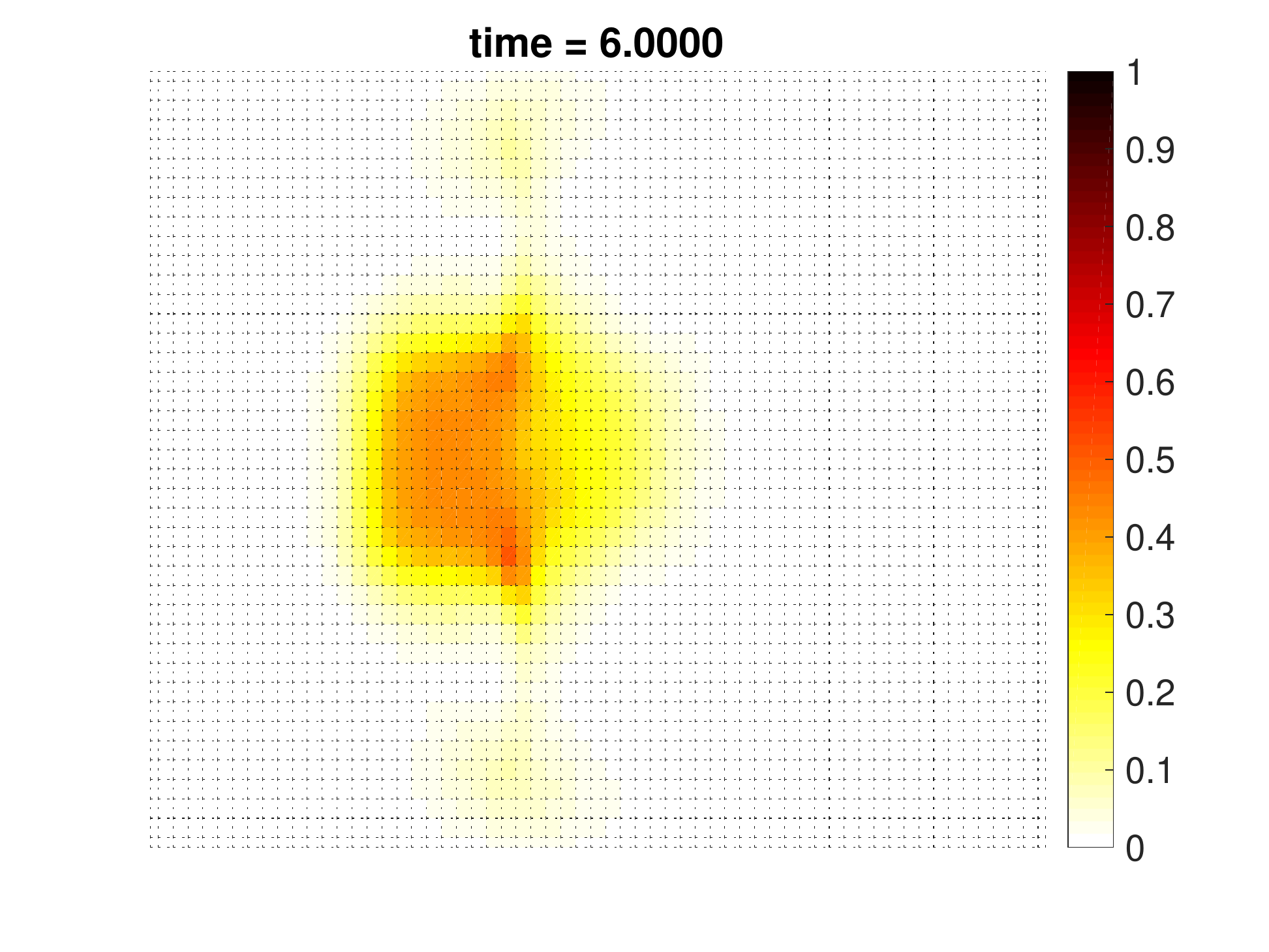}
\linethickness{1pt}
\put(11.7, 69.5){\line(1,0){71}}
\put(11.7, 7.8){\line(0,1){61.9}}
\put(11.7, 8){\line(1,0){71}}
\put(61.2, 7.8){\line(0,1){23}}
\put(61.2, 46.7){\line(0,1){23}}

\put(42.3, 46.5){\line(1,0){9.6}}
\put(42.3, 61.9){\line(1,0){9.6}}
\put(42.5, 46.5){\line(0,1){15.6}}
\put(51.8, 46.3){\line(0,1){15.8}}

\put(42.3, 14.3){\line(1,0){9.6}}
\put(42.3, 29.7){\line(1,0){9.6}}
\put(42.5, 14.3){\line(0,1){15.6}}
\put(51.8, 14.1){\line(0,1){15.8}}

\put(48, 49.3){\line(0,1){9.7}}
\put(51.4, 49.3){\line(0,1){9.3}}
\put(47.9, 49.5){\line(1,0){3.7}}
\put(47.9, 58.8){\line(1,0){3.7}}

\put(48, 17.1){\line(0,1){9.7}}
\put(51.4, 17.1){\line(0,1){9.3}}
\put(47.9, 17.2){\line(1,0){3.7}}
\put(47.9, 26.6){\line(1,0){3.7}}
\end{overpic}
\begin{overpic}[width=0.32\textwidth,grid=false,tics=10]{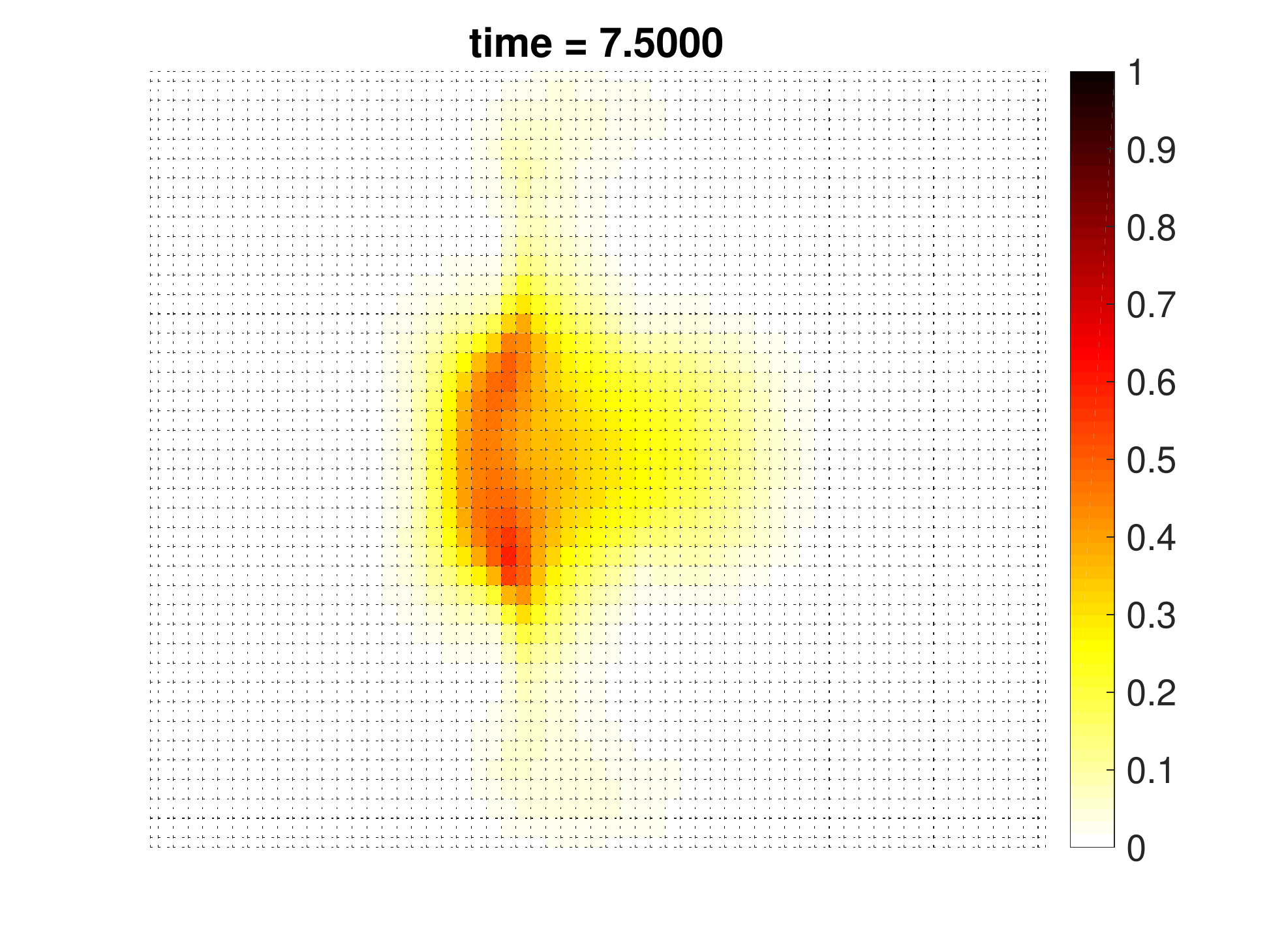}
\linethickness{1pt}
\put(11.7, 69.5){\line(1,0){71}}
\put(11.7, 7.8){\line(0,1){61.9}}
\put(11.7, 8){\line(1,0){71}}
\put(61.2, 7.8){\line(0,1){23}}
\put(61.2, 46.7){\line(0,1){23}}

\put(42.3, 46.5){\line(1,0){9.6}}
\put(42.3, 61.9){\line(1,0){9.6}}
\put(42.5, 46.5){\line(0,1){15.6}}
\put(51.8, 46.3){\line(0,1){15.8}}

\put(42.3, 14.3){\line(1,0){9.6}}
\put(42.3, 29.7){\line(1,0){9.6}}
\put(42.5, 14.3){\line(0,1){15.6}}
\put(51.8, 14.1){\line(0,1){15.8}}

\put(48, 49.3){\line(0,1){9.7}}
\put(51.4, 49.3){\line(0,1){9.3}}
\put(47.9, 49.5){\line(1,0){3.7}}
\put(47.9, 58.8){\line(1,0){3.7}}

\put(48, 17.1){\line(0,1){9.7}}
\put(51.4, 17.1){\line(0,1){9.3}}
\put(47.9, 17.2){\line(1,0){3.7}}
\put(47.9, 26.6){\line(1,0){3.7}}
\end{overpic}
\begin{overpic}[width=0.32\textwidth,grid=false,tics=10]{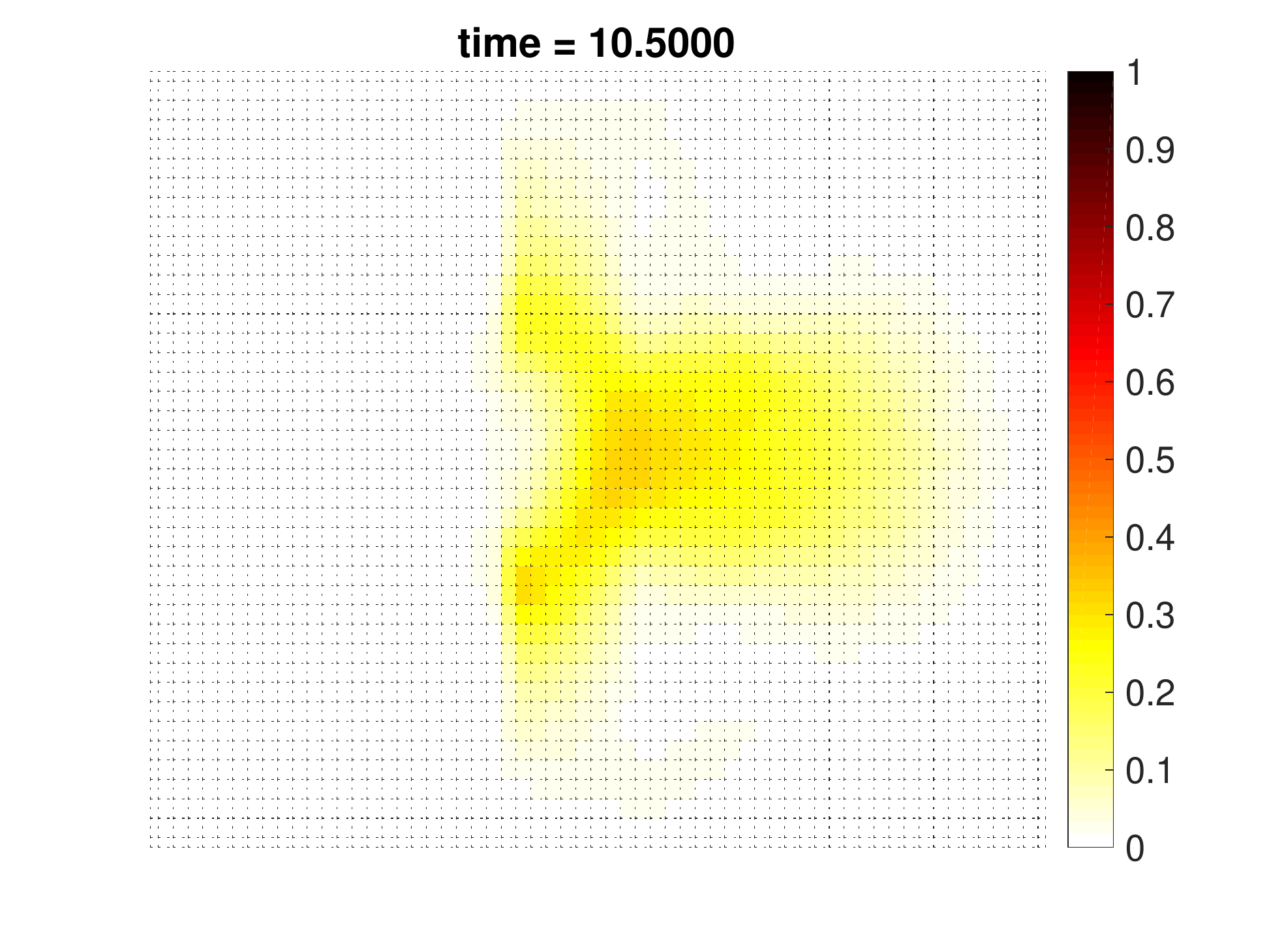}
\linethickness{1pt}
\put(11.7, 69.5){\line(1,0){71}}
\put(11.7, 7.8){\line(0,1){61.9}}
\put(11.7, 8){\line(1,0){71}}
\put(61.2, 7.8){\line(0,1){23}}
\put(61.2, 46.7){\line(0,1){23}}

\put(42.3, 46.5){\line(1,0){9.6}}
\put(42.3, 61.9){\line(1,0){9.6}}
\put(42.5, 46.5){\line(0,1){15.6}}
\put(51.8, 46.3){\line(0,1){15.8}}

\put(42.3, 14.3){\line(1,0){9.6}}
\put(42.3, 29.7){\line(1,0){9.6}}
\put(42.5, 14.3){\line(0,1){15.6}}
\put(51.8, 14.1){\line(0,1){15.8}}

\put(48, 49.3){\line(0,1){9.7}}
\put(51.4, 49.3){\line(0,1){9.3}}
\put(47.9, 49.5){\line(1,0){3.7}}
\put(47.9, 58.8){\line(1,0){3.7}}

\put(48, 17.1){\line(0,1){9.7}}
\put(51.4, 17.1){\line(0,1){9.3}}
\put(47.9, 17.2){\line(1,0){3.7}}
\put(47.9, 26.6){\line(1,0){3.7}}
\end{overpic}
\begin{overpic}[width=0.32\textwidth,grid=false,tics=10]{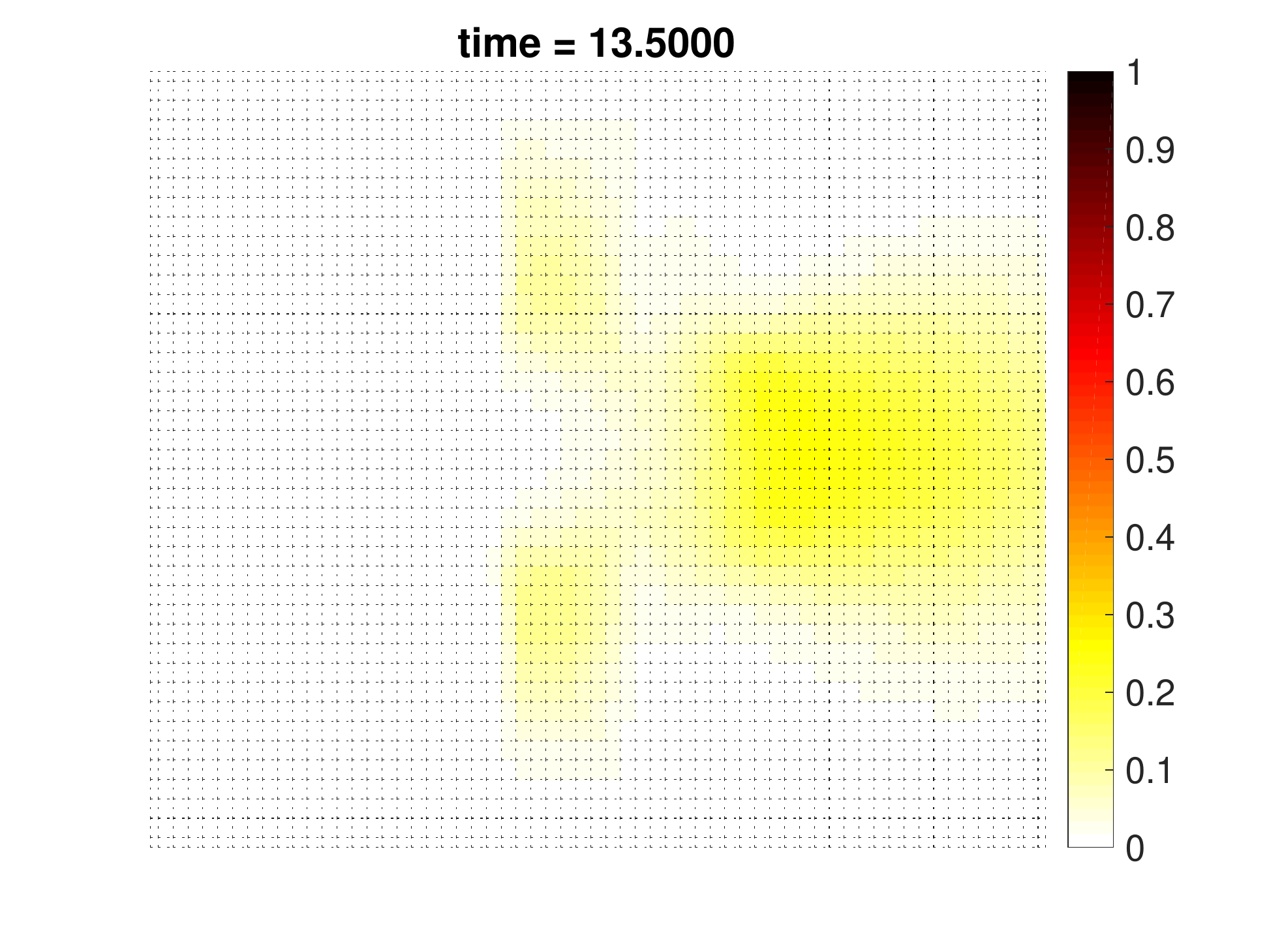}
\linethickness{1pt}
\put(11.7, 69.5){\line(1,0){71}}
\put(11.7, 7.8){\line(0,1){61.9}}
\put(11.7, 8){\line(1,0){71}}
\put(61.2, 7.8){\line(0,1){23}}
\put(61.2, 46.7){\line(0,1){23}}

\put(42.3, 46.5){\line(1,0){9.6}}
\put(42.3, 61.9){\line(1,0){9.6}}
\put(42.5, 46.5){\line(0,1){15.6}}
\put(51.8, 46.3){\line(0,1){15.8}}

\put(42.3, 14.3){\line(1,0){9.6}}
\put(42.3, 29.7){\line(1,0){9.6}}
\put(42.5, 14.3){\line(0,1){15.6}}
\put(51.8, 14.1){\line(0,1){15.8}}

\put(48, 49.3){\line(0,1){9.7}}
\put(51.4, 49.3){\line(0,1){9.3}}
\put(47.9, 49.5){\line(1,0){3.7}}
\put(47.9, 58.8){\line(1,0){3.7}}

\put(48, 17.1){\line(0,1){9.7}}
\put(51.4, 17.1){\line(0,1){9.3}}
\put(47.9, 17.2){\line(1,0){3.7}}
\put(47.9, 26.6){\line(1,0){3.7}}
\end{overpic}
\caption{Configuration 2 with $\alpha=0$ in the effective area: computed density for $t$ = 0, 3, 6, 7.5, 10.5, 13.5 s s. The small rectangle within the effective area represents the real obstacle.}
\label{two_obstacle_0}
\end{figure}

Figure~\ref{middle_obstacle} and \ref{middle_obstacle_0} show the evacuation process 
for configuration 1 when $\alpha=1$ and $\alpha=0$ in the effective area, respectively. 
As explained in Section~\ref{sec:obstacles}, 
when $\alpha=1$ (resp., $\alpha=0$) pedestrians avoid the front (resp., back) 
part of the effective area. Moreover, the shape of the real obstacle is different:
it is square for $\alpha=1$, while for $\alpha=0$
it is slender. These findings are confirmed by Figure~\ref{two_obstacle} and
\ref{two_obstacle_0}, which show the evacuation process 
for configuration 2 when $\alpha=1$ and $\alpha=0$ in the effective area, respectively.

\begin{figure}[h!]
\centering
\begin{overpic}[width=0.35\textwidth,grid=false,tics=10]{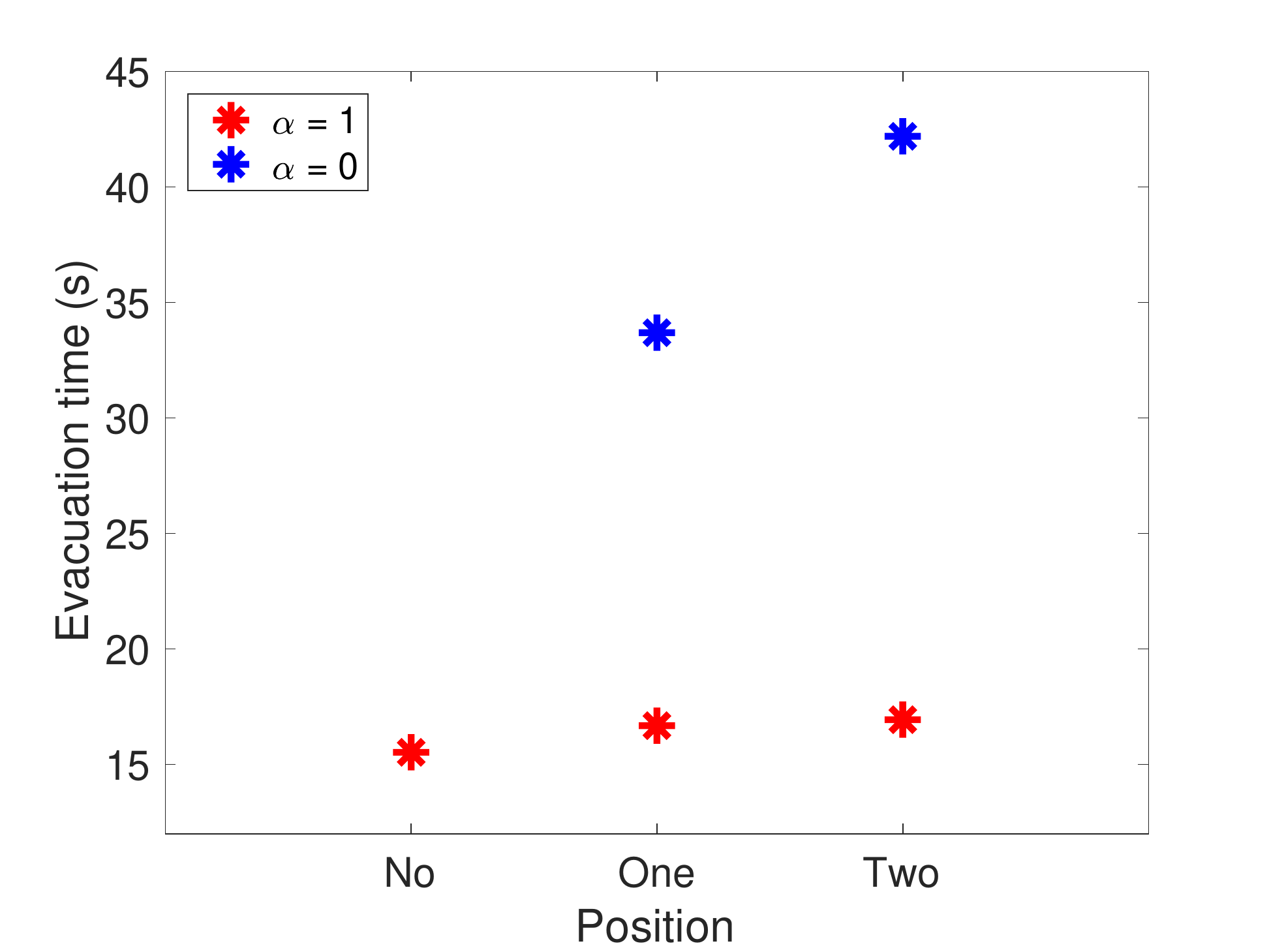}
\end{overpic}
\caption{Evacuation times for the room 
with no obstacles ($\alpha =1$ everywhere in the domain), for room with one and two obstacles
with $\alpha=1$ and $\alpha=0$ in the effective area.}
\label{obstacle_graph}
\end{figure}

Finally, Figure~\ref{obstacle_graph} compares the evacuation times for the room 
with no obstacles ($\alpha =1$ everywhere in the domain), and for configurations 1 and 2 
for $\alpha=1$ and $\alpha=0$ in the effective area. 
Obviously, the shortest evacuation time is for the room with no obstacles
and overall good quality of the environment. 
The evacuation time is slightly larger when there is one or two obstacles
and $\alpha$ is equal to 1 in the effective area. 
When we compare Figure~\ref{middle_obstacle} and \ref{middle_obstacle_0} at $t = 13.5$ s, 
it seems that roughly the same number of people is left in the room. However, 
the small group of people in the effective area in Figure~\ref{middle_obstacle_0} has
velocity modulus close 0 since $\alpha$ is equal to $0$ there and thus the overall
evacuation process takes longer. The same happens in Figure~\ref{two_obstacle_0}.

\subsection{Lane formation}\label{sec:lane}

In this section, we assess the ability of our model to reproduce an 
empirically observed phenomenon: 
formation of lanes in a corridor when two or more groups of pedestrians have opposite 
walking directions \cite{Helbing2002}.
We consider a computational domain of length $L = 20$ m  and width $H = 5$ m
and we impose periodic boundary conditions on the short edges. 
The reference quantities we use for this test are: $D=5\sqrt{17}$ m, $V_{M}=2$ m/s, and $\rho_{M}=7$ per/m$^2$. 
See Section~\ref{sec:mathematicl_model}. The reference time is thus $T_{M}=5\sqrt{17}/2$ s.

\begin{figure}[h!]
\centering
\begin{overpic}[width=.49\textwidth, height=0.11\textwidth, grid=false,tics=10]{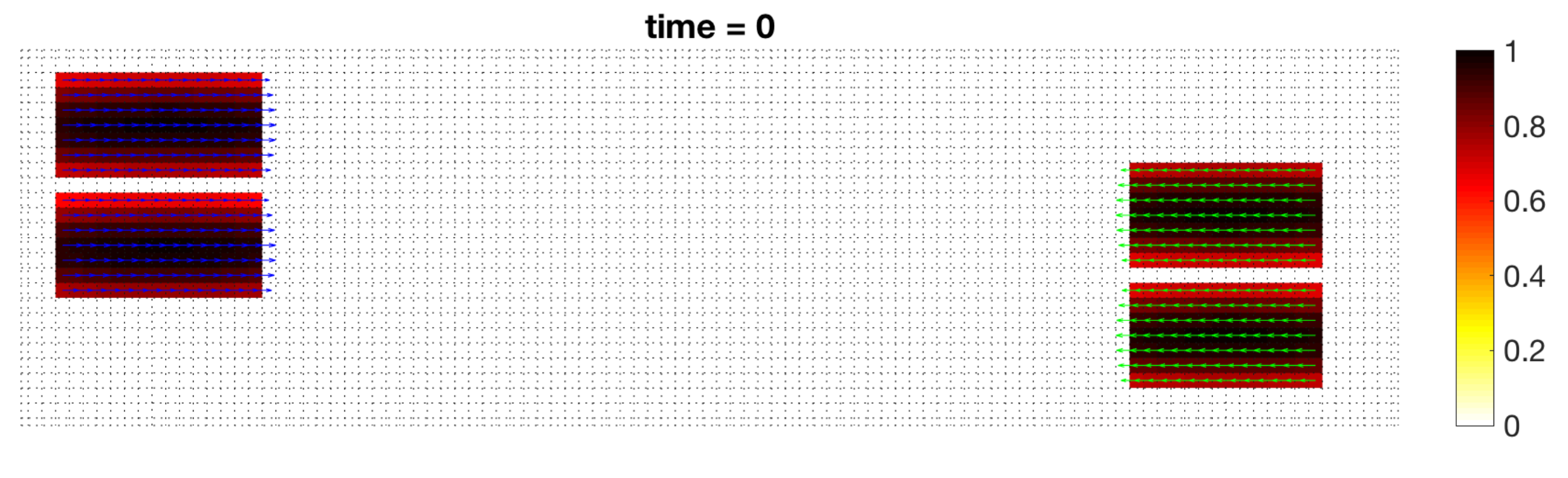}
\linethickness{1pt}
\put(1.2, 2){\line(1,0){88.3}}
\put(1.2, 20){\line(1,0){88.3}}
\end{overpic} 
\begin{overpic}[width=.49\textwidth, height=0.11\textwidth, grid=false,tics=10]{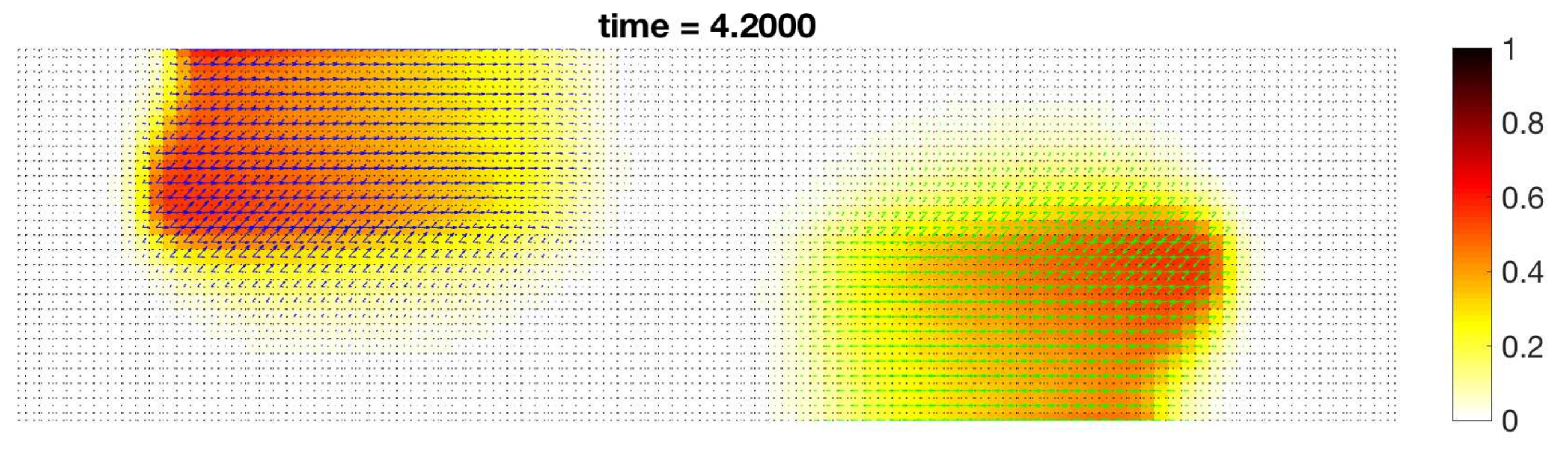}
\linethickness{1pt}
\put(1.2, 2){\line(1,0){88.3}}
\put(1.2, 20){\line(1,0){88.3}}
\end{overpic}
\begin{overpic}[width=.49\textwidth, height=0.11\textwidth,  grid=false,tics=10]{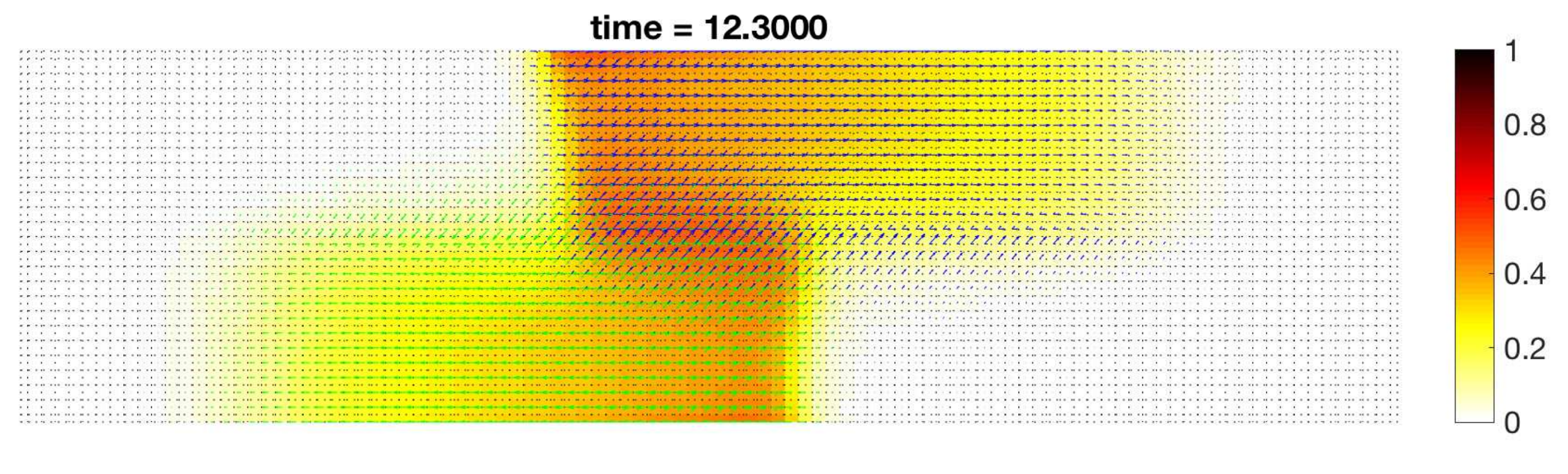}
\linethickness{1pt}
\put(1.2, 2){\line(1,0){88.3}}
\put(1.2, 20){\line(1,0){88.3}}
\end{overpic}
\begin{overpic}[width=.49\textwidth, height=0.11\textwidth,  grid=false,tics=10]{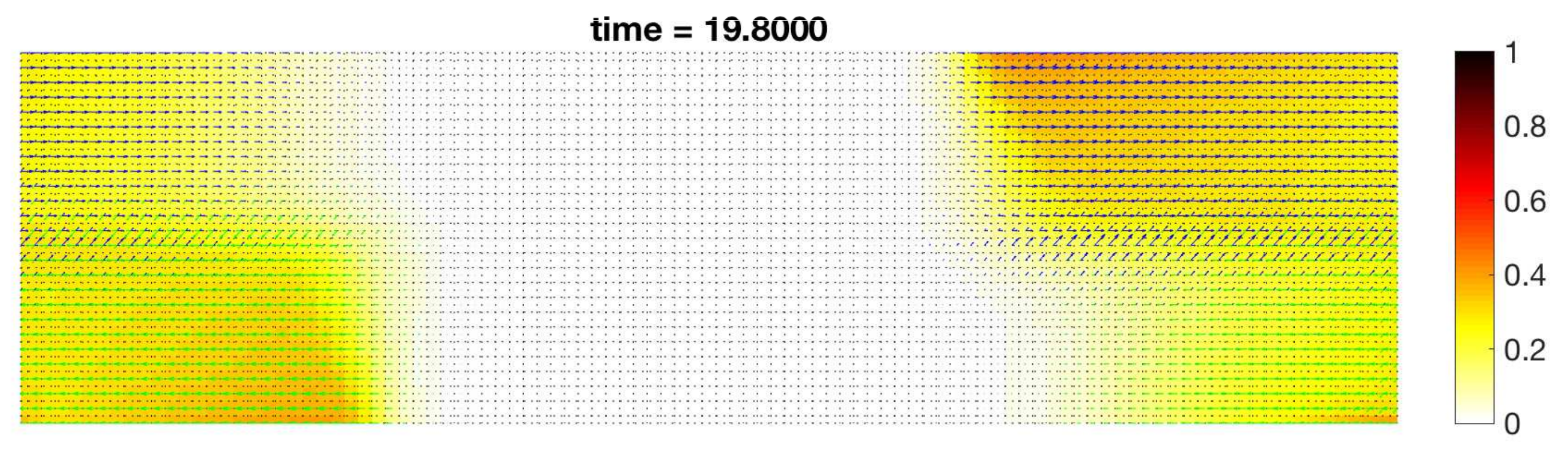}
\linethickness{1pt}
\put(1.2, 2){\line(1,0){88.3}}
\put(1.2, 20){\line(1,0){88.3}}
\end{overpic}
\begin{overpic}[width=.49\textwidth, height=0.11\textwidth, grid=false,tics=10]{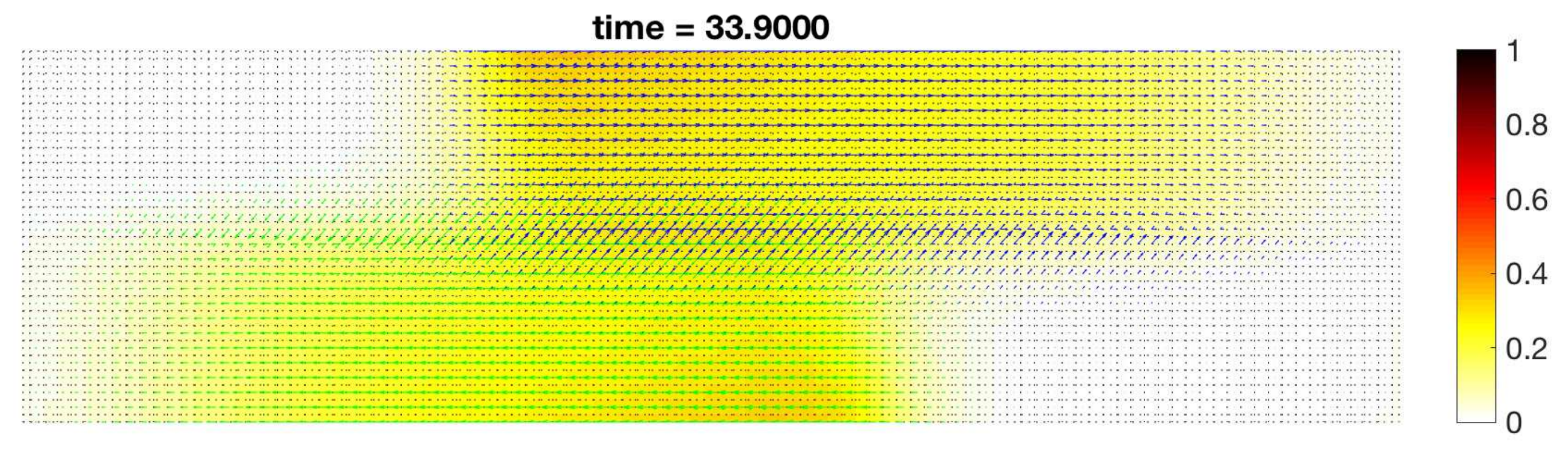}
\linethickness{1pt}
\put(1.2, 2){\line(1,0){88.3}}
\put(1.2, 20){\line(1,0){88.3}}
\end{overpic}
\begin{overpic}[width=.49\textwidth, height=0.11\textwidth, grid=false,tics=10]{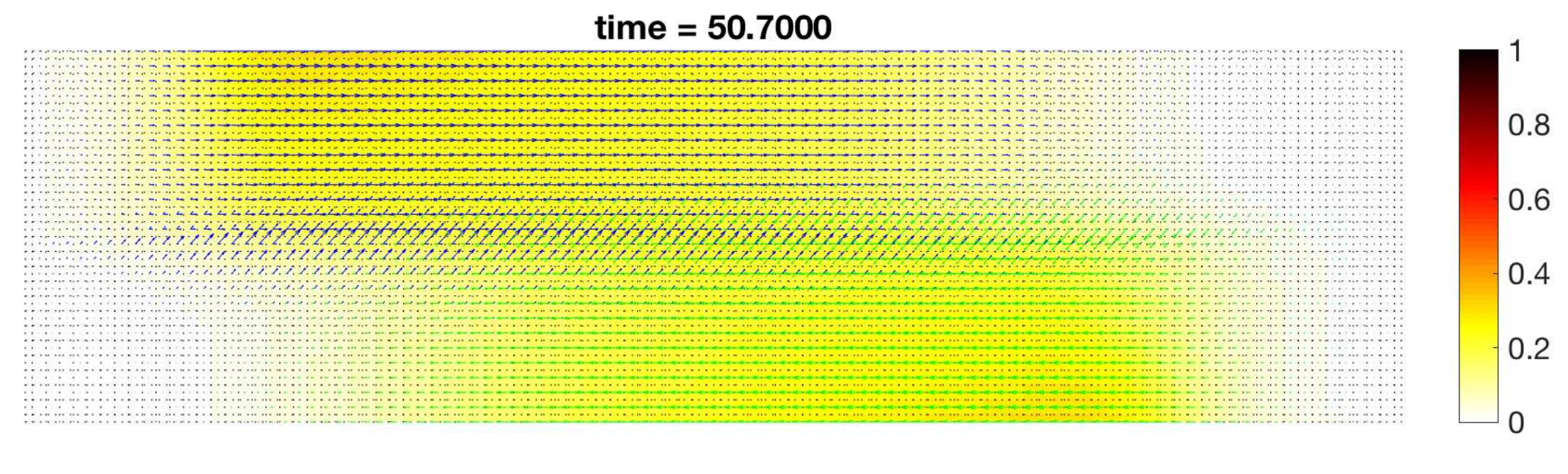}
\linethickness{1pt}
\put(1.2, 2){\line(1,0){88.3}}
\put(1.2, 20){\line(1,0){88.3}}
\end{overpic}
\begin{overpic}[width=.49\textwidth, height=0.11\textwidth,  grid=false,tics=10]{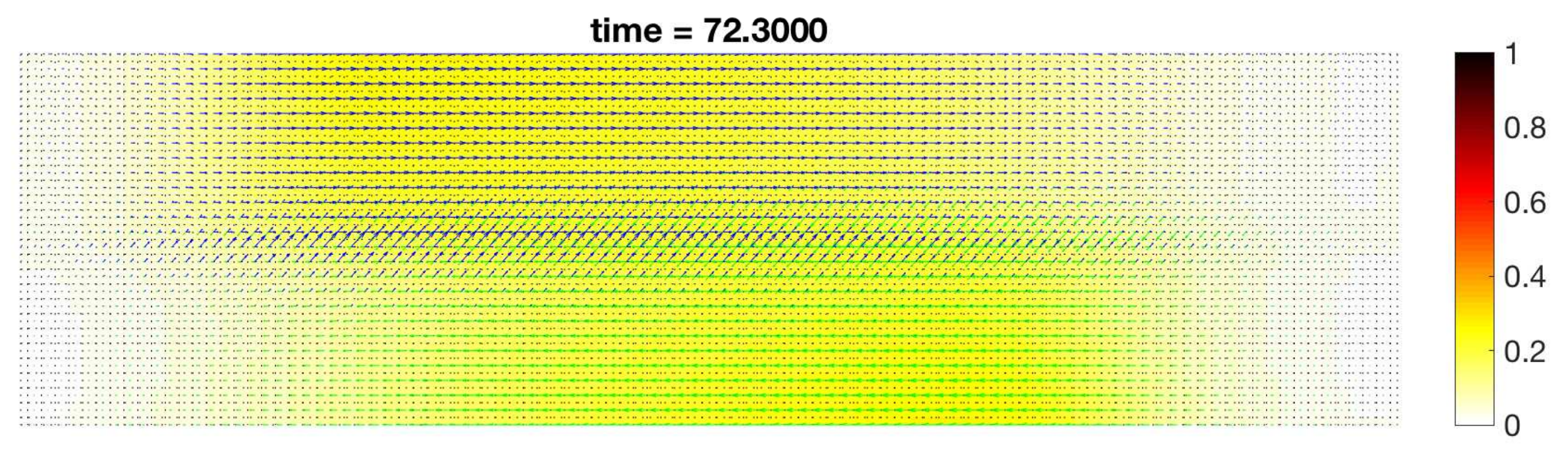}
\linethickness{1pt}
\put(1.2, 2){\line(1,0){88.3}}
\put(1.2, 20){\line(1,0){88.3}}
\end{overpic}
\caption{The movement process of 98 pedestrians grouped into four clusters with opposite initial direction $\theta_1$ and $\theta_5$ in the periodic corridor for $t = 0, 4.2, 12.3, 19.8, 33.9, 50.7, 72.3$ s, respectively.}
\label{Peri_Corr_98}
\end{figure}

\begin{figure}[h!]
\centering
\begin{overpic}[width=.49\textwidth,height=0.11\textwidth, grid=false,tics=10]{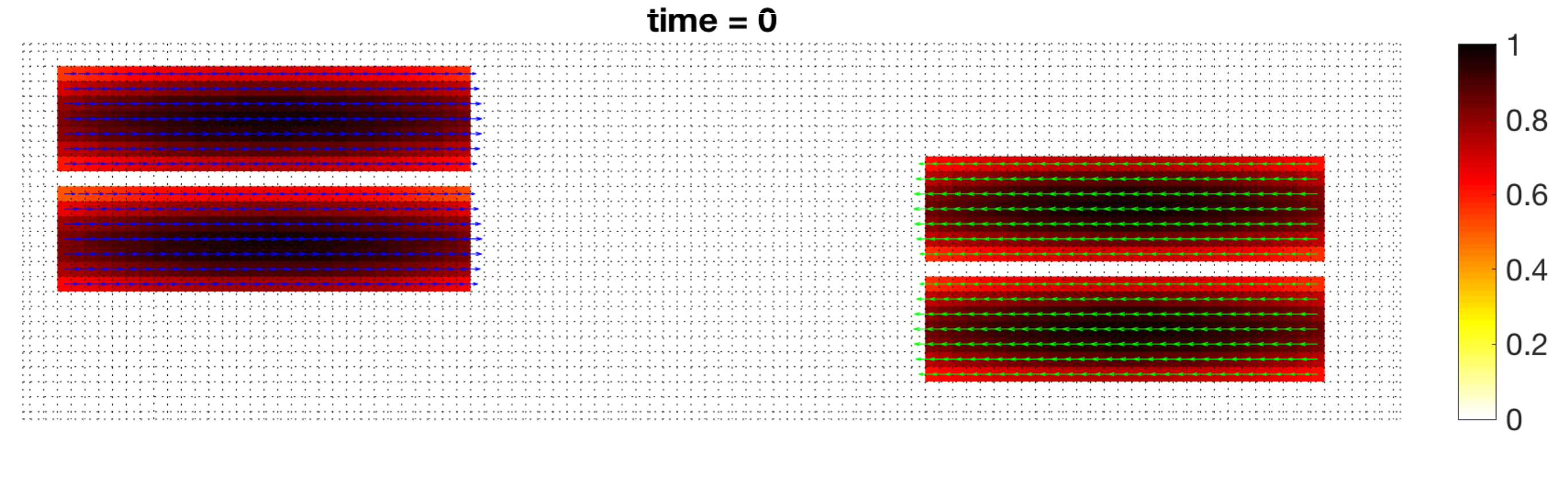}
\linethickness{1pt}
\put(1.2, 2){\line(1,0){88.3}}
\put(1.2, 20.2){\line(1,0){88.3}}
\end{overpic} 
\begin{overpic}[width=.49\textwidth, height=0.11\textwidth, grid=false,tics=10]{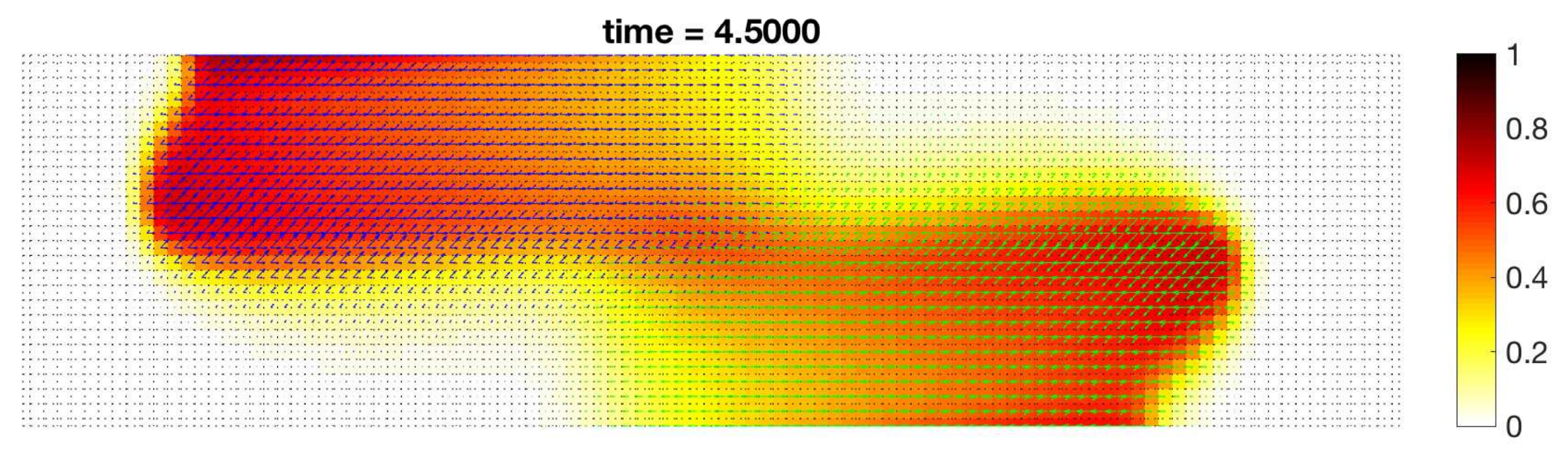}
\linethickness{1pt}
\put(1.2, 2){\line(1,0){88.3}}
\put(1.2, 20){\line(1,0){88.3}}
\end{overpic}
\begin{overpic}[width=.49\textwidth, height=0.11\textwidth,  grid=false,tics=10]{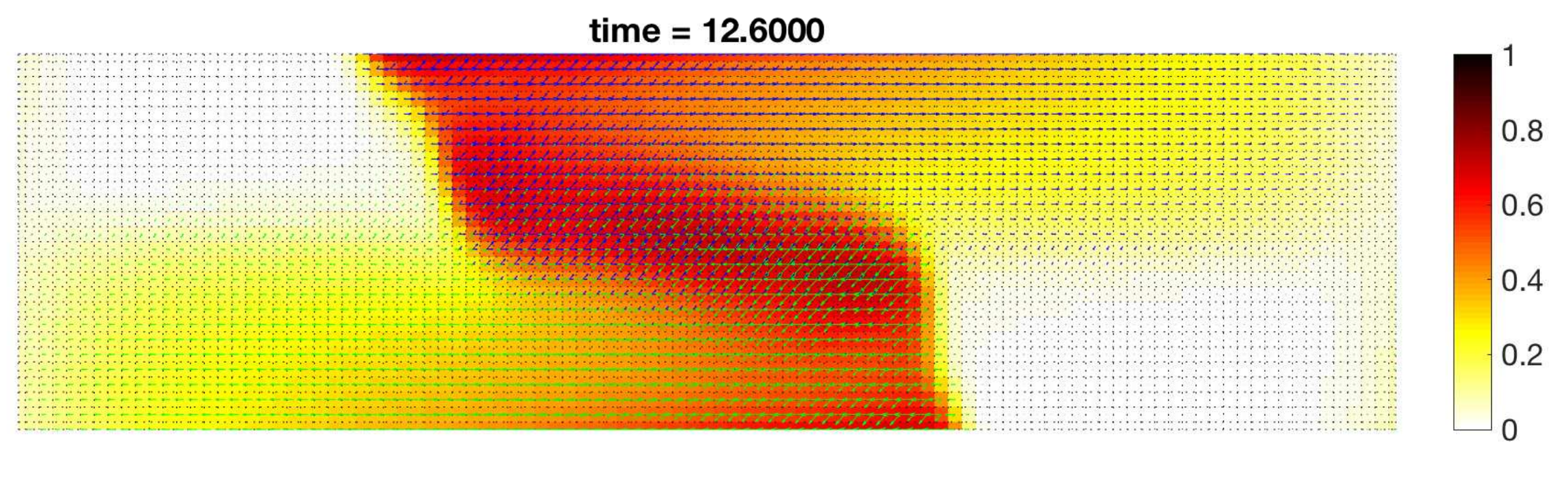}
\linethickness{1pt}
\put(1.2, 2){\line(1,0){88.3}}
\put(1.2, 20){\line(1,0){88.3}}
\end{overpic}
\begin{overpic}[width=.49\textwidth, height=0.11\textwidth,  grid=false,tics=10]{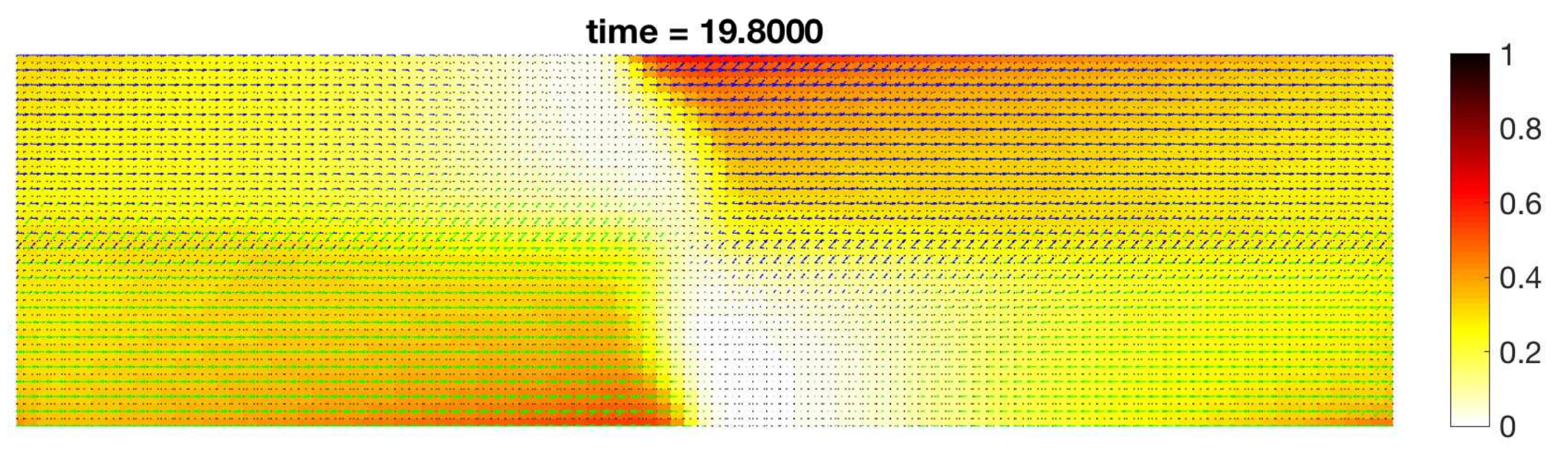}
\linethickness{1pt}
\put(1.2, 2){\line(1,0){88.3}}
\put(1.2, 20){\line(1,0){88.3}}
\end{overpic}
\begin{overpic}[width=.49\textwidth, height=0.11\textwidth, grid=false,tics=10]{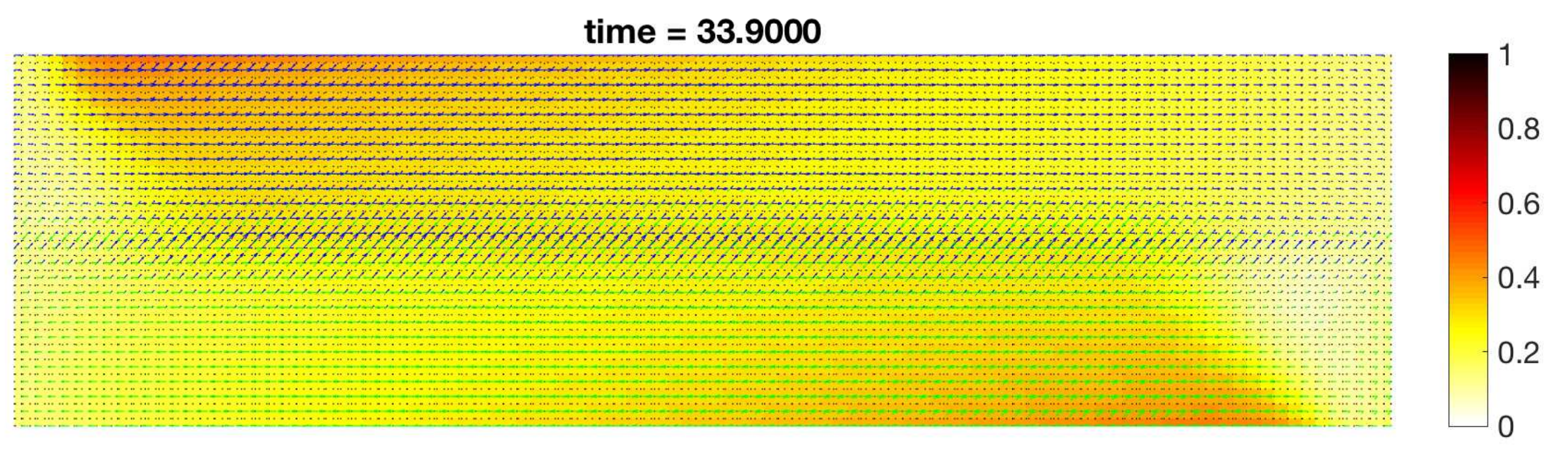}
\linethickness{1pt}
\put(1.2, 2){\line(1,0){88.3}}
\put(1.2, 20){\line(1,0){88.3}}
\end{overpic}
\begin{overpic}[width=.49\textwidth, height=0.11\textwidth, grid=false,tics=10]{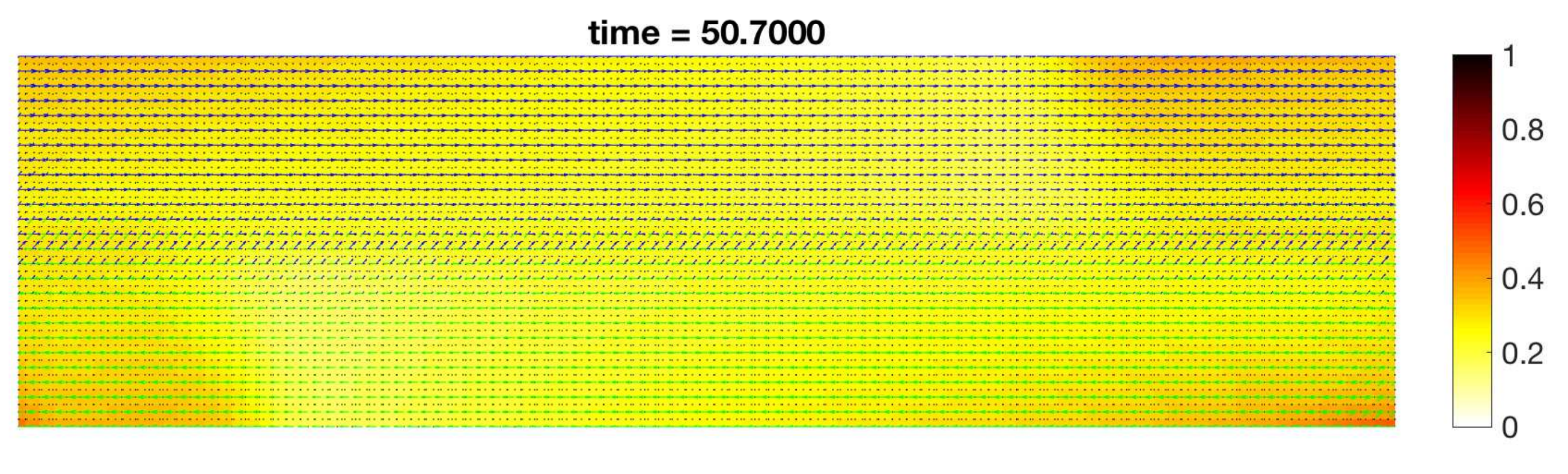}
\linethickness{1pt}
\put(1.2, 2){\line(1,0){88.3}}
\put(1.2, 20){\line(1,0){88.3}}
\end{overpic}
\begin{overpic}[width=.49\textwidth, height=0.11\textwidth,  grid=false,tics=10]{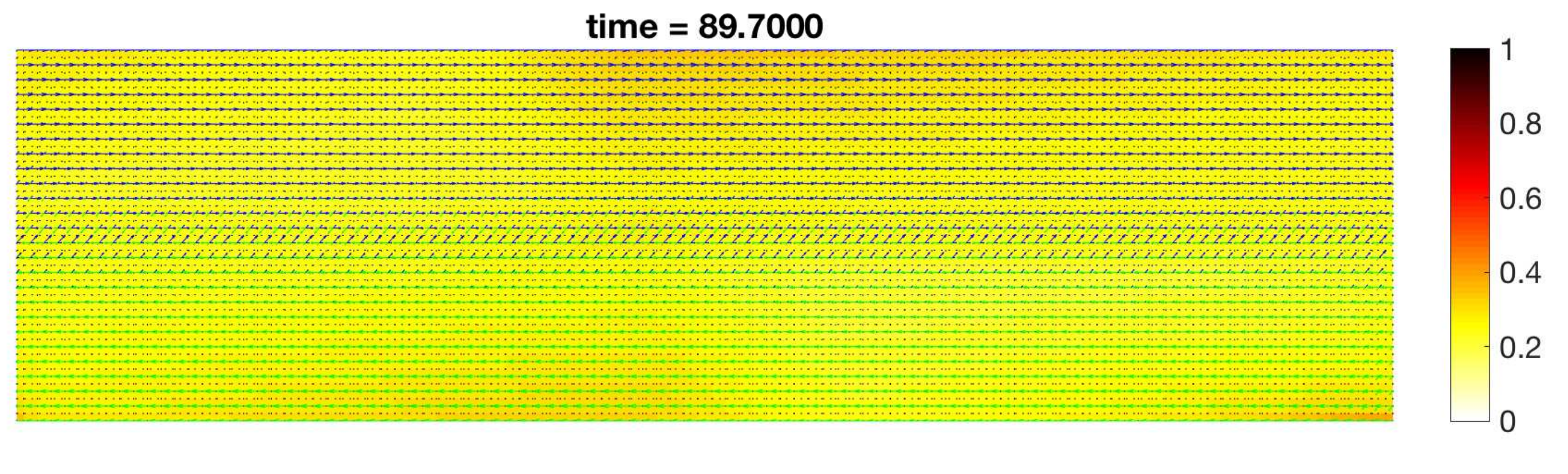}
\linethickness{1pt}
\put(1.2, 2){\line(1,0){88.3}}
\put(1.2, 20){\line(1,0){88.3}}
\end{overpic}
\caption{The movement process of 188 pedestrians grouped into four clusters with initial opposite direction $\theta_1$ and $\theta_5$ in the periodic corridor for $t = 0, 4.5, 12.6, 19.8, 33.9, 50.7, 89.7$ s, respectively.}
\label{Peri_Corr_188}
\end{figure}

We generate a mesh with $\Delta x=\Delta y$ = 0.2 m and set $\Delta t = 0.3$ s. 
We run two tests: in one 98 people are present in the computational domain, while in the
other we increase the number of people to 188. 
Pedestrians are initially distributed into four equal-area rectangular clusters 
with a parabolic density distribution with maximum density $\rho = 1$ in the center and
minimum density $\rho = 0.6$. The groups at the opposite ends of the corridor
initially move with opposite initial directions $\theta_1$ and $\theta_5$. 
See the top left panel in Figure~$\ref{Peri_Corr_98}$ and $\ref{Peri_Corr_188}$.
The rest of Figure~$\ref{Peri_Corr_98}$ and $\ref{Peri_Corr_188}$ shows the evolution
of the pedestrian dynamics. From both pictures, we see that 
pedestrians try to avoid contact by changing the direction, which leads to sorting, separation,
and lane formation.

\section{Conclusion}\label{sec:concl}

We considered a kinetic theory approach to model pedestrian dynamics in bounded domain 
and adapted it to handle obstacles. 
For the numerical approximation of the solution to our model, we applied the Lie splitting scheme
which breaks the problem into two pure advection problems and a problem involving the interaction 
with the environment and other pedestrians.

Several test cases have been considered in order to show the
ability of the model to reproduce qualitatively:
\begin{itemize}
\item[-] evacuation from a room with one exit, without and with obstacles;
\item[-] evacuation from a room with two exits and no obstacles;
\item[-] lane formation. 
\end{itemize}
In the case of the room with two exits and no obstacles, we also presented a quantitative comparison
with experimental data. Numerical results and experimental data are in 
very good agreement for medium and medium-to-large groups of people. 
With the confidence in the model given by the experimental validation, we performed
numerical tests to study evacuation for different scenarios in terms of exit sizes, 
obstacle shapes, and velocity moduli. 

\section*{Acknowledgements}
This work has been partially supported by NSF through grant DMS-1620384.

\bibliographystyle{plain}
\bibliography{crowd_dyn}

\begin{thebibliography}{10}

\bibitem{Agnelli2015}
    \newblock J. P. Agnelli, F. Colasuonno and D. Knopoff,
    \newblock  A kinetic theory approach to the dynamics of crowd evacuation from bounded domains,
     \newblock \emph{Mathematical Models and Methods in Applied Sciences}, \textbf{25} (2015), 109--129.

\bibitem{ANTONINI2006667}
    \newblock G. Antonini, M. Bierlaire and M. Weber,
    \newblock  Discrete choice models of pedestrian walking behavior,
     \newblock \emph{Transportation Research Part B: Methodological}, \textbf{40} (2006), 667--687.
     
\bibitem{ASANO2010842}
    \newblock M. Asano, T. Iryo and M. Kuwahara,
    \newblock  Microscopic pedestrian simulation model combined with a tactical model for route choice behaviour,
     \newblock \emph{Transportation Research Part C: Emerging Technologies}, \textbf{18} (2010), 842--855.     
  
  \bibitem{bandini2009}
    \newblock S. Bandini, S. Manzoni and G. Vizzari,
    \newblock Agent based modeling and simulation: An informatics perspective,
     \newblock \emph{Journal of Artificial Societies and Social Simulation}, \textbf{12} (2009).    
  
\bibitem{Bellomo2011383}
    \newblock N. Bellomo and A. Bellouquid,
    \newblock  On the modeling of crowd dynamics: Looking at the beautiful shapes of swarms,
     \newblock \emph{Networks and Heterogeneous Media}, \textbf{6} (2011), 383--399.        

\bibitem{Bellomo2013_new}
    \newblock N. Bellomo, A. Bellouquid and D. Knopoff,
    \newblock  From the microscale to collective crowd dynamics,
     \newblock \emph{Society for Industrial and Applied Mathematics Multiscale Modeling and Simulation}, \textbf{11} (2013), 943--963.          

\bibitem{Bellomo2017_book}
    \newblock N. Bellomo, A. Bellouquid, L. Gibelli and N. Outada,
    \newblock \emph{A Quest Towards a Mathematical Theory of Living Systems},
     \newblock  Modeling and Simulation in Science, Engineering and Technology, Birkhauser, 2017.   
    
\bibitem{Bellomo2011}
    \newblock N. Bellomo and C. Dogbe,
    \newblock On the modeling of traffic and crowds: A survey of models, speculations, and perspectives,
     \newblock \emph{Society for Industrial and Applied Mathematics Review}, \textbf{53} (2011), 409--463.          

\bibitem{Bellomo2015_new}
    \newblock N. Bellomo and L. Gibelli,
    \newblock Toward a mathematical theory of behavioral-social dynamics for pedestrian crowds,
     \newblock \emph{Mathematical Models and Methods in Applied Sciences}, \textbf{25} (2015), 2417--2437.   

\bibitem{Bellomo2016_new}
    \newblock N. Bellomo and L. Gibelli,
    \newblock Behavioral crowds: {Modeling and Monte Carlo} simulations toward validation,
     \newblock \emph{Computers and Fluids}, \textbf{141} (2016), 13--21.   

\bibitem{Bellomo2019_new}
    \newblock N. Bellomo, L. Gibelli and N. Outada,
    \newblock On the interplay between behavioral dynamics and social interactions in human crowds,
     \newblock \emph{Kinetic and Related Models}, \textbf{12} (2019), 397--409.                        
     
\bibitem{Bellomo2013}
    \newblock N. Bellomo, D. Knopoff and J. Soler,
    \newblock On the difficult interplay between life, ``Complexity'', and mathematical sciences,
     \newblock \emph{Mathematical Models and Methods in Applied Sciences}, \textbf{23} (2013), 1861--1913.    

\bibitem{Bellomo2012}
    \newblock N. Bellomo, B. Piccoli and A. Tosin,
    \newblock Modeling crowd dynamics from a complex system viewpoint,
     \newblock \emph{Mathematical Models and Methods in Applied Sciences}, \textbf{22} (2012).    

\bibitem{Blue1999}
    \newblock V. J. Blue and J. L. Adler,
    \newblock Cellular automata microsimulation of bidirectional pedestrian flows,
     \newblock \emph{Journal of Transportation Research Record}, \textbf{1678} (1999), 135--141.        
 
 \bibitem{Blue2000}
    \newblock V. J. Blue and J. L. Adler,
    \newblock Cellular automata microsimulation for modeling bi-directional pedestrian walkways,
     \newblock \emph{Transportation Research Part B: Methodological}, \textbf{35} (2000), 293--312.       
     
\bibitem{BURSTEDDE2001507}
    \newblock C. Burstedde, K. Klauck, A. Schadschneider, J. Zittartz,
    \newblock Simulation of pedestrian dynamics using a two-dimensional cellular automaton,
     \newblock \emph{Physica A: Statistical Mechanics and its Applications}, \textbf{295} (2001), 507--525.    

\bibitem{CHOORAMUN20121685}
    \newblock N. Chooramun, P. J. Lawrence and E. R. Galea,
    \newblock An agent based evacuation model utilising hybrid space discretisation,
     \newblock \emph{Safety Science}, \textbf{50} (2012), 1685--1694.      
   
\bibitem{Chraibi2011425}
    \newblock M. Chraibi, U. Kemloh, A. Schadschneider and A. Seyfried,
    \newblock Force-based models of pedestrian dynamics,
     \newblock \emph{Networks and Heterogeneous Media}, \textbf{6} (2011), 1556--1801.       

\bibitem{Chraibi2019}
    \newblock M. Chraibi, A. Tordeux, A. Schadschneider and A. Seyfried,
    \newblock \emph{Modelling of Pedestrian and Evacuation Dynamics},
     \newblock  Encyclopedia of Complexity and Systems Science Series, 2019, 649--669.       

\bibitem{DAI20132202}
    \newblock J. Dai, X. Li and L. Liu,
    \newblock Simulation of pedestrian counter flow through bottlenecks by using an agent-based model,
     \newblock \emph{Physica A: Statistical Mechanics and its Applications}, \textbf{392} (2013), 2202--2211.      

\bibitem{Dijkstra}
    \newblock J. Dijkstra, J. Jessurun and H. Timmermans,
    \newblock A multi-agent cellular automata model of pedestrian movement,
     \newblock \emph{Pedestrian and Evacuation Dynamics}, (2001), 173--181.  
     
 \bibitem{Einarsson2005}
     \newblock B. Einarsson,
     \newblock \emph{Accuracy and Reliability in Scientific Computing},
     \newblock  Society for Industrial and Applied Mathematics, 2005.         
 
  \bibitem{glowinski2003finite}
     \newblock R. Glowinski,
     \newblock \emph{Finite Element Methods for Incompressible Viscous Flow},
     \newblock  Handbook of numerical analysis, \textbf{9}. North-Holland, Amsterdam, 2003.     
     
\bibitem{BS:BS3830360405}
    \newblock D. Helbing,
    \newblock A mathematical model for the behavior of pedestrians,
     \newblock \emph{Behavioral Science}, \textbf{36} (1991), 298--310.        

 \bibitem{Helbing2002}
    \newblock D. Helbing, I. J. Farkas, P. Molnar, and T. Vicsek,
    \newblock Simulation of pedestrian crowds in normal and evacuation situations,
     \newblock \emph{Pedestrian and Evacuation Dynamics}, \textbf{21} (2002), 21--58.  
   
\bibitem{Helbing1995}
    \newblock D. Helbing and P. Molnar,
    \newblock Social force model for pedestrian dynamics,
     \newblock \emph{Physical Review E}, \textbf{51} (1998), 4282--4286.    
     
\bibitem{1367-2630-1-1-313}
    \newblock D. Helbing and T. Vicsek,
    \newblock Optimal self-organization,
     \newblock \emph{New Journal of Physics}, \textbf{1} (1999).         
                      
\bibitem{HUGHES2002507}
    \newblock R. L. Hughes,
    \newblock A continuum theory for the flow of pedestrians,
     \newblock \emph{Transportation Research Part B: Methodological}, \textbf{36} (2002), 507--535.    
 
\bibitem{JOHANSSON_HELBING}
    \newblock A. Johansson,  D. Helbing and P. K. Shukla,
    \newblock Specification of the social force pedestrian model by evolutionary adjustment to video tracking data,
     \newblock \emph{Advances in Complex Systems}, \textbf{10} (2007), 271--288.        

\bibitem{leveque1992numerical}
     \newblock R. J. LeVeque,
     \newblock \emph{Numerical Methods for Conservation Laws},
     \newblock 2$^{nd}$ edition,  Springer Science \& Business Media, 1992.
       
 \bibitem{Li_2012}
    \newblock X. Li, X. Yan, X. Li and J. Wang,
    \newblock Using cellular automata to investigate pedestrian conflicts with vehicles in crosswalk at signalized intersection,
     \newblock \emph{Discrete Dynamics in Nature and Society}, (2012).    
 
 \bibitem{6701214}
    \newblock S. Liu, S. Lo, J. Ma and W. Wang,
    \newblock An agent-based microscopic pedestrian flow simulation model for pedestrian traffic problems,     
    \newblock \emph{IEEE Transactions on Intelligent Transportation Systems}, \textbf{15} (2014), 992--1001.   
          
  \bibitem{Moussad2755}
     \newblock J. Moussa{\"\i}d, D. Helbing, S. Garnier, A. Johansson, M. Combe and G. Theraulaz,
     \newblock Experimental study of the behavioural mechanisms underlying self-organization in human crowds,
     \newblock \emph{Proceedings of the Royal Society B: Biological Sciences}, \textbf{276} (2009), 2755--2762.      
     
 \bibitem{Schadschneider2011}
     \newblock A. Schadschneider, W. Klingsch, H. Kluepfel, T. Kretz, C. Rogsch, A. Seyfried,
     \newblock Evacuation dynamics: Empirical results, modeling and applications,
     \newblock \emph{Extreme Environmental Events: Complexity in Forecasting and Early Warning}, (2011), 517--550.   
  
\bibitem{Schadschneider2011545}
   \newblock  A. Schadschneider and  A. Seyfried, 
    \newblock Empirical results for pedestrian dynamics and their implications for modeling,
     \newblock \emph{Networks and Heterogeneous Media}, \textbf{6} (2011), 545--560.    

 \bibitem{Seyfried2009}
     \newblock A. Seyfried, O. Passon, B. Steffen, M. Boltes, T. Rupprecht and W. Klingsch,
     \newblock New insights into pedestrian flow through bottlenecks,
     \newblock \emph{Transportation Science}, \textbf{43} (2009), 395--406.   
         
 \bibitem{5773492}
     \newblock A. Shende, M. P. Singh and P. Kachroo,
     \newblock Optimization-Based feedback control for pedestrian evacuation from an exit corridor,
     \newblock \emph{IEEE Transactions on Intelligent Transportation Systems}, \textbf{12} (2011), 1167--1176.      

  \bibitem{Steffen2009}
     \newblock B. Steffen and A.Seyfried,
     \newblock Methods for measuring pedestrian density, flow, speed and direction with minimal scatter,
     \newblock \emph{Physica A: Statistical Mechanics and its Applications}, \textbf{389} (2009), 1902--1910.     

\bibitem{TurnerPenn}
    \newblock A. Turner and A. Penn,
    \newblock Encoding Natural Movement as an Agent-Based System: An Investigation into Human Pedestrian Behaviour in the Built Environment,
     \newblock \emph{Environment and Planning B: Planning and Design}, \textbf{29} (2002), 473--490. 
 
\bibitem{Kemloh}
    \newblock A. U. K. Wagoum, A. Tordeux and W. Liao,
    \newblock Understanding human queuing behaviour at exits: an empirical study,
     \newblock \emph{Royal Society Open Science}, \textbf{4} (2017). 
       
 \bibitem{Ward150703}
     \newblock J. A. Ward, A. J. Evans and N. S. Malleson,
     \newblock Dynamic calibration of agent-based models using data assimilation,
     \newblock \emph{Royal Society Open Science}, \textbf{3} (2016).      

 \bibitem{Zhang2011}
     \newblock J. Zhang, W. Klingsch, A. Schadschneider and A. Seyfried,
     \newblock Transitions in pedestrian fundamental diagrams of straight corridors and T-junctions,
     \newblock \emph{Journal of Statistical Mechanics: Theory and Experiment}, \textbf{6} (2011).             
     
 \bibitem{6248013}
     \newblock B. Zhou, X. Wang and X. Tang,
     \newblock Understanding collective crowd behaviors: Learning a Mixture model of Dynamic pedestrian-Agents,
     \newblock \emph{2012 IEEE Conference on Computer Vision and Pattern Recognition}, (2012), 2871--2878.        
         

\end{thebibliography}

\end{document}